\input amstex
\documentstyle{amsppt}
\loadmsbm

\nologo

\TagsOnRight

\NoBlackBoxes

\define\acc{\operatorname{acc}}

\define\supp{\operatorname{supp}}
\define\Lip{\operatorname{Lip}}
\define\diam{\operatorname{diam}}
\define\dist{\operatorname{dist}}

\def\floor{\mathbin{\hbox{\vrule height1.2ex width0.8pt depth0pt
        \kern-0.8pt \vrule height0.8pt width1.2ex depth0pt}}}

\font\letter=cmss10

\font\normal=cmss10 scaled 700

\font\bignormal=cmss10 scaled 1300

\define\Haus{\text{\normal H}}

\define\scon{\text{\normal con}}

\define\sdyncon{\text{\normal dyn-con}}

\define\dyn{\text{\normal dyn}}

\define\erg{\text{\normal erg}}

\define\radius{\text{\normal rad}}

\define\distance{\text{\letter d}}

\define\LDistance{\text{\letter L}}

\define\X{\text{\bignormal X}}

\font\tenscr=callig15 scaled 800
\font\sevenscr=callig15
\font\fivescr=callig15
\skewchar\tenscr='177 \skewchar\sevenscr='177 \skewchar\fivescr='177
\newfam\scrfam \textfont\scrfam=\tenscr \scriptfont\scrfam=\sevenscr
\scriptscriptfont\scrfam=\fivescr
\def\scr#1{{\fam\scrfam#1}}

\font\tenscri=suet14
\font\sevenscri=suet14
\font\fivescri=suet14
\skewchar\tenscri='177 \skewchar\sevenscri='177 \skewchar\fivescri='177
\newfam\scrifam \textfont\scrifam=\tenscri \scriptfont\scrifam=\sevenscri
\scriptscriptfont\scrifam=\fivescri
\def\scri#1{{\fam\scrifam#1}}

\hsize = 6 true in
\vsize = 9 true in

\topmatter
\title
Dynamical multifractal zeta-functions,
multifractal pressure
and
fine
multifractal spectra
\endtitle
\endtopmatter

\centerline{\smc L\. Olsen}
\centerline{Department of Mathematics}
\centerline{University of St\. Andrews}
\centerline{St\. Andrews, Fife KY16 9SS, Scotland}
\centerline{e-mail: {\tt lo\@st-and.ac.uk}}

\topmatter
\abstract{We introduce 
multifractal pressure
and
dynamical 
multifractal zeta-functions
 providing precise information 
of  a very general class of multifractal spectra, including, for example, the 
fine
multifractal spectra of self-conformal 
 measures 
 and the fine multifractal spectra of ergodic Birkhoff averages of continuous functions.
}
\endabstract
\endtopmatter

\bigskip

\centerline{\smc Contents}
\medskip
\roster
 \item"1." Introduction.
 \item"2." The setting, Part 1: 
Self-conformal sets and self-conformal masures.
 \item"3." The setting, Part 2:
Pressure and dynamical zeta-functions.
 \item"4." Motivation of the main results.
 \item"5." Statements of the  main results.
 \item"6." Applications:
 Multifractal spectra of measures
 and
  multifractal spectra of ergodic averages.
 \item"7." 
 Proofs. Preliminary results:
 the modified multifractal pressure.
   \item"8."
   Proof of Theorem 5.3.
   \item"9."
   Proof of Theorem 5.5.   
 \item"" References
\endroster

\bigskip

\footnote""
{
\!\!\!\!\!\!\!\!
}

\leftheadtext{L\. Olsen}

\rightheadtext{Dynamical multifractal zeta-functions}

\heading{1. Introduction.}\endheading

For a Borel measure $\mu$ on $\Bbb R^{d}$ 
and
a positive number $\alpha$,
let us consider
the set  of
those points
$x$ in $\Bbb R^{d}$ for which the measure
$\mu(B(x,r))$ of the ball
$B(x,r)$ with center $x$ and radius $r$ behaves like
$r^{\alpha}$ for small $r$,
i\.e\. the set
 $$
   \Bigg\{
   x\in K
   \,\Bigg|\,
   \lim_{r\searrow 0}
   \frac{\log\mu(B(x,r))}{\log r}
   =
   \alpha
   \Bigg\}\,.
   \tag1.1
 $$
If the intensity of the measure $\mu$ varies very widely, it may
happen that the sets in (1.1)
display a
fractal-like character for a range of values of $\alpha$. In this case
it is natural to study
the Hausdorff dimensions of the sets in (1.1)
 as $\alpha$
varies.
We therefore define the  fine multifractal spectrum of $\mu$
by
 $$
 f_{\mu}(\alpha)
 =
 \dim_{\Haus}
   \Bigg\{
   x\in K
   \,\Bigg|\,
   \lim_{r\searrow 0}
   \frac{\log\mu(B(x,r))}{\log r}
   =
   \alpha
   \Bigg\}\,.
 \tag1.2
$$
where $\dim_{\Haus}$ denotes the Hausdorff dimension;
here and below
we use the following convention,
namely, we define the Hausdorff
of the empty set to be $-\infty$, i\.e\. we put
 $$
 \dim_{\Haus}\varnothing
 =
 -
 \infty\,.
 $$
The second main ingredient in multifrctal analysis is the  
Renyi dimensions.
Renyi
dimensions quantify the varying intensity of a measure by analyzing its moments at different scales. 
Formally, 
for $q\in\Bbb R$,
the $q$'th
Renyi dimensions 
$\tau_{\mu}(q)$
of $\mu$ is defined by
 $$
 \tau_{\mu}(q)
 =
 \lim_{r\searrow 0}
 \frac{\dsize \log\int\limits_{\,\,\,K}\mu(B(x,r))^{q-1}\,d\mu(x)}{-\log r}\,,
 $$
provided the limit exists. 
One of the main problems
in multifractal analysis is to 
 understand the
 multifractal spectrum and the Renyi dimensions,
 and
their relationship with each other.
During the past 20 years
there has been an enormous interest 
in
computing the multifractal spectra of measures
in 
the mathematical literature
and within the last
15 years the multifractal spectra of various classes of measures
in Euclidean space $\Bbb R^{d}$ 
exhibiting some degree of self-similarity have been computed 
rigorously, 
see
the textbooks [Fa2,Pe]
and the references therein.

Dynamical zeta-functions were introduced by Artin \& Mazur
in the mid 1960's [ArMa]
based on an analogy with the number theoretical zeta-functions
associated with a function field
over a finite ring.
Subsequently Ruelle [Rue1,Rue2]
associated
zeta-functions
to 
certain statistical mechanical models in one dimensions.
Motivated by the 
powerful techniques
provided by 
the use of Artin-Mazur zeta-functions in dynamical systems
and Ruelle
zeta-functions in dynamical systems,
Lapidus and collaborators
(see the intriguing books by Lapidus \& van Frankenhuysen [Lap-vF1,Lap-VF2]
and the references therein)
have recently
introduced and pioneered
to use of
zeta-functions in fractal geometry.
Inspired by this,
within the past 2-3 years
several authors have 
paralleled this development
by
introducing
zeta-functions
into  multifractal geometry.
For example,
 in
[Bak,MiOl]
the authors
introduced
multifractal zeta-functions 
tailored to study
multifractal spectra
of self-conformal measures,
and 
in
[Le-VeMe,Ol4,Ol5]
the authors 
introduced 
multifractal zeta-functions 
designed
 to study
the 
multifractal Renyi dimensions
of self-conformal measures.
In addition, we note that Lapidus and collaborators have introduced 
various intriguing multifractal zeta-functions [LapRo,LapLe-VeRo]. 
However, the multifractal 
zeta-functions 
in [LapRo,LapLe-VeRo] serve
 very different purposes and are significantly different from 
the multifractal zeta-functions introduced in 
this paper and in
[Bak,Le-VeMe,MiOl,Ol4,Ol5].

It has been a major challenge 
to 
introduce 
and develop 
a
natural and meaningful
theory of
multifractal zeta-functions
paralleling
the existing
 powerful 
theory
 of 
dynamical zeta-functions introduced and developed by Ruelle [Rue1,Rue2]
and others, see for example,  
the surveys and books 
[Bal1,Bal2,ParPo1,ParPo2]
and the references therein.
The purpose of this paper 
is to propose such a theory.
In particular,
we introduce
a family
of multifractal zeta-functions
motivated
by
the definition
of Ruelle's
dynamical zeta-functions.
Whereas
the zeta-functions in
[Bak,Ol4,Ol5]
were designed to study the
multifractal Renyi
dimensions,
the zeta-functions in this paper 
(and the zeta-functions in [MiOl])
are tailored to
 provide
precise information of very general classes of multifractal spectra, 
including, for example, the multifractal 
spectra of self-conformal measures and the multifractal spectra of ergodic 
Birkhoff averages 
of continuous functions, see Section 6.

The framework developed in this paper will be 
formulated 
in the setting of 
self-conformal sets
and
self-conformal measures.
For this reason we will now briefly recall the definition 
of self-conformal constructions.

\bigskip

%
%

\heading{2. The setting, Part 1:
Self-conformal sets and self-conformal measures.}\endheading

{\bf 2.1. Notation from symbolic dynamics.}
We first recall the 
 notation
and terminology  from symbolic dynamics 
that will be used in this paper
Fix a positive integer $N$.
 Let
$\Sigma
 =
 \{1,\ldots,N\}$ and 
 for a positive integer $n$,
 write
 $$
 \aligned
 \Sigma^{n}
&=
\{1,\dots,N\}^{n}\,,\\
 \Sigma^{*}
&=
\bigcup_{m}\Sigma^{m}\,,\\
 \Sigma^{\Bbb N}
&=
\{1,\dots,N\}^{\Bbb N}\,,
 \endaligned
 $$
i\.e\. $\Sigma^{n}$ is the family of all
strings
$\bold i=i_{1}\ldots i_{n}$
of length $n$ 
with $i_{j}\in\{1,\ldots,N\}$;
$\Sigma^{*}$ is the family of all finite strings
$\bold i=i_{1}\ldots i_{m}$
with $m\in\Bbb N$ and $i_{j}\in\{1,\ldots,N\}$;
and 
$\Sigma^{\Bbb N}$ is the family of all
infinite
strings
$\bold i=i_{1}i_{2}\ldots $
with $i_{j}\in\{1,\ldots,N\}$.
For an infinite string 
$\bold i=i_{1}i_{2}\ldots\in\Sigma^{\Bbb N} $
and a positive integer $n$, we will write
$\bold i|n
=
i_{1}\ldots i_{n}$.
In addition, for
a positive integer $n$
and
a finite string 
$\bold i=i_{1}\ldots i_{n}\in\Sigma^{n} $
with length equal to  $n$,
 we will write
$|\bold i|
=
n$, and we let $[\bold i]$ denote the cylinder 
generated by $\bold i$, i\.e\.
  $$
  [\bold i]
 =
 \Big\{
 \bold j\in\Sigma^{\Bbb N}
 \,\Big|\,
 \bold j|n=\bold i
 \Big\}\,.
 $$
Also, let $S:\Sigma^{\Bbb N}\to\Sigma^{\Bbb N}$ denote the shift map, i\.e\.
 $$
 S(i_{1}i_{2}\ldots)
 =
 i_{2}i_{3}\ldots\,.
 $$

%

%

\bigskip

{\bf 2.2. Self-conformal sets and self-conformal measures.}
Next, we recall the definition of
 self-conformal 
(and self-similar) sets and measures. 
A conformal iterated function system with probabilities 
is a list
$\big(
 \,
 V
 \,,\,
 X
 \,,\,
 (S_{i})_{i=1,\ldots,N}
 \,
 \big)$
where
\roster
\item"$\bullet$"
$V$ is an open, connected subset of $\Bbb R^{d}$.
\item"$\bullet$" 
$X$ is a compact set with $X\subseteq V$
and
$X^{\circ\,-}=X$.
\item"$\bullet$" 
$S_{i}:V\to V$ is a contractive
$C^{1+\gamma}$ diffeomorphism with
$0<\gamma<1$ such that
 $S_{i}X
 \subseteq
 X$
for all $i$. 
\item"$\bullet$"
The Conformality Condition:
For each $x\in V$, we have that
$(DS_{i})(x)$ is a contractive similarity map, i\.e\.
there exists
$r_{i}(x)\in(0,1)$ such that
$|(DS_{i})(x)u-(DS_{i})(x)v|
 =
 r_{i}(x)|u-v|$
for all $u,v\in\Bbb R^{d}$;
here $(DS_{i})(x)$ denotes the derivative of $S_{i}$ at $x$.
\endroster

\noindent
It follows from [Hu] that there exists a unique
non-empty compact set $K$ with $K\subseteq X$
 such that
 $$
 K
 =
 \bigcup_{i}\,
 S_{i}K\,.
 \tag2.1
 $$
The set $K$
is called the self-conformal set
associated with the list 
$\big(
 \,
 V
 \,,\,
 X
 \,,\,
 (S_{i})_{i=1,\ldots,N}
 \,
 \big)$;
in particular,
 if each map $S_{i}$ is a contracting similarity, then the 
 set $K$ is called the 
 self-similar set
 associated with the list 
$\big(
 \,
 V
 \,,\,
 X
 \,,\,
 (S_{i})_{i=1,\ldots,N}
 \,
 \big)$.
In addition,
if
$(p_{i})_{i=1,\ldots,N}$ is a probability vector
then it follows from [Hu]
that there is a unique
probability measure $\mu$ with
$\supp\mu=K$ such that
 $$
 \mu
 =
 \sum_{i}\,
 p_{i}\,\mu\circ S_{i}^{-1}\,.
 \tag2.2
 $$
The measure $\mu$
is called the self-conformal measure
associated with the list
$\big(
 \,
 V
 \,,\,
 X
 \,,\,
 (S_{i})_{i=1,\ldots,N}
  \,,\,\allowmathbreak
 (p_{i})_{i=1,\ldots,N}
 \,
 \big)$;
 if each map $S_{i}$ is a contracting similarity,
 then the
 measure $\mu$
is called the self-similar measure
associated with the list
$\big(
 \,
 V
 \,,\,
 X
 \,,\,
 (S_{i})_{i=1,\ldots,N}
  \,,\,\allowmathbreak
 (p_{i})_{i=1,\ldots,N}
 \,
 \big)$.
We will frequently assume that the list
$\big(
 \,
 V
 \,,\,
 X
 \,,\,
 (S_{i})_{i=1,\ldots,N}\allowmathbreak
 \,
 \big)$
 satisfies
  the Open Set Condition defined below.
Namely,
the list 
$\big(
 \,
 V
 \,,\,
 X
 \,,\,
 (S_{i})_{i=1,\ldots,N}
 \,
 \big)$
satisfies the Open Set Condition (OSC) if
there exists an
open, non-empty and bounded set $O$
with $O\subseteq X$
and
$S_{i}O
 \subseteq
 O$
for all $i$ such that
$S_{i}O
 \,\cap\,
 S_{j}O
 =
 \varnothing$
for all $i,j$ with 
$i\not=j$.

For
$\bold i=i_{1}\ldots i_{n}\in\Sigma^{*}$, we will write
 $$
\aligned
p_{\bold i}
&=
p_{i_{1}}\cdots p_{i_{n}}\,,\\
S_{\bold i}
&=
S_{i_{1}}\cdots S_{i_{n}}\,,\\
K_{\bold i}
&=
S_{\bold i}K\,.
\endaligned
\tag2.3
$$

Next,
 we 
define 
the natural projection map $\pi:\Sigma^{\Bbb N}\to K$
 by
 $$
 \Big\{\,\pi(\bold i)\,\Big\}
 =
 \bigcap_{n}K_{\bold i|n}
 \tag2.4
 $$
for $\bold i\in\Sigma^{\Bbb N}$.
Finally, we define the scaling map
$\Lambda:\Sigma^{\Bbb N}\to\Bbb R$
by
 $$
 \Lambda(\bold i)
 =
 \log
 |DS_{i_{1}}(\pi S\bold i)|
 \tag2.5
 $$
for $\bold i=i_{1}i_{2}\ldots\in\Sigma^{\Bbb N}$.

\bigskip

\heading{3. The setting, Part 2: 
Pressure and dynamical zeta-functions.}\endheading

Throughout 
 this section, 
and in the remaining parts of the paper, we will used the
following notation.
 Namely,
if $(a_{n})_{n}$ is a sequence of 
complex numbers
and if $f$ is the power
series
defined by $f(z)=\sum_{n}a_{n}z^{n}$ for $z\in\Bbb C$,
then we will denote the 
radius
of convergence of $f$ by $\sigma_{\radius}(f)$, i\.e\. we 
write
 $$
 \sigma_{\radius}(f)
 =
 \text{
 \lq\lq the radius of convergence of $f$"
 }\,.
 $$

Our
definitions
and
 results are
 motivated by the 
 notion of pressure
 from the thermodynamic formalism
  and
 the
 dynamical zeta-functions introduced by
 Ruelle [Rue1,Rue2]; see, also 
 [Bal1,\allowlinebreak
 Bal2,\allowlinebreak
 ParPo1,\allowlinebreak
 ParPo2].
 In addition,
Bowen's formula 
expressing the 
Hausdorff dimension of 
a self-conformal set
in terms of the pressure 
(or the dynamical zeta-function)
of the scaling map $\Lambda$ in (2.5)
also plays a leitmotif 
in our work.
Because of this we
now
recall the definition 
of pressure and dynamical zeta-function,
and the statement of Bowen's formula.
Let
$\varphi:\Sigma^{\Bbb N}\to\Bbb R$
be a
continuous function.
 The
pressure
of 
 $\varphi$
 is defined 
by
 $$
 \align
 P(\varphi)
 &=
\lim_{n}
  \,\,
  \frac{1}{n}
 \,\,
 \log
   \sum
  \Sb
  |\bold i|=n
  \endSb
  \,\,
  \sup_{\bold u\in[\bold i]}
  \,\,
   \exp
 \,\,
 \sum_{k=0}^{n-1}\varphi S^{k}\bold u\,,
 \tag3.1
  \endalign
 $$
 see
 [Bo2]
or
[ParPo2];
we note that
it is well-known that the limit in (3.1) exists.
Also,
the dynamical zeta-function of $\varphi$ is defined by
$$
  \zeta^{\dyn}(\varphi;z)
  =
   \sum_{n}
   \,\,
   \frac{z^{n}}{n}
 \left(
   \sum
  \Sb
 |\bold i|=n
  \endSb
  \,\,
  \sup_{\bold u\in[\bold i]}
  \,\,
 \exp
 \,\,
 \sum_{k=0}^{n-1}\varphi S^{k}\bold u
\right)
\tag3.2
 $$
for those complex numbers $z$ for which the series converge, see
[ParPo2].
We now  list two easily established and well-known
properties of 
the
pressure
$P(\varphi)$
and of 
the radius of convergence
$\sigma_{\radius}\big(\,\zeta^{\dyn}(\varphi;\cdot)\,\big)$
of
the power-series  
$\zeta^{\dyn}(\varphi;\cdot)$.
While both results are well-known and easily proved
(see, for example, [Bar,Fa2]), 
we have decided to list them 
since 
they play an important part 
in 
the discussion of our results.

\bigskip

\proclaim{Theorem A (see, for example, [Bar,Fa2]). Radius of convergence}
Fix a continuous function
$\varphi:\Sigma^{\Bbb N}\to\Bbb R$.
Then we have
  $$
  -
  \log\sigma_{\radius}\big(\,\zeta^{\dyn}(\varphi;\cdot)\,\big)
  =
  P(\varphi)\,.
  $$
\endproclaim

\bigskip

\proclaim{Theorem B (see, for example, [Bar,Fa2]). 
Continuity and monotony properties of the pressure}
Fix a  a
continuous function
$\Phi:\Sigma^{\Bbb N}\to\Bbb R$
with $\Phi<0$.
Then the function
 $
 t\to P(t\Phi),
 $
where $t\in\Bbb R$, 
is continuous, strictly decreasing and convex 
with
$\lim_{t\to-\infty}P(t\Phi)=\infty$
and
 $\lim_{t\to\infty}P(t\Phi)=-\infty$.
 In particular, there is a unique real number 
 $s$ such that
  $$
  P(s\Phi)=0\,;
  $$
 alternatively,
 $s$ 
  is the unique real number such that
  $$
  \sigma_{\radius}\big(\,\zeta^{\dyn}(s\Phi;\cdot)\,\big)
  =
  1\,.
  $$
 \endproclaim

\bigskip

\noindent
The main importance of the 
pressure
(for the purpose of this exposition)  is that 
it provides a beautiful formula for the 
Hausdorff dimension of
a self-conformal set satisfying the OSC.
This
result
was first noted by 
[Bo1]
(in the setting of quasi-circles)
and is the content of the next result.

\bigskip

\proclaim{Theorem C (see, for example, [Bar,Fa2]). Bowen's formula}
Let
 $K$ be the 
 self-conformal set defined by (2.1)
 and
 let $\Lambda:\Sigma^{\Bbb N}\to\Bbb R$
 be the scaling function defined by (2.5).
 Let
$s$ be the unique real number
such that
  $$
  P(s\Lambda)=0\,;
  $$
 alternatively,
 $s$ 
  is the unique real number such that
  $$
  \sigma_{\radius}\big(\,\zeta^{\dyn}(s\Lambda;\cdot)\,\big)
  =
  1\,.
  $$  
If the OSC is satisfied, then we have
  $$
  \dim_{\Haus}K=s\,.
  $$ 
 \endproclaim

\bigskip

\noindent
Any meaningful
theory of dynamical multifractal zeta-functions
is likely to
produce
multifractal analogues of Bowen's equation.
We will propose a framework for such a theory in Section 5.
However, before doing so,
we believe that it is useful 
to
illustrate 
the underlying ideas
 in a  simple setting.
 For this reason we will now 
 illustrate 
 how  meaningful
 multifractal dynamical zeta-functions might be defined
for self-conformal measures.

\bigskip

\heading{4. Motivation of  the main results.}\endheading

To illustrate the ideas behind our main 
definitions
in a simple setting, 
we consider the following example
involving self-conformal
measures.
Fix a a conformal iterated function system 
 $(
 \,
 V
 \,,\,
 X
 \,,\,
 (S_{i})_{i=1,\ldots,N}
 \,
 )$
 and a 
 be a probability vector $(p_{1},\ldots,p_{N})$.
We let $K$ denote the self-conformal set associated with the list
 $(
 \,
 V
 \,,\,
 X
 \,,\,
 (S_{i})_{i=1,\ldots,N}
 \,
 )$,
 i\.e\. $K$ is the unique 
 non-empty and compact subset  of $\Bbb R^{d}$
satisfying (2.1),
and we let
 $\mu$ be the
self-conformal measure associated with the list
$(S_{1},\ldots,\allowmathbreak S_{N},\allowmathbreak p_{1},\allowmathbreak \ldots,p_{N})$,
i\.e\.
$\mu$ is
the unique Borel probability measure on $\Bbb R^{d}$ 
satisfying (2.2).
Recall, that the multifractal spectrum of $\mu$ is defined by
 $$
 \align
 f_{\mu}(\alpha)
&=
   \,\dim_{\Haus}
   \left\{x\in K
    \,\left|\,
     \lim_{r\searrow0}
     \frac
     {\log\mu B(x,r)}{\log r}
     =
     \alpha
      \right.
   \right\}\,,
  \endalign
 $$
for $\alpha\in\Bbb R$.
 If the OSC is satisfied,
then
the multifractal spectrum
$f_{\mu}$
is given by the following formula.
Namely,
define $\Phi:\Sigma^{\Bbb N}\to\Bbb R$ by $\Phi(\bold i) =\log p_{i_{1}}$
for
$\bold i=i_{1}i_{2}\ldots\in\Sigma^{\Bbb N}$
and let $\Lambda:\Sigma^{\Bbb N}\to\Bbb R$ denote the 
scaling map in (2.5).
Next, define $\beta:\Bbb R\to\Bbb R$
by
 $$
P\big(
\,
q\Phi+\beta(q)\Lambda
\,
\big)
=
0;
\tag4.1
 $$
 alternatively,
 the function
 $\beta:\Bbb R\to\Bbb R$
is defined by
 $$
 \sigma_{\radius}
 \big(
 \,
 \zeta^{\dyn}(q\Phi+\beta(q)\Lambda;\cdot)
 \,
 \big)
 =
 1\,.
 \tag4.2
 $$
If the OSC is satisfied, then it follows from [CaMa,Pa] that
 $$
 f_{\mu}(\alpha)
 =
 \beta^{*}(\alpha)
 \tag4.3
 $$
for all $\alpha\in\Bbb R$
 where $\beta^{*}$ denotes the Legendre transform of $\beta$;
 recall, that if
 if $\varphi:\Bbb R\to\Bbb R$ is a function,
 then the Legendre transform 
 $\varphi^{*}:\Bbb R\to[\-\infty,\infty]$
 of $\varphi$ is defined
 by
$ \varphi^{*}(x)
  =
  \inf_{y}
  (xy+\varphi(y))$.

 While one may argue that
 (4.1) and (4.3) provide a 
 pressure formula
  for the 
 multifractal spectrum
 $ f_{\mu}(\alpha)$
 of a self-sconformal measure
 (or, alternatively, that
 (4.2) and (4.3)
 provide a zeta-function formula
 for the 
 multifractal spectrum
 $ f_{\mu}(\alpha)$
 of a self-sconformal measure), this formula 
 can hardly be said
 to be 
 in the 
 spirit of Bowen's formula.
Adopting this viewpoint,
for a given $\alpha\in\Bbb R$,
it is natural
to attempt to introduce
dynamical 
multifractal zeta-functions 
$\zeta_{\alpha}^{\dyn\text{-}\scon}$
of self-conformal measures
tailored, for example, to
see the  multifractal  
decomposition sets
 $$
   \Bigg\{x\in K
    \,\Bigg|\,
     \lim_{r\searrow0}
     \frac
     {\log\mu B(x,r)}{\log r}
     =
     \alpha
   \Bigg\}
   $$
more directly
and, as a result of this, hopefully lead to a 
better conceptual understanding of the pressure formula
(4.3).
More precisely,
and 
motivated by Bowen's formula,
for each $\alpha\in\Bbb R$
it seems 
natural 
to expect that any dynamically meaningful
multifractal zeta-function
$\zeta_{\alpha}^{\dyn\text{-}\scon}$
should have the following property:
there
 is a unique
real number
$\,\scr f\,\,(\alpha)$
such that
$$
\sigma_{\radius}\big(
\,
\zeta_{\alpha}^{\dyn\text{-}\scon}(\,\scr f\,\,(\alpha)\,\Lambda;\cdot)
\,
\big)
=
1\,,
$$
and the  number 
 $\,\scr f\,\,(\alpha)$
equals the multifractal spectrum $f_{\mu}(\alpha)$, i\.e\.
$$
f_{\mu}(\alpha)
=
\,\scr f\,\,(\alpha)\,.
$$
Since 
$f_{\mu}(\alpha)$
measures 
the size of the 
set of points 
$x$ for which
$\lim_{r\searrow 0}\frac{\log\mu(B(x,r))}{\log r}=\alpha$
and since
$\frac{\log\mu(B(x,r))}{\log r}$
has the same form as
$\frac{\log p_{\bold i}}{\log \diam K_{\bold i}}$,
it is natural to
define the 
dynamical 
self-conformal multifractal zeta-function $\zeta_{\alpha}^{\dyn\text{-}\scon}(\varphi;\cdot)$ 
of a continuous function $\varphi:\Sigma^{\Bbb N}\to\Bbb R$
by
 $$
 \zeta_{\alpha}^{\dyn\text{-}\scon}(\varphi;z)
  =
   \sum_{n}
   \,\,
   \frac{z^{n}}{n}
 \left(
   \sum
 \Sb
 |\bold i|=n\\
 {}\\
 \frac{\log p_{\bold i}}{\log \diam K_{\bold i}}=\alpha
 \endSb  \,\,
  \sup_{\bold u\in[\bold i]}
  \,\,
 \exp
 \,\,
 \sum_{k=0}^{n-1}\varphi S^{k}\bold u
\right)
\tag4.4
 $$
for those complex numbers $z$ for which the series converges.
The main difference between the classical dynamical zeta-function (3.2) and
its proposed multifractal counter part (4.4) 
is that in (4.4) we only sum over those strings $\bold i$ 
with $|\bold i |=n$
that are multifractally
relevant.
An easy and straight forward calculation,
which we present in Observation 4.1  below,
shows that
if 
there is a unique
real number
$\,\scr f\,\,(\alpha)$
such that
$$
\sigma_{\radius}
\big(
\,
\zeta_{\alpha}^{\dyn\text{-}\scon}(\,\scr f\,\,(\alpha)\,\Lambda;\cdot)
\,
\big)
=
1\,,
$$
then
 this number is less than $f_{\mu}(\alpha)$,
i\.e\.
$$
\,\scr f\,\,(\alpha)
\le
f_{\mu}(\alpha)\,.
\tag4.5
$$

\bigskip

\proclaim{Observation 4.1}
Let $\mu$ 
be the self-conformal measure defined by (2.2)
and let $\Lambda:\Sigma^{\Bbb N}\to\Bbb R$
be the scaling function defined 
by (2.5).
For $\alpha,t\in\Bbb R$,
let
$\zeta_{\alpha}^{\dyn\text{-}\scon}(t\Lambda;\cdot)$
be defined by (4.4).
If 
there is a unique
real number
$\,\scr f\,\,(\alpha)$
such that
$$
\sigma_{\radius}
\big(
\,
\zeta_{\alpha}^{\dyn\text{-}\scon}(\,\scr f\,\,(\alpha)\,\Lambda;\cdot)
\,
\big)
=
1\,,
$$
then
 $$
 \,\scr f\,\,(\alpha)
 =
 \inf
 \Big\{
 t\in\Bbb R
 \,\Big|\,
 \sigma_{\radius}(\,\zeta_{\alpha}^{\dyn\text{-}\scon}(t\Lambda)\,)
 \ge
 1
 \Big\}
$$
and this number is less than $f_{\mu}(\alpha)$,
i\.e\.
$$
\,\scr f\,\,(\alpha)
\le
f_{\mu}(\alpha)\,.
\tag4.5
$$
\endproclaim
\noindent{\it Proof}\newline
\noindent
Indeed, if
$\alpha
\not\in
-\beta'(\Bbb R)$.
then it is well-known that
that 
for all $\bold i\in\Sigma^{*}$, we have
$\frac{\log p_{\bold i}}{\log \diam K_{\bold i}}\not=\alpha$
(see, for example, [Pa])
This implies that if
$\alpha
\not\in
-\beta'(\Bbb R)$,
then the sum
 $$
    \sum
 \Sb
 |\bold i|=n\\
 {}\\
 \frac{\log p_{\bold i}}{\log \diam K_{\bold i}}=\alpha
 \endSb  \,\,
  \sup_{\bold u\in[\bold i]}
  \,\,
 \exp
 \,\,
 \sum_{k=0}^{n-1}t\Lambda S^{k}\bold u
 \tag4.6
$$
is the empty sum
and therefore equal to $0$
 for all $n$ and all $t$,
whence
$\zeta_{\alpha}^{\dyn\text{-}\scon}(t\Lambda;z)
  =
   \sum_{n}
   \frac{z^{n}}{n}
   0
   =
   0$
  for all $z$ and all $t$.
It follows immediately
this
 that if
$\alpha
\not\in
-\beta'(\Bbb R)$,
then
  $\sigma_{\radius}
(
\,
\zeta_{\alpha}^{\dyn\text{-}\scon}(t\Lambda;\cdot)
\,
)
=
\infty$ for all $t$,
whence
 $\,\scr f\,\,(\alpha)
 =
 \inf
 \{
 t\in\Bbb R
 \,|\,
 \sigma_{\radius}(\,\zeta_{\alpha}^{\dyn\text{-}\scon}(t\Lambda)\,)
 \ge
 1
 \}
 =
 -\infty$,
and inequality (4.5) is therefore trivially satisfied.
On the other hand, if
$\alpha
\in
-\beta'(\Bbb R)$,
then it follows from [CaMa,Fa1,Pa]
that there we can find  a (unique)
$q\in\Bbb R$ with 
$f_{\mu}(\alpha)
=
\alpha q+\beta(q)$.
It is also well-known, see, for example, [Bar,Fa2],
that there is a 
constant $c>0$ such that
for all positive integers $n$ and all $\bold i$ with $|\bold i|=n$
and all $\bold u\in[\bold i]$,
we have
$\frac{1}{c}
\le
\frac{|DS_{\bold i}(\pi S^{n}\bold u|)}{\diam K_{\bold i}}
\le 
c$.
This clearly implies that there is a constant $C$ 
such that
for all positive integers $n$ and all $\bold i$ with $|\bold i|=n$
and all $\bold u\in[\bold i]$,
we have
$\alpha q
\log
|DS_{\bold i}(\pi S^{n}\bold u)|
\le
C+\alpha q\log \diam K_{\bold i}$.
Since also
$\sum_{k=0}^{n-1}\Lambda S^{k}\bold u
=
\log |DS_{\bold i}(\pi S^{n}\bold u)|$
for all positive integers $n$ and all $\bold i$ with $|\bold i|=n$
and all $\bold u\in[\bold i]$,
we therefore conclude  that
 $$
 \align
   \sum
 \Sb
 |\bold i|=n\\
 {}\\
 \frac{\log p_{\bold i}}{\log \diam K_{\bold i}}=\alpha
 \endSb  \,\,
 &\sup_{\bold u\in[\bold i]}
  \,\,
 \exp
 \,\,
 \sum_{k=0}^{n-1}f_{\mu}(\alpha)\Lambda S^{k}\bold u\\
&{}\\
 &=
    \sum
 \Sb
 |\bold i|=n\\
 {}\\
 \frac{\log p_{\bold i}}{\log \diam K_{\bold i}}=\alpha
 \endSb  \,\,
  \sup_{\bold u\in[\bold i]}
  \,\,
 \exp
 \,\,
 \sum_{k=0}^{n-1}(\alpha q+\beta(q))\Lambda S^{k}\bold u\\
&{}\\
 &=
   \sum
 \Sb
 |\bold i|=n\\
 {}\\
 \frac{\log p_{\bold i}}{\log \diam K_{\bold i}}=\alpha
 \endSb  \,\,
  \sup_{\bold u\in[\bold i]}
  \,\,
 \exp
 \,\,
 \Bigg(
 \alpha q\sum_{k=0}^{n-1}\Lambda S^{k}\bold u
 +
 \sum_{k=0}^{n-1}\beta(q)\Lambda S^{k}\bold u
 \Bigg)\\
&{}\\
 &=
   \sum
 \Sb
 |\bold i|=n\\
 {}\\
 \frac{\log p_{\bold i}}{\log \diam K_{\bold i}}=\alpha
 \endSb  \,\,
  \sup_{\bold u\in[\bold i]}
  \,\,
 \exp
 \,\,
 \Bigg(
 \alpha q\log |DS_{\bold i}(\pi S^{n}\bold u)|
 +
 \sum_{k=0}^{n-1}\beta(q)\Lambda S^{k}\bold u
 \Bigg)\\
&{}\\ 
&\le
   \sum
 \Sb
 |\bold i|=n\\
 {}\\
 \frac{\log p_{\bold i}}{\log \diam K_{\bold i}}=\alpha
 \endSb  \,\,
  \sup_{\bold u\in[\bold i]}
  \,\,
 \exp
 \,\,
 \Bigg(
 C+\alpha q\log \diam K_{\bold i}
 +
 \sum_{k=0}^{n-1}\beta(q)\Lambda S^{k}\bold u
 \Bigg)\\
&{}\\ 
&=
   \sum
 \Sb
 |\bold i|=n\\
 {}\\
 \frac{\log p_{\bold i}}{\log \diam K_{\bold i}}=\alpha
 \endSb  \,\,
  \sup_{\bold u\in[\bold i]}
  \,\,
 \exp
 \,\,
 \Bigg(
 C+q\log p_{\bold i}
 +
 \sum_{k=0}^{n-1}\beta(q)\Lambda S^{k}\bold u
 \Bigg)\\ 
&{}\\ 
&=
   \sum
 \Sb
 |\bold i|=n\\
 {}\\
 \frac{\log p_{\bold i}}{\log \diam K_{\bold i}}=\alpha
 \endSb  \,\,
  \sup_{\bold u\in[\bold i]}
  \,\,
 \exp
 \,\,
 \Bigg(
 C
 +
  \sum_{k=0}^{n-1}q\Phi S^{k}\bold u
 +
 \sum_{k=0}^{n-1}\beta(q)\Lambda S^{k}\bold u
 \Bigg)\\  
&{}\\ 
&=
e^{C}
   \sum
 \Sb
 |\bold i|=n
 \endSb  \,\,
  \sup_{\bold u\in[\bold i]}
  \,\,
 \exp
 \,\,
 \Bigg(
  \sum_{k=0}^{n-1}(q\Phi +\beta(q)\Lambda)S^{k}\bold u
  \Bigg)\,,
 \endalign
 $$
whence
 $$
\align
&\sum_{n}
 \left|
 \frac{z^{n}}{n}
 \,\,
   \sum
 \Sb
 |\bold i|=n\\
 {}\\
 \frac{\log p_{\bold i}}{\log \diam K_{\bold i}}=\alpha
 \endSb  \,\,
 \sup_{\bold u\in[\bold i]}
  \,\,
 \exp
 \,\,
 \sum_{k=0}^{n-1}f_{\mu}(\alpha)\Lambda S^{k}\bold u
 \right|
 \\
&\qquad\qquad
 \qquad\qquad
 \le
e^{C}
 \sum_{n}
 \left|
 \frac{z^{n}}{n}
 \,\,
   \sum
 \Sb
 |\bold i|=n
 \endSb  \,\,
  \sup_{\bold u\in[\bold i]}
  \,\,
 \exp
 \,\,
 \Bigg(
  \sum_{k=0}^{n-1}(q\Phi +\beta(q)\Lambda)S^{k}\bold u
  \Bigg)
  \right|\,,
\endalign
 $$
for all complex numbers $z$.
We immediately conclude from this that
$1
=
\sigma_{\radius}
 \big(
 \,
 \zeta^{\dyn}(q\Phi+\beta(q)\Lambda;\cdot)
 \,
 \big)
 \le
  \sigma_{\radius}
 \big(
 \,
 \zeta_{\alpha}^{\dyn\text{-}\scon}(f_{\mu}(\alpha)\Lambda;\cdot)
 \,
 \big)$,
 whence 
$\,\scr f\,\,(\alpha)
 =
 \inf
 \{
 t\in\Bbb R
 \,|\,
 \sigma_{\radius}(\,\zeta_{\alpha}^{\dyn\text{-}\scon}(t\Lambda)\,)
 \ge
 1
 \}
 \le 
 f_{\mu}(\alpha)$.
This proves (4.5).
\hfill$\square$

\bigskip

\noindent
However, it is also clear that we, in general, do not
have equality in (4.5).
This is the content of the next observation.

\bigskip

\proclaim{Observation 4.2}
Let $\mu$ 
be the self-conformal measure defined by (2.2)
and let $\Lambda:\Sigma^{\Bbb N}\to\Bbb R$
be the scaling function defined 
by (2.5).
For $\alpha,t\in\Bbb R$,
let
$\zeta_{\alpha}^{\dyn\text{-}\scon}(t\Lambda;\cdot)$
be defined by (4.4).
If 
there is a unique
real number
$\,\scr f\,\,(\alpha)$
such that
$$
\sigma_{\radius}
\big(
\,
\zeta_{\alpha}^{\dyn\text{-}\scon}(\,\scr f\,\,(\alpha)\,\Lambda;\cdot)
\,
\big)
=
1\,,
$$
then
 $$
 \align
 &{}\\
 \,\scr f\,\,(\alpha)
 &=
  -\infty
  <
  0
  <
  f_{\mu}(\alpha)
  \,\,\,\,\\
  &{}\\
 &\qquad
  \text{
  for all 
  except at most countably many 
  $\alpha
  \in 
  -\beta'(\Bbb R)$.
  }\\
  &{}
 \tag4.7
   \endalign
  $$
\endproclaim
\noindent{\it Proof}\newline
\noindent
Indeed, the set 
$\{\frac{\log p_{\bold i}}{\log \diam K_{\bold i}}
 \,|\,
 \bold i\in\Sigma^{*}
 \}$ 
is clearly countable 
(because $\Sigma^{*}$ is countable)
and 
if 
$\alpha\in
 \Bbb R
 \setminus
\{\frac{\log p_{\bold i}}{\log r_{\bold i}}
 \,|\,
 \bold i\in\Sigma^{*}
 \}$, then
the sum
$ \sum_{
 |\bold i|=n\,,\,
 \frac{\log p_{\bold i}}{\log \diam K_{\bold i}}=\alpha
 }
  \sup_{\bold u\in[\bold i]}
 \exp
 \sum_{k=0}^{n-1}t\Lambda S^{k}\bold u$
is the empty sum
and therefore equal to $0$
 for all $n$ and all $t$.
 It follows from this,
 using an argument similar to the 
reasoning
following (4.6),
that
$\,\scr f\,\,(\alpha)
 =
 \inf
 \{
 t\in\Bbb R
 \,|\,
 \sigma_{\radius}(\,\zeta_{\alpha}^{\dyn\text{-}\scon}(t\Lambda)\,)
 \ge
 1
 \}
 =
 -\infty$.
Since it also
follows from [CaMa,Fa1,Pa]
 that
 $f_{\mu}(\alpha)>0$
for all
 $\alpha
 \in 
-\beta'(\Bbb R)$,
 we therefore conclude that
 $ \,\scr f\,\,(\alpha)
=
  -\infty
  <
  0
  <
  f_{\mu}(\alpha)$
  for all 
  except at most countably many 
  $\alpha
  \in 
  -\beta'(\Bbb R)$.
  \hfill$\square$

\bigskip

It follows from the above discussion that
while 
the definition
of $\zeta_{\alpha}^{\dyn\text{-}\scon}(s)$
is \lq\lq natural",
it is not 
does not encode
sufficient
information
allowing us to recover the multifractal spectrum
$f_{\mu}(\alpha)$.
The reason for the strict inequality in (4.7) 
is, of course, clear:
even though
there are no
strings
$\bold i\in\Sigma^{*}$
for which
the ratio
$\frac{\log p_{\bold i}}{\log \diam K_{\bold i}}$ equals $\alpha$
if
$\alpha
\in
-\beta'(\Bbb R)
  \setminus
 \{\frac{\log p_{\bold i}}{\log \diam K_{\bold i}}
 \,|\,
 \bold i\in\Sigma^{*}
 \}$,
there are
nevertheless
many sequences $(\bold i_{n})_{n}$ of strings
$\bold i_{n}\in\Sigma^{*}$
for which the sequence
of
ratios 
$(\frac{\log p_{\bold i_{n}}}{\log \diam K_{\bold i_{n}}})_{n}$ converges to $\alpha$.
In order to capture this,
it is necessary
to 
ensure
that
those
 strings
$\bold i$
for which
the ratio
$\frac{\log p_{\bold i}}{\log \diam K_{\bold i}}$ 
is 
\lq\lq close"
 to $\alpha$
are
also included in the series
defining 
 the multifractal zeta-function.
For this reason, we
modify the definition of
$\zeta_{\alpha}^{\dyn\text{-}\scon}$
and
introduce
a
self-conformal multifractal zeta-function
obtained by
replacing the 
original 
small \lq\lq target" set
$\{\alpha\}$
by a larger
 \lq\lq target" set
$I$
(for example,
 we may
choose the enlarged
 \lq\lq target" set
$I$
to be a non-degenerate  interval centered at $\alpha$).
In order to
make this idea precise we proceed as follows.
 For a closed interval
$I$,
we 
define the self-conformal
multifractal zeta-function $\zeta^{\dyn\text{-}\scon}_{I}$ 
of a continuous function $\varphi:\Sigma^{\Bbb N}\to\Bbb R$
by
 $$
 \zeta_{I}^{\dyn\text{-}\scon}(\varphi;z)
  =
   \sum_{n}
   \,\,
   \frac{z^{n}}{n}
 \left(
   \sum
 \Sb
 |\bold i|=n\\
 {}\\
 \frac{\log p_{\bold i}}{\log \diam K_{\bold i}}\in I
 \endSb  \,\,
  \sup_{\bold u\in[\bold i]}
  \,\,
 \exp
 \,\,
 \sum_{k=0}^{n-1}\varphi S^{k}\bold u
\right)
\tag4.8
 $$
for those complex numbers $z$ for which the series converges.
Observe that if $I=\{\alpha\}$, then
 $$
  \zeta^{\dyn\text{-}\scon}_{I}(\varphi;z)
=
 \zeta^{\dyn\text{-}\scon}_{\alpha}(\varphi;z)\,.
 $$
We can now 
proceed in two equally natural ways.
Either, we can 
consider
a family
of enlarged
\lq\lq target" sets 
shrinking to the original 
main \lq\lq target" $\{\alpha\}$;
this approach will be referred to as the
shrinking target approach.
Or, alternatively,
we can consider
a
fixed enlarge 
\lq\lq target" set
and regard this as our
original main \lq\lq target";
this approach will be referred to as the
fixed target approach.
We now
discuss these approaches in more detail.

\bigskip

\noindent
{\it (1) The shrinking target approach.
}
For a given (small) \lq\lq target"
$\{\alpha\}$,
we
consider
the following family
$\big(\,[\alpha-r,\alpha+r]\,\big)_{r>0}$
of enlarged
\lq\lq target" sets $[\alpha-r,\alpha+r]$
shrinking  to the original 
main \lq\lq target" $\{\alpha\}$ as $r\searrow0$,
and 
attempt
to
relate the limiting behaviour of the 
radius of convergence
of 
$\zeta_{[\alpha-r,\alpha+r]}^{\dyn\text{-}\scon}(t\Lambda;\cdot)$
as $r\searrow 0$
to the multifractal spectrum
$f_{\mu}(\alpha)$ at $\alpha$.
The next result,
which is an application of one of our main results (see Theorem 6.1),
shows that 
the multifractal zeta-functions
$\zeta_{[\alpha-r,\alpha+r]}^{\dyn\text{-}\scon}(t\Lambda;\cdot)$
encode sufficient
information allowing us
to recover 
the multifractal spectra
$f_{\mu}(\alpha)$
by letting 
$r\searrow 0$.

\bigskip

\proclaim{Theorem 4.1. Shrinking targets}
Let
 $\mu$ be the 
 self-conformal measure defined by (2.2)
 and
 let $\Lambda:\Sigma^{\Bbb N}\to\Bbb R$
 be the scaling function defined by (2.5).
For  $\alpha\in\Bbb R$, $r>0$
and $t\in\Bbb R$,
let $\zeta_{[\alpha-r,\alpha+r]}^{\dyn\text{-}\scon}(t\Lambda;\cdot)$ be defined by (4.8).
\roster
\item"(1)"
There is a unique real number
$\,\scr f\,\,(\alpha)$
such that
$$
\lim_{r\searrow 0}
\sigma_{\radius}
\big(
\,
\zeta_{[\alpha-r,\alpha+r]}^{\dyn\text{-}\scon}(\,\scr f\,\,(\alpha)\,\Lambda;\cdot)
\,
\big)
=
1\,,
$$
\item"(2)"
If the OSC is satisfied, then we have
 $$
 \,\scr f\,\,(\alpha)
  =
 f_{\mu}(\alpha)\,.
  $$
\endroster
\endproclaim
\noindent{\it Proof}\newline
\noindent
This result is a special case of Theorem 6.1.
\hfill$\square$

\bigskip

\noindent
We note  that
Theorem 4.1 has 
a very clear resemblance to Bowen's formula in Theorem C.

\bigskip

\noindent
{\it (2) The fixed target approach
}
Alternatively,
we
can keep the
enlarged 
\lq\lq target" set $I$
fixed
and 
attempt
to relate the
radius of convergence
of the
multifractal zeta-function 
$\zeta_{I}^{\dyn\text{-}\scon}(t\Lambda;\cdot)$
associated with the enlarger \lq\lq target" set $I$
to the 
values
 of the multifractal spectrum
$f_{\mu}(\alpha)$ 
for $\alpha\in I$.
Of course, inequality (4.7) shows that
if the \lq\lq target" 
set $I$ is 
\lq\lq too small",
then this is not possible.
However,
if
the enlarger \lq\lq target" set $I$ 
satisfies 
a mild
non-degeneracy
condition, namely condition (4.9), 
guaranteeing
that $I$ is 
sufficiently
\lq\lq big",
then
the next result,
which is also 
an application of one of our main results (see Theorem 6.1),
shows that 
this is possible. 
More
precisely the result
shows that
if
the enlarger \lq\lq target" set $I$ 
satisfies condition 
 (4.9),
then
the multifractal zeta-function
$\zeta_{I}^{\dyn\text{-}\scon}(t\Lambda;\cdot)$
associated with the enlarger \lq\lq target" set $I$
encode sufficient
information allowing us
to recover 
the suprema
$\sup_{\alpha\in I}
 f_{\mu}(\alpha)$
 of the 
 multifractal spectrum $ f_{\mu}(\alpha)$
for $\alpha\in I$.

\bigskip

\proclaim{Theorem 4.2. Fixed targets}
Let
 $\mu$ be the 
 self-conformal measure defined by (2.2)
 and
 let $\Lambda:\Sigma^{\Bbb N}\to\Bbb R$
 be the scaling function defined by (2.5).
 For  a closed interval $I$
and $t\in\Bbb R$,
let 
$\zeta_{I}^{\dyn\text{-}\scon}(t\Lambda;\cdot)$
 be defined by (4.8).
Assume that
 $$
 \overset{\,\circ}\to{I}
 \cap
 \big(-\beta'(\Bbb R)\big)
 \not=
 \varnothing
 \tag4.9
 $$
(where $\overset{\,\circ}\to{I}$ denotes the interior of $I$).
\roster
\item"(1)"
There is a unique real number
$\scr F\,\,(I)$
such that
$$
\sigma_{\radius}
\big(
\,
\zeta_{I}^{\dyn\text{-}\scon}(\scr F\,\,(I)\,\Lambda;\cdot)
\,
\big)
=
1\,,
$$
\item"(2)"
If the OSC is satisfied, the we have
 $$
 \scr F\,\,(I)
 =
 \sup_{\alpha\in I}
 f_{\mu}(\alpha)\,.
 $$
\endroster
 \endproclaim
 \noindent{\it Proof}\newline
\noindent
This result is a special case of Theorem 6.1.
\hfill$\square$

\bigskip

\noindent
As with Theorem 4.1, we also note that
Theorem 4.2 has 
a very clear resemblance to Bowen's formula in Theorem C.

We emphasise that Theorem 4.1 and Theorem 4.2
are
presented
in order to motive this work
and
are special cases
of the
substantially
more general and abstract 
theory 
of dynamical multifractal zeta-function
developed in this paper.

The next section, i\.e\. Section 5, describes the 
general framework
developed in this paper and lists our main results.
In Section 6 we will discuss a number of examples, including,
mixed and non-mixed
multifractal spectra 
of self-conformal measures,
and multifractal spectra of Birkhoff ergodic averages.

  \bigskip


\heading{5. Statements of the main results.}\endheading

We also  denote the family of Borel probability measures on 
$\Sigma^{\Bbb N}$
and the 
 family of shift invariant Borel probability measures on 
$\Sigma^{\Bbb N}$
by $\Cal P(\Sigma^{\Bbb N})$
and
$\Cal P_{S}(\Sigma^{\Bbb N})$, respectively, i\.e\.
we write
 $$
 \align
 \Cal P(\Sigma^{\Bbb N})
&=
 \Big\{
 \mu\,\Big|\,
 \text{
 $\mu$ is a Borel probability measures on 
$\Sigma^{\Bbb N}$
 }
 \Big\}\,,\\
 \Cal P_{S}(\Sigma^{\Bbb N})
&=
 \Big\{
 \mu\,\Big|\,
 \text{
 $\mu$ is a shift invariant Borel probability measures on 
$\Sigma^{\Bbb N}$
 }
 \Big\}\,;\\
 \endalign
 $$
we will always
equip
$\Cal P(\Sigma^{\Bbb N})$ 
and
$\Cal P_{S}(\Sigma^{\Bbb N})$ with the weak topology.
Fix a metric space $X$ and a
continuous map
$U:\Cal P(\Sigma^{\Bbb N})\to X$.
%
%
%
%
%
%
%
%
%
%
 For a positive integer $n$,
 let
$L_{n}:\Sigma^{\Bbb N}\to\Cal P(\Sigma^{\Bbb N})$ be defined 
by
 $$
 L_{n}\bold i
 =
 \frac{1}{n}\sum_{k=0}^{n-1}\delta_{S^{k}\bold i}\,.
 \tag5.1
 $$
 We can now define the 
multifractal pressure and zeta-function
associated with the space $X$ and the map $U$,

 \bigskip
 
 \proclaim{Definition.
 The multifractal pressure 
 $ \underline P_{C}^{U}(\varphi)$
 and
 $ \overline P_{C}^{U}(\varphi)$
associated with the space $X$ and the map $U$}
Let
$\varphi:\Sigma^{\Bbb N}\to\Bbb R$
  be a continuous map.
    For $C\subseteq X$,
we define the lower and upper
mutifractal pressure
of 
 $\varphi$
associated with the space $X$ and the map $U$ and  by
 $$
 \align
  \underline P_{C}^{U}(\varphi)
 &=
 \,
  \liminf_{n}
 \,
  \,\,
  \frac{1}{n}
 \,\,
 \log
   \sum
  \Sb
  |\bold i|=n\\
  {}\\
  UL_{n}[\bold i]\subseteq C
  \endSb
   \,\,
  \sup_{\bold u\in[\bold i]}
  \,\,
 \exp
 \,\,
 \sum_{k=0}^{n-1}\varphi S^{k}\bold u\,,\\
   \overline P_{C}^{U}(\varphi)
 &=
  \limsup_{n}
  \,\,
  \frac{1}{n}
 \,\,
 \log
   \sum
  \Sb
 |\bold i|=n\\
  {}\\
  UL_{n}[\bold i]\subseteq C
  \endSb
 \,\,
  \sup_{\bold u\in[\bold i]}
  \,\,
 \exp
 \,\,
 \sum_{k=0}^{n-1}\varphi S^{k}\bold u\,.
 \endalign
 $$
If $ \underline P_{C}^{U}(\varphi)$ and $ \overline P_{C}^{U}(\varphi)$ coincide, then we write
$P_{C}^{U}(\varphi)$ for their common value, i\.e\. we write
$ P_{C}^{U}(\varphi)
=
 \underline P_{C}^{U}(\varphi)
 =
  \overline P_{C}^{U}(\varphi)$.

\endproclaim

 \bigskip
 
 \proclaim{Definition.
 The dynamical multifractal zeta-function 
 $\zeta_{C}^{\dyn,U}(\varphi;\cdot)$
 associated with the space $X$ and the map $U$}
Let
$\varphi:\Sigma^{\Bbb N}\to\Bbb R$
  be a continuous map.
    For $C\subseteq X$,
we define the 
dynamical
multifractal
 zeta-function $\zeta_{C}^{\dyn,U}(\varphi;\cdot)$
associated with the space $X$ and the map $U$ by
 $$
  \zeta_{C}^{\dyn,U}(\varphi;z)
  =
   \sum_{n}
   \,\,
   \frac{z^{n}}{n}
 \left(
   \sum
  \Sb
 |\bold i|=n\\
  {}\\
  UL_{n}[\bold i]\subseteq C
  \endSb
  \,\,
  \sup_{\bold u\in[\bold i]}
  \,\,
 \exp
 \,\,
 \sum_{k=0}^{n-1}\varphi S^{k}\bold u
\right)
 $$
for those complex numbers $z$ for which the series converges.

\endproclaim

\bigskip

\noindent
{\bf Remark.}
It is clear that if $C=X$, then the 
multifractal 
\lq\lq constraint"
$UL_{n}[\bold i]\subseteq C$
is vacuously
satisfied,
and the 
multifractal pressure and 
dynamical multifractal zeta-function
reduce to the usual pressure and the usual
dynamical zeta-function, i\.e\.
 $$
 \underline P_{X}^{U}(\varphi)
 =
 \overline P_{X}^{U}(\varphi)
 =
  P(\varphi)
 $$
and
 $$
 \zeta_{X}^{\dyn,U}(\varphi;\cdot)
 =
 \zeta^{\dyn}(\varphi;\cdot)\,.
 $$

\bigskip

Before developing
the theory of the multifractal pressure and the 
multifractal 
zeta-functions further
we make to following two simple observations.
Firstly, we note (see Proposition 5.1) that
the expected
relationship between the
the multifractal pressure and the radius of 
convergence of the multifractal zeta-function
holds.
Secondly,
we would expect any dynamically meaningful
 theory of
dynamical multifractal zeta-functions to
lead to
multifractal Bowen formulas.
For this to hold, we must, at the very least,
ensure that
there  are unique solutions to 
the 
relevant
multifractal Bowen equations. i\.e\.
 we must ensure
 that
there is are unique real numbers
 $\,\,\scri f(C)$ and
 $\,\,\scri F\,(C)$
solving the following
 equations, namely,
 $$
 \aligned
 \limsup_{r\searrow 0}\overline P_{B(C,r)}^{U}(\,\,\,\scri f(C)\,\Phi)
 &=
 0\,,\\
 \overline P_{C}^{U}(\,\scri F\,(C)\,\Phi)
 &=
 0\,,
 \endaligned
 \tag5.2
 $$ 
That there are 
unique numbers 
 $\,\,\scri f(C)$ and
 $\,\,\scri F\,(C)$
 satisfying (5.2)
 is our second simple observation (see Proposition 5.2).

\bigskip

\proclaim{Proposition 5.1. Radius of convergence}
Let $X$ be a metric space and let $U:\Cal P(\Sigma^{\Bbb N})\to X$ be 
continuous with respect to the weak topology.
Let $C\subseteq X$ be a  subset of $X$.
Fix a continuous function $\varphi:\Sigma^{\Bbb N}\to\Bbb R$.
We have
 $$
 -
 \log
 \sigma_{\radius}
 \big(
 \,
 \zeta_{C}^{\dyn,U}(\varphi;\cdot)
 \,
 \big)
 =
 \overline P_{C}^{U}(\varphi)\,.
 $$
\endproclaim 
\noindent{\it Proof}\newline
\noindent
This follows immediately from the fact that
if
$(a_{n})_{n}$ is a sequence of complex
numbers
and
$f(z)=\sum_{n}a_{n}z^{n}$, then
$\sigma_{\radius}(f)
 =
 \frac{1}{\limsup_{n}|a_{n}|^{\frac{1}{n}}}$.
\hfill$\square$

\bigskip

\proclaim{Proposition 5.2. 
Continuity and monotonicity
of the multifractal pressure}
Let $X$ be a metric space and let $U:\Cal P(\Sigma^{\Bbb N})\to X$ be 
continuous with respect to the weak topology.
Let $C\subseteq X$ be a  subset of $X$.
Fix a continuous
map
$\Phi:\Sigma^{\Bbb N}\to\Bbb R$
with
$\Phi<0$.
Let $C$ be a subset of $X$.
Then the functions
$t\to
 \limsup_{r\searrow 0}\overline P_{B(C,r)}^{U}(t\Phi)$
and 
$t
 \to
 \overline P_{C}^{U}(t\Phi)$,
 where
 $t\in\Bbb R$, 
are
continuous,
strictly decreasing and convex with
 $\lim_{t\to-\infty}
\limsup_{r\searrow 0}\overline P_{B(C,r)}^{U}(t\Phi)
 =
\infty$
and
$\lim_{t\to\infty}
\limsup_{r\searrow 0}\overline P_{B(C,r)}^{U}(t\Phi)
 =
-\infty$,
and
$\lim_{t\to-\infty}
\overline P_{C}^{U}(t\Phi) 
 =
\infty$
and
$\lim_{t\to\infty}
\overline P_{C}^{U}(t\Phi) 
 =
-\infty$.
In particular, there are unique real numbers
$\,\,\scri f(C)$
and
$\scri F\,(C)$ 
such that
 $$
 \align
 \limsup_{r\searrow 0}\overline P_{B(C,r)}^{U}(\,\,\,\scri f(C)\,\Phi)
 &=
 0\,,
 \tag5.3\\
 \overline P_{C}^{U}(\,\scri F\,(C)\,\Phi)
&=
0\,;
\tag5.4
 \endalign
 $$
alternatively,
$\,\,\scri f(C)$
and
$\scri F\,(C)$ 
 are the  unique real numbers 
 such that
 $$
 \align
 \limsup_{r\searrow 0}
 \sigma_{\radius}
 \big(
 \,
 \zeta_{B(C,r)}^{U}(\,\,\,\scri f(C)\,\Phi;\cdot)
 \,
 \big)
 &=
 1\,,\\
 \sigma_{\radius}
 \big(
 \,
 \zeta_{C}^{U}(\,\scri F\,(C)\,\Phi;\cdot)
 \,
 \big)
&=
1\,.
 \endalign
 $$ 
\endproclaim
\noindent{\it Proof}\newline 
\noindent
This is not difficult to prove and 
for sake of brevity we have decided to omit the proof.
\hfill$\square$

\bigskip

We can now state our main results.
The results are divided into two 
parts
paralleling the discussion in Section 4.2
The first part 
(consisting of
Theorem 5.3 and Corollary 5.4)
presents our results in the shrinking target setting,
and the 
second part 
consisting of 
Theorem 5.5 and Corollary 5.6)
presents our results in the fixed target setting.
More precisely,
 in the
shrinking target
setting,
Theorem 5.3  provide
a
variational principle for the 
multifractal pressure
and 
Corollary 5.4
provide a
variational principle
for
the solution
$\,\,\scri f(C)$ to the multifractal Bowen equation (5.3).

\bigskip

\proclaim{Theorem 5.3. 
The shrinking target variational principle
for the multifractal pressure}
Let $X$ be a metric space and let $U:\Cal P(\Sigma^{\Bbb N})\to X$ be 
continuous with respect to the weak topology.
Let $C\subseteq X$ be a  subset of $X$.
Fix a continuous function $\varphi:\Sigma^{\Bbb N}\to\Bbb R$.

\roster
\item"(1)"
We have
 $$
 \align
 \lim_{r\searrow0}
\underline P_{B(C,r)}^{U}(\varphi)
&=
\overline P_{B(C,r)}^{U}(\varphi)
=
\sup
\Sb
\mu\in\Cal P_{S}(\Sigma^{\Bbb N})\\
{}\\
U\mu\in \overline C
\endSb
\Bigg(
h(\mu)+\int\varphi\,d\mu
\Bigg)\,.\\
\endalign
$$ 
\item"(2)"
We have
 $$
 \align
 \lim_{r\searrow0}
-\log
\sigma_{\radius}\big(\,\zeta_{B(C,r)}^{\dyn,U}(\varphi;\cdot)\,\big)
&=
\sup
\Sb
\mu\in\Cal P_{S}(\Sigma^{\Bbb N})\\
{}\\
U\mu\in \overline C
\endSb
\Bigg(
h(\mu)+\int\varphi\,d\mu
\Bigg)\,.
\endalign
$$ 
\endroster
\endproclaim

\bigskip

\noindent
Theorem 5.3 is proved in Section 8.

\bigskip

\proclaim{Corollary 5.4.
The shrinking target
multifractal Bowen equation}
Let $X$ be a metric space and let $U:\Cal P(\Sigma^{\Bbb N})\to X$ be 
continuous with respect to the weak topology.
Let $C\subseteq X$ be a  subset of $X$.
Fix a continuous function  $\Phi:\Sigma^{\Bbb N}\to\Bbb R$ 
with
$\Phi<0$
and let
$\,\,\scri f(C)$
be the unique real number
such that
 $$
 \align
 \limsup_{r\searrow 0}\overline P_{B(C,r)}^{U}(\,\,\,\scri f(C)\,\Phi)
 &=
 0\,;
 \endalign
 $$ 
alternatively, 
 $\,\,\scri f(C)$
is the unique real number
such that
 $$
 \align
 \limsup_{r\searrow 0}
 \sigma_{\radius}
 \big(
 \,
 \zeta_{B(C,r)}^{U}(\,\,\scri f(C)\,\Phi;\cdot)
 \,
 \big)
 &=
 1\,.
 \endalign
 $$ 
Then
 $$
\scri f(C)
=
\sup
\Sb
\mu\in\Cal P_{S}(\Sigma^{\Bbb N})\\
{}\\
U\mu\in \overline C
\endSb
-
\frac{h(\mu)}{\int\Phi\,d\mu}\,.
$$ 
\endproclaim 
\noindent{\it Proof}\newline
\noindent
It follows from Theorem 5.3 and the definition of
$\,\,\scri f(C)$
that
 $$
 \align
\sup
\Sb
\mu\in\Cal P_{S}(\Sigma^{\Bbb N})\\
{}\\
U\mu\in \overline C
\endSb
\Bigg(
h(\mu)+\,\,\scri f(C)\int\Phi\,d\mu
\Bigg)
&=
 \limsup_{r\searrow 0}\overline P_{B(C,r)}^{U}(\,\,\,\scri f(C)\,\Phi)
 =
0\,.
\tag5.5
\endalign
$$ 
The desired formula for $\,\,\scri f(C)$
follows easily from (5.5).
\hfill$\square$

\bigskip

\newpage

\noindent
Of course, if the set $C$ is 
\lq\lq too small",
then
it follows from the discussion in Section 4.2
that we, in general, cannot expect any
meaningful results 
in the fixed target setting.
However, if the set $C$ satisfies a 
non-degeneracy 
condition guaranteeing that
it is not
\lq\lq too small"
(namely condition (5.6) below),
then 
meaningful results can be obtained
in the fixed target setting.
This is the contents of
Theorem 2.5 and Corollary 2.6 below.
Indeed,
Theorem 2.5 and Corollary 2.6 provide
variational principles for the 
multifractal pressure
and for the solution
$\,\,\scri F\,(C)$ to the multifractal Bowen equation (5.4)
 in the
fixed target
setting.

\bigskip

\proclaim{Theorem 5.5.
The fixed target variational principle
for the multifractal pressure}
Let $X$ be a normed vector space.
Let $\Gamma:\Cal P(\Sigma^{\Bbb N})\to X$
be continuous and affine
and let
$\Delta:\Cal P(\Sigma^{\Bbb N})\to \Bbb R$
be continuous and affine
with
$\Delta(\mu)\not=0$
for all $\mu\in\Cal P(\Sigma^{\Bbb N})$.
Define 
$U:\Cal P(\Sigma^{\Bbb N})\to X$
by
$U=\frac{\Gamma}{\Delta}$.
Let $C$ be a closed and convex subset of $X$ and assume that
 $$
 \overset{\,\circ}\to{C}
 \cap
 \,
 U\big(\,\Cal P_{S}(\Sigma^{\Bbb N})\,\big)
 \not=
 \varnothing\,.
 \tag5.6
 $$
\roster
\item"(1)"
We have
$$
\align
{}\qquad\qquad
  \quad\,\,\,\,
  \,\,\,\,\,\,\,
P_{C}^{U}(\varphi)
&=
\sup
\Sb
\mu\in \Cal P_{S}(\Sigma^{\Bbb N})\\
{}\\
U\mu\in C
\endSb
\Bigg(
h(\mu)+\int\varphi\,d\mu
\Bigg)
=
\sup
\Sb
\mu\in \Cal P_{S}(\Sigma^{\Bbb N})\\
U\mu\in \overset{\,\circ}\to{C}
\endSb
\Bigg(
h(\mu)+\int\varphi\,d\mu
\Bigg)\,.
\endalign
$$

 \item"(2)"
We have
$$
\align
-\log
\sigma_{\radius}\big(\,\zeta_{C}^{\dyn,U}(\varphi;\cdot)\,\big)
&=
\sup
\Sb
\mu\in \Cal P_{S}(\Sigma^{\Bbb N})\\
{}\\
U\mu\in C
\endSb
\Bigg(
h(\mu)+\int\varphi\,d\mu
\Bigg)
=
\sup
\Sb
\mu\in \Cal P_{S}(\Sigma^{\Bbb N})\\
U\mu\in \overset{\,\circ}\to{C}
\endSb
\Bigg(
h(\mu)+\int\varphi\,d\mu
\Bigg)\,.
\endalign
$$
\endroster
\endproclaim

\bigskip

\noindent
Theorem 5.5 is proved in Section 9.

\bigskip

\proclaim{Corollary 5.6.
The fixed target
multifractal Bowen equation}
Let $X$ be a normed vector space.
Let $\Gamma:\Cal P(\Sigma^{\Bbb N})\to X$
be continuous and affine
and let
$\Delta:\Cal P(\Sigma^{\Bbb N})\to \Bbb R$
be continuous and affine
with
$\Delta(\mu)\not=0$
for all $\mu\in\Cal P(\Sigma^{\Bbb N})$.
Define 
$U:\Cal P(\Sigma^{\Bbb N})\to X$
by
$U=\frac{\Gamma}{\Delta}$.
Let $C$ be a closed and convex subset of $X$ and assume that
 $$
 \overset{\,\circ}\to{C}
 \cap
 \,
 U\big(\,\Cal P_{S}(\Sigma^{\Bbb N})\,\big)
 \not=
 \varnothing\,.
 $$
Let $\Phi:\Sigma^{\Bbb N}\to\Bbb R$ be continuous
with
$\Phi<0$.
Let 
$\scri F\,(C)$
be the unique real number
such that
 $$
 \align
 P_{C}^{U}(\,\scri F\,(C)\,\Phi)
 &=
 0\,;
 \endalign
 $$ 
 alternatively, 
 $\,\,\scri f(C)$
is the unique real number
such that
 $$
 \align
 \sigma_{\radius}
 \big(
 \,
 \zeta_{C}^{U}(\,\,\scri F\,(C)\,\Phi;\cdot)
 \,
 \big)
 &=
 1\,.
 \endalign
 $$ 
Then
 $$
\scri F\,(C)
=
\sup
\Sb
\mu\in\Cal P_{S}(\Sigma^{\Bbb N})\\
{}\\
U\mu\in C
\endSb
-
\frac{h(\mu)}{\int\Phi\,d\mu}\,.
$$ 
\endproclaim 
\noindent{\it Proof}\newline
\noindent
The proof is similar to the proof of Corollary 5.4 using
Theorem 5.5 and the definition of
$\scri F\,(C)$.
\hfill$\square$

\bigskip

\noindent
In the next section we will show that in many cases,
 the solutions 
$\,\,\scri f(C)$ 
and
$\,\,\scri F\,(C)$ to the multifractal Bowen equations 
(5.3)
and (5.4) coincide with 
the usual multifractal spectra.

  \bigskip
  
  \newpage


\centerline{\smc 6. Applications: }

\centerline{\smc multifractal spectra of measures}

\centerline{\smc  and}

\centerline{\smc  
multifractal spectra of ergodic Birkhoff averages}

 \medskip

We will now consider several of applications of Theorem 5.3 and Theorem 5.5
to multifractal spectra of measures and ergodic averages.
In particular, we consider the following examples:

\medskip

$\bullet$
 Section 6.1: Multifractal spectra of self-conformal measures.

\medskip

$\bullet$
 Section 6.2: Mixed multifractal spectra of self-conformal measures.

\medskip

$\bullet$
Section 6.3: Multifractal spectra of ergodic Birkhoff averages.

\medskip

{\bf 6.1.
Multifractal 
spectra of 
self-conformal measures.}
Fix a a conformal iterated function system 
 $(
 \,
 V
 \,,\,
 X
 \,,\,
 (S_{i})_{i=1,\ldots,N}
 \,
 )$
 and a 
 be a probability vector $(p_{1},\ldots,p_{N})$.
We let $K$ denote the self-conformal set defined by (2.1),
and we let
 $\mu$ denote
 the
self-conformal measure defined by
 (2.2).
We also recall that
the Hausdorff multifractal spectrum
$f_{\mu}$ of $\mu$ 
is defined
by
 $$
 \align
 f_{\mu}(\alpha)
&=
   \,\dim_{\Haus}
   \left\{x\in K
    \,\left|\,
     \lim_{r\searrow0}
     \frac
     {\log\mu B(x,r)}{\log r}
     =
     \alpha
      \right.
   \right\}\,,
  \endalign
 $$
for $\alpha\in\Bbb R$,
and that the 
multifractal spectrum
$f_{\mu}(\alpha)$
can be computed as follows, see, for example,
 [ArPa,CaMa,Pa].
Define
$\Phi,\Lambda:\Sigma^{\Bbb N}\to\Bbb R$ by
$\Phi(\bold i)=\log p_{i_{1}}$ and 
let
$\lambda:\Sigma^{\Bbb N}\to\Bbb R$
denote the 
scaling map defined in (2.5).
Finally, 
let 
$\beta(q)$
be the unique real number such that
 $$
 P\big(
 \,
 \beta(q)\Lambda
 +
 q\Phi
 \,
 \big)
 =
 0\,;
 \tag6.1
 $$
 alternatively,
 the function
 $\beta:\Bbb R\to\Bbb R$
is defined by
 $$
 \sigma_{\radius}
 \big(
 \,
 \zeta^{\dyn}(q\Phi+\beta(q)\Lambda;\cdot)
 \,
 \big)
 =
 1\,.
 \tag6.2
 $$
The multifractal spectrum $f_{\mu}(\alpha)$
can now be computer as follows. If the 
 OSC is satisfied, then it follows from
[ArPa,CaMa,Pa]
that
 $$
 f_{\mu}(\alpha)
 =
 \beta^{*}(\alpha)\,;
 \tag6.3
 $$
recall, that if
 if $\varphi:\Bbb R\to\Bbb R$ is a function,
 then the Legendre transform 
 $\varphi^{*}:\Bbb R\to[\-\infty,\infty]$
 of $\varphi$ is defined
 by
$ \varphi^{*}(x)
  =
  \inf_{y}
  (xy+\varphi(y))$.

Of course, in general, the limit
 $\lim_{r\searrow0}
     \frac
     {\log\mu B(x,r)}{\log r}$
     may not exist.
     Indeed, recently 
     Barreira \& Schmeling [BaSc]
     (see also
     Olsen \& Winter [OlWi1,OlWi2],
     Xiao, Wu \& Gao [XiWuGa]
     and 
Moran [Mo])
have shown that 
the set of divergence points, 
i\.e\. the set 
of points $x$ for which the 
limit
$\lim_{r\searrow0}
     \frac
     {\log \mu B(x,r)}{\log r}$
does not exist, typically is  highly
\lq\lq visible" and
\lq\lq observable", namely it has full Hausdorff dimension.
More precisely, it follows from 
[BaSc]
that if
the OSC is satisfied and $t$ denotes the Hausdorff 
dimension of $K$, then
 $$
  \Bigg\{
 x\in K
     \,\Bigg|\,
     \text{the expression}
    \,\,
     \frac
     {\log \mu B(x,r)}{\log r}
     \,\,
     \text{diverges as $r\searrow0$}
     \,\,
      \Bigg\}
 =
 \varnothing
 $$
provided $\mu$
is proportional to the $t$-dimensional Hausdorff measure restricted to $K$,
and
 $$
 \dim_{\Haus}
  \Bigg\{
 x\in K
     \,\Bigg|\,
     \text{the expression}
    \,\,
     \frac
     {\log \mu B(x,r)}{\log r}
     \,\,
     \text{diverges as $r\searrow0$}
     \,\,
      \Bigg\}
 =
 \dim_{\Haus} K
 $$ 
provided $\mu$
is 
not
proportional to the $t$-dimensional Hausdorff measure restricted to $K$.
This
 suggests that the set 
 of divergence points
 has a surprising rich and complex  
 fractal structure,
 and in order to explore this more
 carefully
  Olsen \& Winter [OlWi1,OlWi2]
 introduced various
 generalised multifractal spectra functions designed to 
 \lq\lq see"
 different sets of divergence points.
 In order to define these spectra 
 we introduce the following notation.
 If $M$ 
 is a  metric space
 and
 $\varphi:(0,\infty)\to M$ is a function, then we write
 $\acc_{r\searrow 0}f(r)$
 for the set of accumulation 
 points of $f$ as $r\searrow 0$, i\.e\.
  $$
  \underset {r\searrow0}\to\acc\,\,\varphi(r)
  =
  \Big\{x\in M
  \,\Big|\,
  \text{$x$ is an accumulation point of $f$ as $r\searrow 0$}
  \Big\}\,.
  $$
  In [OlWi1]
Olsen \& Winter
introduced and investigated
 the generalised
 Hausdorff multifractal spectrum
 $F_{\mu}$ of $\mu$ 
 defined
by
 $$
 \align
F_{\mu}(C)
&=
   \,\dim_{\Haus}
   \left\{x\in K
    \,\left|\,
     \,\,
     \underset {r\searrow0}\to\acc
 \,
     \frac
     {\log\mu B(x,r)}{\log r}
  \subseteq
 C
      \right.
      \,
   \right\}
  \endalign
 $$
 for 
 $C\subseteq \Bbb R$.
Note that the generalised spectrum
is a genuine extension of the
traditional multifractal spectrum 
$f_{\mu}(\alpha)$, namely if
$C=\{\alpha\}$ is a singleton consisting of the point $\alpha$, 
then clearly
$F_{\mu}(C)=f_{\mu}(\alpha)$.
There 
is a natural
divergence point analogue of Theorem A.
Indeed,
the following divergence point analogue of Theorem A
was first obtained by
Moran [Mo] and
Olsen \& Winter [OlWi1],
and later in a less restrictive setting by
Li, Wu \& Xiong [LiWuXi]
(see also [Ca,Vo]
for earlier
but related results in 
a
slightly different setting).

 \bigskip

\proclaim{Theorem D [LiWuXi,Mo,OlWi1]}
Let $\mu$ be the self-conformal measure defined by (2.2).
Let $C$ be a closed subset of $\Bbb R$.
If the
 OSC is satisfied, then we have 
 $$
 F_{\mu}(C)
 =
 \sup_{\alpha\in C}\beta^{*}(\alpha)\,.
 $$
\endproclaim

As a first application of 
Theorem 5.3, Corollary 5.4, Theorem 5.5 and 
Corollary 5.6 
we obtain a 
dynamical multifractal zeta-function
with an associated
 Bowen equation whose solution
equals the 
generalised 
multifractal spectrum $F_{\mu}(C)$
of a self-conformal measure $\mu$.
This is the content of the next theorem.

\bigskip

\proclaim{Theorem 6.1. 
Dynamical multifractal zeta-functinons for 
multifractal spectra of self-conformal measures}
Let
$(p_{1},\ldots,p_{N})$ be a probability vector, and let
$\mu$ denote the 
self-conformal measure associated with the list
$\big(
 \,
 V
 \,,\,
 X
 \,,\,
 (S_{i})_{i=1,\ldots,N}
  \,,\,\allowmathbreak
 (p_{i})_{i=1,\ldots,N}
 \,
 \big)$, i\.e\.
$\mu$ is the unique probability measure such that
$\mu
 =
 \sum_{i}p_{i}\mu\circ S_{i}^{-1}$.

For $C\subseteq\Bbb R$ and 
an continuous function 
$\varphi:\Sigma^{\Bbb N}\to\Bbb R$,   we
define the dynamical
self-conformal multifractal zeta-function by
  $$
  \zeta_{C}^{\sdyncon}(\varphi;z)
  =
  \sum_{n}
  \frac{z^{n}}{n}
  \left(
  \sum
  \Sb
  |\bold i|=n\\
  {}\\
  \frac{\log p_{\bold i}}{\log\diam K_{\bold i}}
 \in  C
  \endSb
  \sup_{\bold u\in[\bold i]}
  \exp
  \sum_{k=0}^{n-1}
  \varphi S^{k}\bold u
  \right)
  $$
Let $\Lambda$ be defined by (2.5) and let $\beta$ be defined by (6.1) (or, alternatively, by (6.2)).

\noindent
 {\rm (1)} Assume 
 that $C\subseteq\Bbb R$ is closed.

 \roster

  \item"(1.1)"
 There is a unique real number 
 $\,\,\scr f\,\,(C)$ such that
  $$
  \lim_{r\searrow0}
  \sigma_{\radius}
  \big(
  \,
  \zeta_{B(C,r)}^{\sdyncon}(\,\,\scr f\,\,(C)\,\Lambda;\cdot)
  \,
  \big)
  =
  1\,.
  $$
It $\alpha\in\Bbb R$ and 
$C=\{\alpha\}$, then we will write
  $\,\,\scr f\,\,(\alpha)
  =
 \,\,\scr f\,\,(C)$.

\item"(1.2)"
We have
 $$
 \scr f\,\,(C)
 =
 \sup_{\alpha\in C}\beta^{*}(\alpha)\,.
 $$
  
   \item"(1.3)"
 If the OSC is satisfied, then we have
 $$
\align
\quad
  \scr f\,\,(C)
 &=
 F_{\mu}(C)
 =
 \dim_{\Haus}
 \Bigg\{
 x\in K
 \,\Bigg|\,
 \,\underset{r\searrow 0}\to\acc
 \frac{\log\mu(B(x,r))}{\log r}
 \subseteq
 C
 \Bigg\}\,.\\
 \endalign
 $$
 In particular, 
 if the OSC is satisfied and $\alpha\in\Bbb R$,
 then we have
 $$
\align
\quad
  \scr f\,\,(\alpha)
 &=
 f_{\mu}(\alpha)
 =
 \dim_{\Haus}
 \Bigg\{
 x\in K
 \,\Bigg|\,
 \lim_{r\searrow 0}
  \frac{\log\mu(B(x,r))}{\log r}
  =
 \alpha
 \Bigg\}\,.\\
 \endalign
 $$

 \endroster
 
 \bigskip
 
 \noindent
 {\rm (2)} Assume that $C\subseteq\Bbb R$ is a closed 
interval
 with
 $\overset{\,\circ}\to{C}
  \cap
  \big(
  -\beta'(\Bbb R)
  \big)
  \not=
  \varnothing$.

\roster

 \item"(2.1)"
 There is a unique real number 
 $\scr F\,\,(C)$ such that
  $$
  \sigma_{\radius}
  \big(
  \,
  \zeta_{C}^{\sdyncon}(\scr F\,\,(C)\,\Lambda;\cdot)
  \,
  \big)
  =
  1\,.
  $$

\item"(2.2)"
We have
 $$
 \scr F\,\,(C)
 =
 \sup_{\alpha\in C}\beta^{*}(\alpha)\,.
 $$

 \item"(2.3)"
 If the OSC is satisfied then
 $$
\align
\quad
  \scr F\,\,(C)
 &=
 F_{\mu}(C)
 =
 \dim_{\Haus}
 \Bigg\{
 x\in K
 \,\Bigg|\,
 \,\underset{r\searrow 0}\to\acc
 \frac{\log\mu(B(x,r))}{\log r}
 \subseteq
 C
 \Bigg\}\,.\\
 \endalign
 $$
\endroster
 
\endproclaim

\noindent{\it Proof}\newline
This follows immediately from 
the more general Theorem 6.2 in Section 6.2 by putting $M=1$.
\hfill$\square$

\bigskip

{\bf 6.2.
Mixed multifractal 
spectra of 
self-conformal measures.}
Recently
mixed
(or simultaneous)
multifractal spectra
have generated an enormous interest in
the mathematical literature, see
[BaSa,Mo,Ol2,Ol3].
Indeed, previous
 results  (for example, (6.3) and Theorem D) only considered
the scaling behaviour  
of a single measure.
Mixed
multifractal analysis investigates
the {\it simultaneous} scaling behaviour
of finitely many 
measures.
Mixed
multifractal analysis
thus
combines local 
characteristics which depend simultaneously on various different 
aspects of the underlying dynamical system,
and provides the basis for a significantly
better understanding of 
the underlying dynamics.
We will now make these ideas precise.
For $m=1,\ldots, M$, let
$(p_{m,1},\ldots,p_{m,N})$ be a probability vector, and let
$\mu_{m}$ denote the 
self-conformal measure associated with the list
$\big(
 \,
 V
 \,,\,
 X
 \,,\,
 (S_{i})_{i=1,\ldots,N}
  \,,\,\allowmathbreak
 (p_{m,i})_{i=1,\ldots,N}
 \,
 \big)$, i\.e\.
$\mu_{m}$ is the unique probability measure such that
 $$
 \mu_{m}
 =
 \sum_{i}p_{m,i}\mu_{m}\circ S_{i}^{-1}\,.
 \tag6.4
 $$
The mixed multifractal spectrum 
$f_{\pmb\mu}$
of the list $\pmb\mu=(\mu_{1},\ldots,\mu_{M})$ is defined by
 $$
f_{\pmb\mu}(\pmb\alpha)
=
\dim_{\Haus}
 \Bigg\{
 x\in K
 \,\Bigg|\,
 \lim_{r\searrow 0}
 \Bigg(
 \frac{\log\mu_{1}(B(x,r))}{\log r}
 ,\ldots,
 \frac{\log\mu_{M}(B(x,r))}{\log r}
 \Bigg)
 =
 \pmb\alpha
 \Bigg\}
 $$
for $\pmb\alpha\in\Bbb R^{M}$.  
Of course, it is also possible to define 
generalised
mixed 
multifractal spectra
designed to
\lq\lq see"
different sets of divergence points.
Namely, we
define
the generalised  mixed
 Hausdorff multifractal spectrum
$F_{\pmb\mu}$
of the 
list $\pmb\mu=(\mu_{1},\ldots,\mu_{M})$
by
  $$
  F_{\pmb\mu}(C)
  =
\dim_{\Haus}
 \Bigg\{
 x\in K
 \,\Bigg|\,
 \,\underset{r\searrow 0}\to\acc
 \Bigg(
 \frac{\log\mu_{1}(B(x,r))}{\log r}
 ,\ldots,
 \frac{\log\mu_{M}(B(x,r))}{\log r}
 \Bigg)
 \subseteq
 C
 \Bigg\}
 $$
for  $C\subseteq\Bbb R^{M}$. 
Again we  
note that the generalised mixed multifractal spectrum
is a genuine extensions of the
traditional mixed multifractal spectrum
$F_{\mu}(\pmb\alpha)$, namely, if
$C=\{\pmb\alpha\}$ is a singleton consisting of the point $\pmb\alpha$, 
then clearly
$F_{\mu}(C)=f_{\mu}(\pmb\alpha)$.
Assuming the OSC,
 the  generalised mixed multifractal spectrum
$F_{\mu}(C)$ can be computed [Mo,Ol2].
 In order to state the result from [Mo,Ol2], 
 we introduce the following 
 definitions.
Define $\Lambda,\Phi_{m}:\Sigma^{\Bbb N}\to\Bbb R$ for $m=1,\ldots,M$ by
 $\Lambda(\bold i)=\log |DS_{i_{1}}(\pi S\bold i)|$
and
$\Phi_{m}(\bold i)=\log p_{m,i_{1}}$ 
for
$\bold i=i_{1}i_{2}\ldots\in\Sigma^{\Bbb N}$,
and
write $\pmb\Phi=(\Phi_{1},\ldots,\Phi_{M})$.
For $\bold x,\bold y\in\Bbb R^{M}$, we let
$\langle\bold x|\bold y\rangle$ denote the 
usual inner product of  $\bold x$ and $\bold y$,
and
define 
$\beta:\Bbb R^{M}\to\Bbb R$ by
 $$
 0
 =
 P\big(
 \,
 \beta(\bold q)\Lambda
 +
 \langle\bold q|\pmb\Phi\rangle
 \,
 \big)\,;
 \tag6.5
 $$
alternatively, the function $\beta:\Bbb R^{M}\to\Bbb R$
is defined by
 $$
 \sigma_{\radius}
 \big(
 \,
 \zeta^{\dyn}(\langle \bold q\,|\,\pmb\Phi\rangle+\beta(\bold q)\Lambda;\cdot)
 \,
 \big)
 =1\,.
 \tag6.6
 $$ 
Finally,
if $\varphi:\Bbb R^{M}\to\Bbb R$ is a function, we define the Legendre transform
$\varphi^{*}:\Bbb R^{M}\to[-\infty,\infty]$ of $\varphi$ by
 $$
 \varphi^{*}(\bold x)
 =
 \inf_{\bold y}(\langle\bold x|\bold y\rangle+\varphi(\bold y))\,.
 $$
The  generalised mixed multifractal spectra
$f_{\pmb\mu}$
and
$F_{\pmb\mu}$ 
are now given by the following theorem.

\bigskip

\proclaim{Theorem E [Mo,Ol2]}
Let $\mu_{1},\ldots,\mu_{M}$ be defined by (4.1)
and
let $C\subseteq\Bbb R^{M}$ be a
closed  set.
Put
$\pmb\mu
=
(\mu_{1},\ldots,\mu_{M})$.
If the OSC is satisfied,
then we have 
 $$
 F_{\pmb\mu}(C)
 =
 \sup_{\pmb\alpha\in C}\beta^{*}(\pmb\alpha)\,.
 $$
In particular,
if the OSC is satisfied
and $\pmb\alpha\in\Bbb R^{M}$, then we have
 $$
 f_{\pmb\mu}(\pmb\alpha)
 =
 \beta^{*}(\pmb\alpha)\,.
 $$
\endproclaim

\bigskip

As a second application of 
Theorem 5.3, Corollary 5.4, Theorem 5.5 and 
Corollary 5.6 
we obtain a 
dynamical multifractal zeta-function
with an associated
 Bowen equation whose solution
equals the 
generalised
mixed 
multifractal spectrum $F_{\pmb\mu}(C)$
of a list $\pmb\mu$ of self-conformal measures.
The is the content of the next theorem.

\bigskip

\proclaim{Theorem 6.2. 
Multifractal zeta-functinons for mixed
multifractal spectra of self-conformal measures}
For $m=1,\ldots, M$, let
$(p_{m,1},\ldots,p_{m,N})$ be a probability vector, and let
$\mu_{m}$ denote the 
self-conformal measure associated with the list
$\big(
 \,
 V
 \,,\,
 X
 \,,\,
 (S_{i})_{i=1,\ldots,N}
  \,,\,\allowmathbreak
 (p_{m,i})_{i=1,\ldots,N}
 \,
 \big)$, i\.e\.
$\mu_{m}$ is the unique probability measure such that
$\mu_{m}
 =
 \sum_{i}p_{m,i}\mu_{m}\circ S_{i}^{-1}$.

For $C\subseteq\Bbb R^{M}$ and 
an continuous function \
$\varphi:\Sigma^{\Bbb N}\to\Bbb R$,   we
define the dynamical
self-conformal multifractal zeta-function by
  $$
  \zeta_{C}^{\sdyncon}(\varphi;z)
  =
  \sum_{n}
  \frac{z^{n}}{n}
  \left(
  \sum
  \Sb
  |\bold i|=n\\
  {}\\
  \big(
  \frac{\log p_{1,\bold i}}{\log\diam K_{\bold i}}
  ,
  \ldots
  ,
  \frac{\log p_{M,\bold i}}{\log\diam K_{\bold i}}
  \big)
 \in  C
  \endSb
  \sup_{\bold u\in[\bold i]}
  \exp
  \sum_{k=0}^{n-1}
  \varphi S^{k}\bold u
  \right)
  $$
Let $\Lambda$ be defined by (2.5) and let $\beta$ be defined by (6.5) (or, alternatively, by (6.6)).

\noindent
 {\rm (1)} Assume 
 that $C\subseteq\Bbb R^{M}$ is closed.

 \roster

  \item"(1.1)"
 There is a unique real number 
 $\,\,\scr f\,\,(C)$ such that
  $$
  \lim_{r\searrow0}
  \sigma_{\radius}
  \big(
  \,
  \zeta_{B(C,r)}^{\sdyncon}(\,\,\scr f\,\,(C)\,\Lambda;\cdot)
  \,
  \big)
  =
  1\,.
  $$
It $\pmb\alpha\in\Bbb R^{M}$ and 
$C=\{\pmb\alpha\}$, then we will write
  $\,\,\scr f\,\,(\pmb\alpha)
  =
 \,\,\scr f\,\,(C)$.

\item"(1.2)"
We have
 $$
 \scr f\,\,(C)
 =
 \sup_{\pmb\alpha\in C}\beta^{*}(\pmb\alpha)\,.
 $$
  
   \item"(1.3)"
 If the OSC is satisfied, then we have
 $$
\align
\quad
  \scr f\,\,(C)
 &=
 F_{\pmb\mu}(C)
 =
 \dim_{\Haus}
 \Bigg\{
 x\in K
 \,\Bigg|\,
 \,\underset{r\searrow 0}\to\acc
 \Bigg(
 \frac{\log\mu_{1}(B(x,r))}{\log r}
 ,\ldots,
 \frac{\log\mu_{M}(B(x,r))}{\log r}
 \Bigg)
 \subseteq
 C
 \Bigg\}\,.\\
 \endalign
 $$
 In particular, 
 if the OSC is satisfied and $\pmb\alpha\in\Bbb R^{M}$,
 then we have
 $$
\align
\quad
  \scr f\,\,(\pmb\alpha)
 &=
 f_{\pmb\mu}(\pmb\alpha)
 =
 \dim_{\Haus}
 \Bigg\{
 x\in K
 \,\Bigg|\,
 \lim_{r\searrow 0}
 \Bigg(
 \frac{\log\mu_{1}(B(x,r))}{\log r}
 ,\ldots,
 \frac{\log\mu_{M}(B(x,r))}{\log r}
 \Bigg)
 =
 \pmb\alpha
 \Bigg\}\,.\\
 \endalign
 $$

 \endroster
 
 \bigskip
 
 \noindent
 {\rm (2)} Assume that $C\subseteq\Bbb R^{M}$ is closed 
 and convex
 with
 $\overset{\,\circ}\to{C}
  \cap
  \big(
  -\nabla \beta(\Bbb R^{M})
  \big)
  \not=
  \varnothing$.

\roster

 \item"(2.1)"
 There is a unique real number 
 $\scr F\,\,(C)$ such that
  $$
  \sigma_{\radius}
  \big(
  \,
  \zeta_{C}^{\sdyncon}(\scr F\,\,(C)\,\Lambda;\cdot)
  \,
  \big)
  =
  1\,.
  $$

\item"(2.2)"
We have
 $$
 \scr F\,\,(C)
 =
 \sup_{\pmb\alpha\in C}\beta^{*}(\pmb\alpha)\,.
 $$

 \item"(2.3)"
 If the OSC is satisfied then
 $$
\align
\quad
  \scr F\,\,(C)
 &=
 F_{\pmb\mu}(C)
=
 \dim_{\Haus}
 \Bigg\{
 x\in K
 \,\Bigg|\,
 \,\underset{r\searrow 0}\to\acc
 \Bigg(
 \frac{\log\mu_{1}(B(x,r))}{\log r}
 ,\ldots,
 \frac{\log\mu_{M}(B(x,r))}{\log r}
 \Bigg)
 \subseteq
 C
 \Bigg\}\,.\\
 \endalign
 $$
\endroster
 
\endproclaim

\bigskip

We  will now prove Theorem 6.2.
Recall  that the function $\Lambda:\Sigma^{\Bbb N}\to\Bbb R$ is defined by
$\Lambda(\bold i)=\log|DS_{i_{1}}(\pi S\bold i)|$ for $\bold i=i_{1}i_{2}\ldots\in\Sigma^{\Bbb N}$.
Also,
recall that
$\pmb\Phi=(\Phi_{1},\ldots,\Phi_{M})$
where
$\Phi_{m}:\Sigma^{\Bbb N}\to\Bbb R$ is defined by
$\Phi_{m}(\bold i)=\log p_{m,i_{1}}$ 
for
$\bold i=i_{1}i_{2}\ldots\in\Sigma^{\Bbb N}$.
We now introduce the following definitions.
For $\mu\in\Cal P(\Sigma^{\Bbb N})$, write
$\int\pmb\Phi\,d\mu
=
(\int\Phi_{1}\,d\mu,\ldots,\int\Phi_{M}\,d\mu)$,
and
define
$\Gamma:\Cal P(\Sigma^{\Bbb N})\to\Bbb R^{M}$
and
$\Delta:\Cal P(\Sigma^{\Bbb N})\to\Bbb R$ by
 $$
 \align
 \Gamma(\mu)
&=
\int\pmb\Phi\,d\mu\,,\,\,\,\,
\Delta(\mu)
 =
\int\Lambda\,d\mu\,.
\endalign
 $$
 Observe that the maps $\Gamma$ and $\Delta$ are 
 affine and continuous.
Finally, define 
$U:\Cal P(\Sigma^{\Bbb N})\to\Bbb R^{M}$
by $U=\frac{\Gamma}{\Delta}$, i\.e\. 
 $$
 U\mu
 =
 \frac{\Gamma(\mu)}{\Delta(\mu)}
 =
 \frac{\int\pmb\Phi\,d\mu}{\int\Lambda\,d\mu}\,,
 \tag6.7
 $$
and note that if $\bold i\in\Sigma^{*}$, then
 $$
 UL_{|\bold i|}[\bold i]
 =
 \Bigg\{
  \Bigg(
 \,
 \frac{\log p_{1,\bold i}}{\log |DS_{\bold i}(\pi\bold u)|}
 \,
 ,
 \ldots
 ,
 \,
 \frac{\log p_{M,\bold i}}{\log |DS_{\bold i}(\pi\bold u)|} 
 \,
 \Bigg)
 \,\Bigg|\,
 \bold u\in\Sigma^{\Bbb N}
 \Bigg\}\,.
 $$
 It therefore follows that
  $$
  \zeta_{C}^{\dyn,U}(\varphi;z)
  =
  \sum_{n}
  \,\,
  \frac{z^{n}}{n}
  \,\,
  \left(
  \sum
  \Sb
   |\bold i|=n\\
   {}\\
   \forall
 \bold u\in\Sigma^{\Bbb N}
 \,\,:\,\,
  \big(
 \,
 \frac{\log p_{1,\bold i}}{\log |DS_{\bold i}(\pi\bold u)|}
 \,
 ,
 \ldots
 ,
 \,
 \frac{\log p_{M,\bold i}}{\log |DS_{\bold i}(\pi\bold u)|} 
 \,
 \big)
\in
C
  \endSb
  \sup_{\bold u\in[\bold i]}
  \,\,
  \exp
  \sum_{k=0}^{n-1}
  \varphi S^{k}\bold u
  \right)
  \,.
  \tag6.8
  $$ 
In order to prove Theorem 6.2, we first prove that 
radii of convergence of 
the zeta-functions
$\zeta_{C}^{\dyn,U}(\varphi;\cdot)$
and
$\zeta_{C}^{\dyn\text{-}\scon}(\varphi;\cdot)$
are comparable;
this is the context of Proposition 6.5.
However, in order to prove Proposition 6.5 we first
prove two small auxiliary
results, namely, Proposition 6.3 and Proposition 6.4.

\bigskip

\proclaim{Proposition 6.3}
Let $U$
be defined by (6.7).
Let $\Lambda$ be defined by (2.5)
and let
 $\beta$
 be defined by (6.5).
Let $C\subseteq \Bbb R^{M}$ be a closed set
and
let $t$ be the unique real number such that
 $$
 \sigma_{\radius}
 \big(
 \,
 \zeta_{C}^{\dyn,U}(t\Lambda;\cdot)
 \,
 \big)
 =
 1\,.
 $$
Then 
 $$
 \sup
 \Sb
 \mu\in\Cal P_{S}(\Sigma^{\Bbb N})\\
 {}\\
 U\mu\in C
 \endSb
 \Bigg(
 h(\mu)+t\int\Lambda\,d\mu
 \Bigg)
=
0\,,
$$
and we have
 $$
 t
 =
 \sup
 \Sb
 \mu\in\Cal P_{S}(\Sigma^{\Bbb N})\\
 {}\\
 U\mu\in C
 \endSb
 -
 \frac{h(\mu)}{\int\Lambda\,d\mu}
 =
 \sup_{\pmb\alpha\in C}
 \beta^{*}(\pmb\alpha)\,.
$$
\endproclaim
\noindent{\it Proof}\newline

\noindent
We first note that
it follows immediately from
Theorem 5.5 that
 $$
 \sup
 \Sb
 \mu\in\Cal P_{S}(\Sigma^{\Bbb N})\\
 {}\\
 U\mu\in C
 \endSb
 \Bigg(
 h(\mu)+t\int\Lambda\,d\mu
 \Bigg)
 =
 \overline P_{C}^{U}(t\Lambda) 
 =
 -\log
  \sigma_{\radius}
 \big(
 \,
 \zeta_{C}^{\dyn,U}(t\Lambda;\cdot)
 \,
 \big)
 =
 0\,.
 $$
It therefore suffices to 
prove the following three inequalities, namely
 $$
 \align
  \sup_{\pmb\alpha\in C}
 \beta^{*}(\pmb\alpha)
&\le
 \sup
 \Sb
 \mu\in\Cal P_{S}(\Sigma^{\Bbb N})\\
 {}\\
 U\mu\in C
 \endSb
 -
 \frac{h(\mu)}{\int\Lambda\,d\mu}\,,
 \tag6.9\\
 \sup
 \Sb
 \mu\in\Cal P_{S}(\Sigma^{\Bbb N})\\
 {}\\
 U\mu\in C
 \endSb
 -
 \frac{h(\mu)}{\int\Lambda\,d\mu}
&\le
  t\,,
  \tag6.10\\
  t
  &\le
 \sup_{\pmb\alpha\in C}
 \beta^{*}(\pmb\alpha)\,.
 \tag6.11
 \endalign
 $$

{\it Proof of (6.9).}
For $s\in\Bbb R$ and $\bold q\in\Bbb R^{M}$,
let $\mu_{s,\bold q}$ denote the
Gibbs state of $s\Lambda+\langle\bold q|\pmb\Phi\rangle$.
We now prove the following three claims.

\medskip

{\it Claim 1. For all $\bold q$, we have
$\frac
 {\int\pmb\Phi\,d\mu_{\beta(\bold q),\bold q}}
 {\int\Lambda\,d\mu_{\beta(\bold q),\bold q}} 
 =
 -\nabla\beta(\bold q)$.}

\noindent
{\it Proof of Claim 1.}
Define $F:\Bbb R\times\Bbb R^{M}\to\Bbb R$ by
$F(s,\bold q)
 =
 P
 \big(
 \,
 s\Lambda+\langle\bold q|\pmb\Phi\rangle
 \,
 \big)$
for $s\in\Bbb R$ and $\bold q\in\Bbb R^{M}$.
It follows from [Rue1] that
$F$ is real analytic with
 $$
 \align
 \nabla F(s,\bold q)
&=
\Bigg(
\,
 \int\Lambda\,d\mu_{s,\bold q}
 \,,\,
\int\pmb\Phi\,d\mu_{s,\bold q} 
\,
\Bigg)\,.
\tag6.12
 \endalign
 $$
Next, 
since
$0
 =
 F(\, \beta(\bold q) \,,\,\bold q  \,)$
for all $\bold q$, 
it follows from (6.12)
and
an application of the chain rule that
$\bold 0
=
\int\Lambda\,d\mu_{\beta(\bold q),\bold q}
\,\,
\nabla\beta(\bold q)
\,\,
+
\,\,
\int\pmb\Phi\,d\mu_{\beta(\bold q),\bold q}$
for all $\bold q$.
This clearly implies that
 $-\nabla\beta(\bold q)
 =
 \frac
 {\int\pmb\Phi\,d\mu_{\beta(\bold q),\bold q}}
 {\int\Lambda\,d\mu_{\beta(\bold q),\bold q}} $
for all $\bold q$.
This completes the proof of Claim 1.

\medskip

{\it Claim 2. For all $\bold q$, we have
$ -
 \frac{h(\mu_{\beta(\bold q),\bold q})}{\int\Lambda\,d\mu_{\beta(\bold q),\bold q}}
 \ge
 \beta^{*}(-\nabla\beta(\bold q))$.}

\noindent
{\it Proof of Claim 2.}
Since $\mu_{\beta(\bold q),\bold q}$ is a Gibbs state
of 
$\beta(\bold q)\Lambda+\langle\bold q|\pmb\Phi\rangle$
and
$P(\, \beta(\bold q)\Lambda+\langle\bold q|\pmb\Phi\rangle \,)=0$,
we deduce that
$0
=
P(\,\beta(\bold q)\Lambda+\langle\bold q|\pmb\Phi\rangle\,)
=
h(\mu_{\beta(\bold q),\bold q})
+
\int(\beta(\bold q)\Lambda+\langle\bold q|\pmb\Phi\rangle)\,d\mu_{\beta(\bold q),\bold q}
=
h(\mu_{\beta(\bold q),\bold q})
+
\beta(\bold q)\int\Lambda\,d\mu_{\beta(\bold q),\bold q}
+
\langle\bold q|\int\pmb\Phi\,d\mu_{\beta(\bold q),\bold q}\rangle$.
Hence
 $$
 \align
 -
 \frac{h(\mu_{\beta(\bold q),\bold q})}{\int\Lambda\,d\mu_{\beta(\bold q),\bold q}}
&=
 \beta(\bold q)
 +
 \frac
 {\langle\bold q|\int\pmb\Phi\,d\mu_{\beta(\bold q),\bold q}\rangle}
 {\int\Lambda\,d\mu_{\beta(\bold q),\bold q}}
 =
 \beta(\bold q)
 +
 \Bigg\langle
 \bold q
 \Bigg|
 {\dsize \frac{\int\pmb\Phi\,d\mu_{\beta(\bold q),\bold q}}{\int\Lambda\,d\mu_{\beta(\bold q),\bold q}}}
 \Bigg\rangle\,.
 \tag6.13
 \endalign
 $$
Combining Claim 1 and (6.13) now yields
 $$
 \align
 -
 \frac{h(\mu_{\beta(\bold q),\bold q})}{\int\Lambda\,d\mu_{\beta(\bold q),\bold q}}
&=
 \beta(\bold q)
  +
 \langle
 \bold q
 |
 -\nabla\beta(\bold q)
 \rangle
 \ge
 \inf_{\bold r} 
 \big(
 \,
  \beta(\bold r)
  +
 \langle
 \bold r
 |
 -\nabla\beta(\bold q)
 \rangle
 \,
 \big)
 =
 \beta^{*}(-\nabla\beta(\bold q)) 
 \endalign
 $$
for all $\bold q$.
This completes the proof of Claim 2.
 
 \medskip

{\it Claim 3. For all $\pmb\alpha\in\Bbb R^{M}$, we have
 $\beta^{*}(\pmb\alpha)
 \le
 \sup
 \Sb
 \mu\in\Cal P_{S}(\Sigma^{\Bbb N})\\
 U\mu=\pmb\alpha
 \endSb
 -
 \frac{h(\mu)}{\int\Lambda\,d\mu}$.}

\noindent
{\it Proof of Claim 3.}
If $ \beta^{*}(\pmb\alpha)=-\infty$, then
the statement is clear. Hence, we may assume that
$ \beta^{*}(\pmb\alpha)>-\infty$.
In this case it follows from the convexity of $\beta$
that there is a point $\bold q_{\pmb\alpha}\in\Bbb R^{M}$
such that
$\pmb\alpha=-\nabla\beta(\bold q_{\pmb\alpha})$, see [Ro].
It therefore follows from Claim 1 that the measure 
$\mu_{\beta(\bold q_{\pmb\alpha}),\bold q_{\pmb\alpha}}$
satisfies
$U\mu_{\beta(\bold q_{\pmb\alpha}),\bold q_{\pmb\alpha}}
=
\frac
 {\int\pmb\Phi\,d\mu_{\beta(\bold q_{\pmb\alpha}),\bold q_{\pmb\alpha}}}
 {\int\Lambda\,d\mu_{\beta(\bold q_{\pmb\alpha}),\bold q_{\pmb\alpha}}} 
 =
 -\nabla\beta(\bold q_{\pmb\alpha})
=
\pmb\alpha$, whence, using Claim 2,
 $$
 \sup
 \Sb
 \mu\in\Cal P_{S}(\Sigma^{\Bbb N})\\
 {}\\
 U\mu=\pmb\alpha
 \endSb
 -
 \frac{h(\mu)}{\int\Lambda\,d\mu}
 \ge
 -
 \frac{h(\mu_{\beta(\bold q_{\pmb\alpha}),\bold q_{\pmb\alpha}})}
 {\int\Lambda\,d\mu_{\beta(\bold q_{\pmb\alpha}),\bold q_{\pmb\alpha}}}
 \ge
 \beta^{*}(-\nabla\beta(\bold q_{\pmb\alpha}))
 =
 \beta^{*}(\pmb\alpha)\,.
$$
This completes the proof of Claim 3.

\medskip

We can now prove the required inequality.
Indeed, it follows immediately from Claim 3 that
 $$
 \align
 \sup_{\pmb\alpha\in C}
 \beta^{*}(\pmb\alpha)
 &\le
 \sup_{\pmb\alpha\in C}
 \,\,
 \sup
 \Sb
 \mu\in\Cal P_{S}(\Sigma^{\Bbb N})\\
 {}\\
 U\mu=\pmb\alpha
 \endSb
 -
 \frac{h(\mu)}{\int\Lambda\,d\mu}
 =
  \sup
 \Sb
 \mu\in\Cal P_{S}(\Sigma^{\Bbb N})\\
 {}\\
 U\mu\in C
 \endSb
 -
 \frac{h(\mu)}{\int\Lambda\,d\mu}\,.
 \endalign
 $$
This completes the proof of (6.9).

{\it Proof of (6.10).}
Fix $\mu\in \Cal P_{S}(\Sigma^{\Bbb N})$ with $U\mu\in C$.
It follows from the defnition of $t$ that
$h(\mu)+t\int\Lambda\,d\mu\le0$, and since $\Lambda<0$, 
we therefore conclude that
$-\frac{h(\mu)}{\int\Lambda\,d\mu}\le t$.
Taking supremum over all 
$\mu\in \Cal P_{S}(\Sigma^{\Bbb N})$ with $U\mu\in C$
in this inequality
now gives the desired result.
This completes the proof of (6.10).

{\it Proof of (6.11).}
Fix $\mu\in\Cal P_{S}(\Sigma^{\Bbb N})$ with
$U\mu\in C$.
Next, let $\bold q\in\Bbb R^{M}$.
It now follows from the definition of $t$ and $\beta(\bold q)$ that
We now have
 $$
 \align
\sup
 \Sb
 \nu\in\Cal P_{S}(\Sigma^{\Bbb N})\\
 {}\\
 U\nu\in C
 \endSb
\Bigg(
h(\nu)+t\int\Lambda\,d\nu
\Bigg)
&=
 0
 =
 P(\,\beta(\bold q)\Lambda+\langle\bold q|\pmb\Phi\rangle\,)\,.
 \tag6.14
\endalign
$$
 Also, using the variational principle (see [Wa]) we  conclude  that
$
P(\,\beta(\bold q)\Lambda+\langle\bold q|\pmb\Phi\rangle\,)
=
\sup_{\nu\in\Cal P(\Sigma^{\Bbb N})}
(\,h(\nu)+\int(\beta(\bold q)\Lambda+\langle\bold q|\pmb\Phi\rangle)\,d\nu\,)
\ge
h(\mu)+\int(\beta(\bold q)\Lambda+\langle\bold q|\pmb\Phi\rangle)\,d\mu
=
h(\mu)+\beta(\bold q)\int\Lambda\,d\mu+\langle\bold q|\int\pmb\Phi\,d\mu\rangle$.
We deduce from this and (6.14) that
$$
 \align
\sup
 \Sb
 \nu\in\Cal P_{S}(\Sigma^{\Bbb N})\\
 {}\\
 U\nu\in C
 \endSb
\Bigg(
h(\nu)+t\int\Lambda\,d\nu
\Bigg)
&\ge
 h(\mu)+\beta(\bold q)\int\Lambda\,d\mu
 +
 \Bigg\langle
 \bold q\Bigg|\int\pmb\Phi\,d\mu
 \Bigg\rangle\,.
 \tag6.15
\endalign
$$
Next, observe that
$U\mu=\frac{\int\pmb\Phi\,d\mu}{\int\Lambda\,d\mu}$, whence
$\int\pmb\Phi\,d\mu=\int\Lambda\,d\mu\,U\mu$, and it therefore follows from (6.15) that
 $$
 \align
\sup
 \Sb
 \nu\in\Cal P_{S}(\Sigma^{\Bbb N})\\
 {}\\
 U\nu\in C
 \endSb
\Bigg(
h(\nu)+t\int\Lambda\,d\nu
\Bigg)
&\ge
 h(\mu)+\beta(\bold q)\int\Lambda\,d\mu
 +
 \Bigg\langle
 \bold q\Bigg|\int\Lambda\,d\mu\,U\mu
 \Bigg\rangle\\
&=
 h(\mu)
 +
 \big(
 \beta(\bold q) 
 +
 \langle\bold q|U\mu
 \rangle
 \big)
 \int\Lambda\,d\mu\,.
 \tag6.16
\endalign
$$
Taking supremum over all $\bold q$  in (6.16)
and using the fact that $\Lambda<0$ now gives
  $$
 \align
\sup
 \Sb
 \nu\in\Cal P_{S}(\Sigma^{\Bbb N})\\
 {}\\
 U\nu\in C
 \endSb
\Bigg(
h(\nu)+t\int\Lambda\,d\nu
\Bigg)
&\ge
 \sup_{\bold q}
 \Bigg(
 h(\mu)
 +
 \big(
 \beta(\bold q) 
 +
 \langle\bold q|U\mu
 \rangle
 \big)
 \int\Lambda\,d\mu
 \Bigg)\\
&=
 h(\mu)
 +
 \inf_{\bold q}
 \big(
 \beta(\bold q) 
 +
 \langle\bold q|U\mu
 \rangle
 \big)
 \int\Lambda\,d\mu\\
&=
 h(\mu)
 +
 \beta^{*}(U\mu)
 \int\Lambda\,d\mu\,.
 \tag6.17
\endalign
$$
By assumption
$U\mu\in C$, whence
$\beta^{*}(U\mu)\le\sup_{\pmb\alpha\in C}\beta^{*}(\pmb\alpha)$. 
It follows from this 
and 
the inequality
$\Lambda<0$ that
$\beta^{*}(U\mu)\int\Lambda\,d\mu
\ge
(\sup_{\pmb\alpha\in C}\beta^{*}(\pmb\alpha))\int\Lambda\,d\mu$,
and we therefore conclude from (6.17) that
 $$
 \align
\sup
 \Sb
 \nu\in\Cal P_{S}(\Sigma^{\Bbb N})\\
 {}\\
 U\nu\in C
 \endSb
\Bigg(
h(\nu)+t\int\Lambda\,d\nu
\Bigg)
&\ge
  h(\mu)
 +
\bigg(\sup_{\pmb\alpha\in C}\beta^{*}(\pmb\alpha)\bigg)
\int\Lambda\,d\mu\,.
\tag6.18
\endalign
$$
Finally, taking supremum over all $\mu\in\Cal P_{S}(\Sigma^{\Bbb N})$
with $U\mu\in C$
in (6.18) yields
 $$
 \align
\sup
 \Sb
 \nu\in\Cal P_{S}(\Sigma^{\Bbb N})\\
 {}\\
 U\nu\in C
 \endSb
\Bigg(
h(\nu)+t\int\Lambda\,d\nu
\Bigg)
&\ge
\sup
 \Sb
 \mu\in\Cal P_{S}(\Sigma^{\Bbb N})\\
 {}\\
 U\mu\in C
 \endSb
\Bigg(
  h(\mu)
 +
 \bigg(\sup_{\pmb\alpha\in C}\beta^{*}(\pmb\alpha)\bigg)
 \int\Lambda\,d\mu
 \Bigg)\,.
 \tag6.19
\endalign
$$
Since $\Lambda<0$, we now deduce from inequality (6.19) that
$t\le \sup_{\pmb\alpha\in C}\beta^{*}(\pmb\alpha)$.
This completes the proof of (6.11).
\hfill$\square$

\bigskip

\proclaim{Proposition 6.4}
Let $U$
be defined by (6.7).
Let $\Lambda$ be defined by (2.5)
and let
 $\beta$
 be defined by (6.5).
Then
$-\nabla\beta(\Bbb R^{M})\subseteq U(\,\Cal P_{S}(\Sigma^{\Bbb N})\,)$.
\endproclaim
\noindent{\it Proof}\newline
\noindent
This follows from Claim 1 in the proof of Proposition 6.3.
\hfill$\square$

\bigskip

\proclaim{Proposition 6.5}
Let $U$
be defined by (6.7).
Fix a continuous function
$\varphi:\Sigma^{\Bbb N}\to\Bbb R$.
\roster
\item"(1)"
There is a sequence $(\Delta_{n})_{n}$ with $\Delta_{n}>0$ 
and 
$\Delta_{n}\to 0$ such that
for all
closed subsets $W$ of $\Bbb R^{M}$
and
for
all
$n\in\Bbb N$, $\bold i\in\Sigma^{n}$ and $\bold u\in\Sigma^{\Bbb N}$, we have
 $$
 \align
\quad\quad\quad\,\,\,
 \dist
 \Bigg(
 \,
  \Bigg(
 \,
&\frac{\log p_{1,\bold i}}{\log |DS_{\bold i}(\pi\bold u)|}
 \,
,
 \ldots
 ,
 \,
 \frac{\log p_{M,\bold i}}{\log |DS_{\bold i}(\pi\bold u)|} 
 \,
 \Bigg)
 \,,\,
 W
 \,
 \Bigg)\\
&\le
\dist
\Bigg(
 \,
  \Bigg(
 \,
 \frac{\log p_{1,\bold i}}{\log \diam K_{\bold i}}
 \,
,
 \ldots
 ,
 \,
 \frac{\log p_{M,\bold i}}{\log \diam K_{\bold i}} 
 \,
 \Bigg)
 \,,\,
 W
 \,
 \Bigg)
\,\,
\,\,\,\,\,\,\,\,\,\,
+
\,\,
\Delta_{n}\,,
\tag6.20\\
&{}\\
\dist
\Bigg(
 \,
  \Bigg(
 \,
&\frac{\log p_{1,\bold i}}{\log \diam K_{\bold i}}
 \,
,
 \ldots
 ,
 \,
 \frac{\log p_{M,\bold i}}{\log \diam K_{\bold i}} 
 \,
 \Bigg)
 \,,\,
 W
 \,
 \Bigg)
\\
&\le
\dist
 \Bigg(
 \,
  \Bigg(
 \,
 \frac{\log p_{1,\bold i}}{\log |DS_{\bold i}(\pi\bold u)|}
 \,
,
 \ldots
 ,
 \,
 \frac{\log p_{M,\bold i}}{\log |DS_{\bold i}(\pi\bold u)|} 
 \,
 \Bigg)
 \,,\,
 W
 \,
 \Bigg)
 \,\,
+
\,\,
\Delta_{n}\,.
\tag6.21
\endalign
$$

\item"(2)"
Let $W$ be a closed subset 
of $\Bbb R^{M}$.
For all $r>0$, we have
 $$
 \align
 \sigma_{\radius}
 \big(
 \,
 \zeta_{B(W,r)}^{\dyn,U}(\varphi;\cdot)
 \,
 \big)
&\le
 \sigma_{\radius}
 \big(
 \,
 \zeta_{W}^{\dyn\text{-}\scon}(\varphi;\cdot)
 \,
 \big)\,,
 \tag6.22\\
 \sigma_{\radius}
 \big(
 \,
 \zeta_{B(W,r)}^{\dyn\text{-}\scon}(\varphi;\cdot)
 \,
 \big)
&\le
\sigma_{\radius}
 \big(
 \,
 \zeta_{W}^{\dyn,U}(\varphi;\cdot)
 \,
 \big)\,.
 \tag6.23
 \endalign
 $$

\item"(3)"
Assume that $C\subseteq\Bbb R^{M}$ is closed.
Then we have
 $$
 \lim_{r\searrow0}
 \sigma_{\radius}
 \big(
 \,
 \zeta_{B(C,r)}^{\dyn\text{-}\scon}(\varphi;\cdot) 
 \,
 \big)
 =
 \lim_{r\searrow0}
 \sigma_{\radius}
 \big(
 \,
 \zeta_{B(C,r)}^{\dyn,U}(\varphi;\cdot)
 \,
 \big)\,.
$$

\item"(4)"
Assume that $C\subseteq\Bbb R^{M}$ is closed
and
 convex 
 with
$\overset{\circ}\to{C}\cap \,\big( -\nabla\beta(\Bbb R^{M})\big)\not=\varnothing$.
Then we have
 $$
 \sigma_{\radius}
 \big(
 \,
 \zeta_{C}^{\dyn\text{-}\scon}(\varphi;\cdot) 
 \,
 \big)
 =
 \sigma_{\radius}
 \big(
 \,
 \zeta_{C}^{\dyn,U}(\varphi;\cdot)
 \,
 \big)\,.
$$

\endroster
\endproclaim
\noindent{\it Proof}\newline
\noindent
(1)
It is well-known and follows from the 
Principle of Bounded Distortion 
(see, for example, [Bar,Fa2])
that there is a constant $c>0$ such that
for all integers $n$
and all $\bold i$ with $|\bold i|=n$
and all $\bold u,\bold v\in[\bold i]$,
 we have
$\frac{1}{c}
\le
\frac{|DS_{\bold i}(\pi S^{n}\bold u)|}{\diam K_{\bold i}}
\le
c$
and
$\frac{1}{c}
\le
\frac{|DS_{\bold i}(\pi S^{n}\bold u)|}{|DS_{\bold i}(\pi S^{n}\bold v)|}
\le
c$.
It is not difficult to see that the desired result follows from this.

\noindent
(2)
Fix $r>0$.
Let $(\Delta_{n})_{n}$ be the sequence from (1).
Since $\Delta_{n}\to 0$, we can find
a positive integer
$N_{r}$
such that if $n\ge N_{r}$, then $\Delta_{n}< r$.
Consequently, using (6.21) in Part (1),
for all $n\ge N_{r}$, we have
$$
\align
\sum
 \Sb
 |\bold i|=n\\
 {}\\
UL_{n}[\bold i]\subseteq W
\endSb
  \sup_{\bold u\in[\bold i]}
  \exp
  \sum_{k=0}^{n-1}
  \varphi S^{k}\bold u
 & 
 =
 \qquad
 \sum
 \Sb
 |\bold i|=n\\
 {}\\
 \forall
 \bold u\in\Sigma^{\Bbb N}
 \,\,:\,\,
  \big(
 \,
 \frac{\log p_{1,\bold i}}{\log |DS_{\bold i}(\pi\bold u)|}
 \,
 ,
 \ldots
 ,
 \,
 \frac{\log p_{M,\bold i}}{\log |DS_{\bold i}(\pi\bold u)|} 
 \,
 \big)
\in
W
 \endSb
 \qquad
  \sup_{\bold u\in[\bold i]}
  \exp
  \sum_{k=0}^{n-1}
  \varphi S^{k}\bold u\\
 & 
   =
   \sum
 \Sb
 |\bold i|=n\\
 {}\\
 \forall
 \bold u\in\Sigma^{\Bbb N}
 \,\,:\,\,
 \dist
 \big(
 \,
  \big(
 \,
 \frac{\log p_{1,\bold i}}{\log |DS_{\bold i}(\pi\bold u)|}
 \,
 ,
 \ldots
 ,
 \,
 \frac{\log p_{M,\bold i}}{\log |DS_{\bold i}(\pi\bold u)|} 
 \,
 \big)
 \,,\,
 W
 \,
 \big)
 \,
 =
 0
 \endSb
 \,
  \sup_{\bold u\in[\bold i]}
  \exp
  \sum_{k=0}^{n-1}
  \varphi S^{k}\bold u\\
   &
   \le
   \quad\,\,
   \sum
  \Sb
 |\bold i|=n\\
 {}\\
 \dist
 \big(
 \,
  \big(
 \,
 \frac{\log p_{1,\bold i}}{\log \diam K_{i}}
 \,
 ,
 \ldots
 ,
 \,
 \frac{\log p_{M,\bold i}}{\log \diam K_{i}} 
 \,
 \big)
 \,,\,
 W
 \,
 \big)
 \,
 \le
 \,
 0
 +
 \Delta_{r}
 \endSb
 \quad\,\,\,
  \sup_{\bold u\in[\bold i]}
  \exp
  \sum_{k=0}^{n-1}
  \varphi S^{k}\bold u\\
   &
      \le
   \qquad\,
   \sum
  \Sb
 |\bold i|=n\\
 {}\\
 \dist
 \big(
 \,
  \big(
 \,
 \frac{\log p_{1,\bold i}}{\log \diam K_{i}}
 \,
 ,
 \ldots
 ,
 \,
 \frac{\log p_{M,\bold i}}{\log \diam K_{i}} 
 \,
 \big)
 \,,\,
 W
 \,
 \big)
 \,
 <
 \,
 r 
 \endSb
 \qquad\,\,
  \sup_{\bold u\in[\bold i]}
  \exp
  \sum_{k=0}^{n-1}
  \varphi S^{k}\bold u\\
   & 
   =
   \qquad\quad\,
   \sum
  \Sb
 |\bold i|=n\\
 {}\\
  \big(
 \,
 \frac{\log p_{1,\bold i}}{\log \diam K_{i}}
 \,
 ,
 \ldots
 ,
 \,
 \frac{\log p_{M,\bold i}}{\log \diam K_{i}} 
 \,
 \big)
\in
B(W,r)
 \endSb
 \qquad\quad\,\,\,
  \sup_{\bold u\in[\bold i]}
  \exp
  \sum_{k=0}^{n-1}
  \varphi S^{k}\bold u\,.  \\
  \tag6.24
  \endalign
  $$
A similar argument using (6.20) in Part 1 shows that
 $$
 \sum
  \Sb
 |\bold i|=n\\
 {}\\
  \big(
 \,
 \frac{\log p_{1,\bold i}}{\log \diam K_{i}}
 \,
 ,
 \ldots
 ,
 \,
 \frac{\log p_{M,\bold i}}{\log \diam K_{i}} 
 \,
 \big)
\in
B(W,r)
 \endSb
  \sup_{\bold u\in[\bold i]}
  \exp
  \sum_{k=0}^{n-1}
  \varphi S^{k}\bold u
  \le
  \sum
 \Sb
 |\bold i|=n\\
 {}\\
UL_{n}[\bold i]\subseteq W
\endSb
  \sup_{\bold u\in[\bold i]}
  \exp
  \sum_{k=0}^{n-1}
  \varphi S^{k}\bold u\,.
  \tag6.25
   $$
The desired results follow immediately from inequalities
(6.24) and (6.25).

\noindent
(3)
This result follows easily from Part (2).

\noindent
(4)
\lq\lq$\ge$"
It follows from (6.22), Theorem 5.3 and Theorem 5.5  that
 $$
 \align
 -
 \log
 \sigma_{\radius}
 \big(
 \,
 \zeta_{C}^{\dyn\text{-}\scon}(\varphi;\cdot)
 \,
 \big)
&\le
\liminf_{r\searrow 0}
 -
 \log
 \sigma_{\radius}
 \big(
 \,
 \zeta_{B(C,r)}^{\dyn,U}(\varphi;\cdot)
 \,
 \big)
 \qquad
 \text{[by (6.22)]}\\ 
&=
\sup
\Sb
\mu\in\Cal P_{S}(\Sigma^{\Bbb N})\\
{}\\
U\mu\in C
\endSb
\Bigg(
h(\mu)+\int\varphi\,d\mu
\Bigg)
  \qquad\quad\,\,
 \text{[by Theorem 5.3]}\\
&=
  -
 \log
 \sigma_{\radius}
 \big(
 \,
 \zeta_{C}^{\dyn,U}(\varphi;\cdot)
 \,
 \big)\,.
   \qquad\qquad
   \,\,\,\,\,
 \text{[by Theorem 5.5]}
  \endalign
 $$
It follows from this inequality that
$\sigma_{\radius}
 \big(
 \,
 \zeta_{C}^{\dyn\text{-}\scon}(\varphi;\cdot)
 \,
 \big)
\ge
 \sigma_{\radius}
 \big(
 \,
 \zeta_{C}^{\dyn,U}(\varphi;\cdot)
 \,
 \big)$.

\noindent
\lq\lq$\le$"
For $\varepsilon>0$, write
$I(C,\varepsilon)
=
 \{
 x\in C
 \,|\,
 \dist(x,\partial C)\ge \varepsilon
 \}$.

Next, fix
$\varepsilon>0$
and note that
if
$r>0$ with $2r<\varepsilon$, 
then it follows from 
(6.23) applied to $W=B(I(C,\varepsilon),r)$
that
 $$
  -
 \log
 \sigma_{\radius}
 \big(
 \,
 \zeta_{B(\,B(I(C,\varepsilon),r)\,,\,r)}^{\dyn\text{-}\scon}(\varphi;\cdot)
 \,
 \big)
 \ge 
 -
 \log
\sigma_{\radius}
 \big(
 \,
 \zeta_{B(I(C,\varepsilon),r)}^{\dyn,U}(\varphi;\cdot)
 \,
 \big)
 \,.
 \tag6.26
$$
However, for  $r>0$
 with $2r<\varepsilon$
 it follows from the convexity of $C$ that
$B(\,B(I(C,\varepsilon),r)\,,\,r)
\subseteq
B(I(C,\varepsilon),2r)
\subseteq 
C$, 
whence
$  \sigma_{\radius}
 \big(
 \,
 \zeta_{C}^{\dyn\text{-}\scon}(\varphi;\cdot)
 \,
 \big)
 \le
 \sigma_{\radius}
 \big(
 \,
 \zeta_{B(\,B(I(C,\varepsilon),r)\,,\,r)}^{\dyn\text{-}\scon}(\varphi;\cdot)
 \,
 \big)$,
and so
 $-
\log
 \sigma_{\radius}
 \big(
 \,
 \zeta_{C}^{\dyn\text{-}\scon}(\varphi;\cdot)
 \,
 \big)
 \ge
 -
 \log
  \sigma_{\radius}
 \big(
 \,
 \zeta_{B(\,B(I(C,\varepsilon),r)\,,\,r)}^{\dyn\text{-}\scon}(\varphi;\cdot)
 \,
 \big)$.
We conclude from this
and (6.26) that
if $r>0$ with $2r<\varepsilon$, then
 $$
 -
 \log
 \sigma_{\radius}
 \big(
 \,
 \zeta_{C}^{\dyn\text{-}\scon}(\varphi;\cdot)
 \,
 \big)
 \ge
 -
 \log
\sigma_{\radius}
 \big(
 \,
 \zeta_{B(I(C,\varepsilon),r)}^{\dyn,U}(\varphi;\cdot)
 \,
 \big)\,.
 \tag6.27
$$
Next, 
since 
$I(C,\varepsilon)$ is closed, 
it follows from (6.27) and  Theorem 5.3
that 
if $\varepsilon>0$,
then
 $$
 \align
 -
 \log
 \sigma_{\radius}
 \big(
 \,
 \zeta_{C}^{\dyn\text{-}\scon}(\varphi;\cdot)
 \,
 \big)
&\ge
 \limsup_{r\searrow 0}
 -
 \log
\sigma_{\radius}
 \big(
 \,
 \zeta_{B(I(C,\varepsilon),r)}^{\dyn,U}(\varphi;\cdot)
 \,
 \big)\\
 &=
\sup
\Sb
\mu\in\Cal P_{S}(\Sigma^{\Bbb N})\\
{}\\
U\mu\in I(C,\varepsilon)
\endSb
\Bigg(
h(\mu)+\int\varphi\,d\mu
\Bigg)\,.
\tag6.27
 \endalign
 $$
Taking supremum over all $\varepsilon>0$ in (6.27) gives 
 $$
 \align
  -
 \log
 \sigma_{\radius}
 \big(
 \,
 \zeta_{C}^{\dyn\text{-}\scon}(\varphi;\cdot)
 \,
 \big)
&\ge
\sup_{\varepsilon>0}
\,\,
\sup
\Sb
\mu\in\Cal P_{S}(\Sigma^{\Bbb N})\\
{}\\
U\mu\in I(C,\varepsilon)
\endSb
\Bigg(
h(\mu)+\int\varphi\,d\mu
\Bigg)\\
&=
\sup
\Sb
\mu\in\Cal P_{S}(\Sigma^{\Bbb N})\\
{}\\
U\mu\in\bigcup_{\varepsilon>0} I(C,\varepsilon)
\endSb
\Bigg(
h(\mu)+\int\varphi\,d\mu
\Bigg)\\
 &=
\sup
\Sb
\mu\in\Cal P_{S}(\Sigma^{\Bbb N})\\
U\mu\in\overset{\,\circ}\to{C}
\endSb
\Bigg(
h(\mu)+\int\varphi\,d\mu
\Bigg)\,.
\qquad\qquad
\qquad\,\,
\text{[since $\cup_{\varepsilon>0} I(C,\varepsilon)=\overset{\,\circ}\to{C}$]}\\
\tag6.28
\endalign
$$

Now note that  
it follows from Proposition 6.4
that
$-\nabla\beta(\Bbb R^{M})\subseteq U(\,\Cal P_{S}(\Sigma^{\Bbb N})\,)$.
Since
$\overset{\,\circ}\to{C}
\cup
(-\nabla\beta(\Bbb R^{M}))
\not=
\varnothing$, we therefore deduce that
$\overset{\,\circ}\to{C}
\cup
U(\,\Cal P_{S}(\Sigma^{\Bbb N})\,)
\not=
\varnothing$, 
and an application of Theorem 5.5 
now gives
 $$
 \align
\sup
\Sb
\mu\in\Cal P_{S}(\Sigma^{\Bbb N})\\
U\mu\in\overset{\,\circ}\to{C}
\endSb
\Bigg(
h(\mu)+\int\varphi\,d\mu
\Bigg)
&=
 -
 \log
 \sigma_{\radius}
 \big(
 \,
 \zeta_{C}^{\dyn,U}(\varphi;\cdot)
 \,
 \big)\,.
 \tag6.29
\endalign
$$

Finally,
combining (6.28) and (6.29) yields
 $-
 \log
 \sigma_{\radius}
 \big(
 \,
 \zeta_{C}^{\dyn\text{-}\scon}(\varphi;\cdot)
 \,
 \big)
\ge
 -
 \log
 \sigma_{\radius}
 \big(
 \,
 \zeta_{C}^{\dyn,U}(\varphi;\cdot)
 \,
 \big)$.
It follows from this inequality that
$\sigma_{\radius}
 \big(
 \,
 \zeta_{C}^{\dyn\text{-}\scon}(\varphi;\cdot)
 \,
 \big)
\le
 \sigma_{\radius}
 \big(
 \,
 \zeta_{C}^{\dyn,U}(\varphi;\cdot)
 \,
 \big)$.
 \hfill$\square$

\bigskip

We can now prove Theorem 6.2.

\bigskip

\noindent{\it Proof of Theorem 6.2}\newline
(1.1) and (2.1):
The statements in Part (1.1) and Part (2.1)
of
Theorem 6.2 follow immediately from 
Proposition 5.2 
and
Proposition 6.5.

\noindent
(1.2) and (2.2):
The statements in Part (1.2) and Part (2.2)
of
Theorem 6.2 follow immediately from 
Part (1.1) and Part (2.1)
using
Corollary 5.4, Corollary 5.6
and
Proposition 6.3.

\noindent
(1.3) and (2.3):
The statements in Part (1.3) and Part (2.3)
of
Theorem 4.2 follow immediately from 
Part (1.2) and Part (2.2)
using
Theorem E.
\hfill$\square$

\bigskip

{\bf 6.3. Multifractal spectra of ergodic Birkhoff averages.}
We first fix $\gamma\in(0,1)$ and define  the metric
 $\distance_{\gamma}$ on $\Sigma^{\Bbb N}$
 as follows.
For $\bold i,\bold j\in\Sigma^{\Bbb N}$ 
with $\bold i\not=\bold j$,
we will write
$\bold i\wedge\bold j$ for the longest common 
prefix of 
$\bold i$ and $\bold j$
(i\.e\.
$
\bold i\wedge\bold j
=
\bold u
$
where $\bold u$ is the unique element in $\Sigma^{*}$
for which there 
are $\bold k,\bold l\in\Sigma^{\Bbb N}$
with
$\bold k=k_{1}k_{2}\ldots$
and
$\bold l=l_{1}l_{2}\ldots$
such that
$k_{1}
 \not= l_{1}$,
 $\bold i
 =
 \bold u\bold k$
and 
$\bold j
 =
 \bold u\bold l$). 
The metric
$\distance_{\gamma}$ is now defined by
 $$
 \distance_{\gamma}(\bold i,\bold j)
 =
 \cases
 0
&\quad 
 \text{if $\bold i=\bold j$;}\\
 \gamma^{|\bold i\wedge\bold j|}
&\quad 
 \text{if $\bold i\not=\bold j$,}
 \endcases
 $$
for $\bold i,\bold j\in\Sigma^{\Bbb N}$;
throughout this section, we equip
 $\Sigma^{\Bbb N}$
with the metric  $\distance_{\gamma}$ 
and continuity and Lipschitz properties of functions $f:\Sigma^{\Bbb N}\to\Bbb R$
from $\Sigma^{\Bbb N}$ to $\Bbb R$
will  always
refer to
the metric  $\distance_{\gamma}$.
Multifractal analysis of Birkhoff
averages has received significant interest
during
the past 10 years, see, for example,
[BaMe,FaFe,FaFeWu,FeLaWu,Oli,Ol3,OlWi2].
The 
multifractal
spectrum
$F_{f}^{\erg}$
of ergodic Birkhoff averages of a continuous function
$f:\Sigma^{\Bbb N}\to\Bbb R$ is defined by
 $$
F_{f}^{\erg}(\alpha)
=
 \dim_{\Haus}
 \pi
 \Bigg\{
 \bold i\in\Sigma^{\Bbb N}
 \,\Bigg|\,
 \lim_{n}
 \frac{1}{n}\sum_{k=0}^{n-1}f(S^{k}\bold i)
 =
 \alpha
 \Bigg\}
$$
for $\alpha\in\Bbb R$.
One of the main problems
in
multifractal analysis of Birkhoff
averages
is the detailed study of the multifractal
spectrum
$F_{f}^{\erg}$.
For example,
Theorem D below is
proved
in different settings and at various levels of generality
in [FaFe,FaFeWu,FeLaWu,Oli,Ol3,OlWi2].

\bigskip

\proclaim{Theorem F [FaFe,FaFeWu,FeLaWu,Oli,Ol3,OlWi2]}
Let $f:\Sigma^{\Bbb N}\to\Bbb R$ be a Lipschitz function.
Let $\Lambda:\Sigma^{\Bbb N}\to\Bbb R$ be defined by (2.5).
Let $C$ be a closed subset of $\Bbb R$.
If the OSC is satisfied, then
 $$
 \dim_{\Haus}
 \pi
 \Bigg\{
 \bold i\in\Sigma^{\Bbb N}
 \,\Bigg|\,
 \,\underset{n}\to\acc
 \frac{1}{n}\sum_{k=0}^{n-1}f(S^{k}\bold i)
\subseteq
C
 \Bigg\}
 =
 \sup
 \Sb
 \mu\in\Cal P_{S}(\Sigma^{\Bbb N})\\
 \int f\,d\mu\in C
 \endSb
 -
 \frac{h(\mu)}{\int\Lambda\,d\mu}\,. 
 $$
In particular, if the OSC is satisfied and $\alpha\in\Bbb R$, then we have 
 $$
  \dim_{\Haus}
 \pi
 \Bigg\{
 \bold i\in\Sigma^{\Bbb N}
 \,\Bigg|\,
 \lim_{n}
 \frac{1}{n}\sum_{k=0}^{n-1}f(S^{k}\bold i)
 =
 \alpha
 \Bigg\}
 =
 \sup
 \Sb
 \mu\in\Cal P_{S}(\Sigma^{\Bbb N})\\
 \int f\,d\mu=\alpha
 \endSb
 -
 \frac{h(\mu)}{\int\Lambda\,d\mu}\,.
 $$ 
\endproclaim

\bigskip

\noindent
As a third 
 application of 
Theorem 2.1
we obtain a 
zeta-function
whose
abscissa of convergence equals the 
multifractal
spectrum
$F_{f}^{\erg}$
of ergodic Birkhoff averages of a Lipschitz function
$f$.
This is the content of the next theorem.

\bigskip

\proclaim{Theorem 6.6.
Multifractal zeta-functinons for
multifractal spectra of of ergodic Birkhoff averages}
Let $f:\Sigma^{\Bbb N}\to\Bbb R$ be a Lipschitz function.
For $C\subseteq\Bbb R$ and 
an continuous function 
$\varphi:\Sigma^{\Bbb N}\to\Bbb R$,   we
define the dynamical
ergodic multifractal zeta-function by
  $$
  \zeta_{C,f}^{\dyn\text{-}\erg}(\varphi;z)
  =
  \sum_{n}
  \frac{z^{n}}{n}
  \left(
  \sum
  \Sb
  |\bold i|=n\\
  {}\\
 \frac{1}{n}\sum_{k=0}^{n-1}f(S^{k}\overline{\bold i})
 \in  C
  \endSb
  \sup_{\bold u\in[\bold i]}
  \exp
  \sum_{k=0}^{n-1}
  \varphi S^{k}\bold u
  \right)\,,
  $$
where we write
$\overline{\bold i}=\bold i\bold i\bold i\ldots$
for  $\bold i\in\Sigma^{*}$.
Let 
 $\Lambda:\Sigma^{\Bbb N}\to\Bbb R$ be defined by (2.5).
 Assume that $C\subseteq\Bbb R$ is closed.

 \roster
  \item"(1)"
 There is a unique real number 
 $\,\,\scr f\,\,(C)$ such that
  $$
  \lim_{r\searrow0}
  \sigma_{\radius}
  \big(
  \,
  \zeta_{B(C,r),f}^{\dyn\text{-}\erg}(\,\,\scr f\,\,(C)\,\Lambda;\cdot)
  \,
  \big)
  =
  1\,.
  $$
It $\alpha\in\Bbb R$ and 
$C=\{\alpha\}$, then we will write
  $\,\,\scr f\,\,(\alpha)
  =
 \,\,\scr f\,\,(C)$.

\item"(2)"
We have
$$
\,\,\scr f\,\,(C)
 =
 \sup_{\alpha\in C}
 \,
 \sup
 \Sb
 \mu\in\Cal P_{S}(\Sigma^{\Bbb N})\\
 \int f\,d\mu=\alpha
 \endSb
 -
 \frac{h(\mu)}{\int\Lambda\,d\mu}\,.
 $$

\item"(3)"
If the OSC is satisfied, then we have
$$
\align
\quad
\,\,\scr f\,\,(C)
 &=
 \dim_{\Haus}
 \pi
 \Bigg\{
 \bold i\in\Sigma^{\Bbb N}
 \,\Bigg|\,
 \,\underset{n}\to\acc
 \frac{1}{n}\sum_{k=0}^{n-1}f(S^{k}\bold i)
\subseteq
C
 \Bigg\}\,.\\
 \endalign
 $$
In particular, if the OSC is satisfied and $\alpha\in\Bbb R$, then we have 
$$
\align
\quad
\,\,\scr f\,\,(\alpha)
 &=
 \dim_{\Haus}
 \pi
 \Bigg\{
 \bold i\in\Sigma^{\Bbb N}
 \,\Bigg|\,
 \lim_{n}
 \frac{1}{n}\sum_{k=0}^{n-1}f(S^{k}\bold i)
 =
 \alpha
 \Bigg\}\,.
 \endalign
 $$
\endroster
\endproclaim

\bigskip

We will now prove  Theorem 6.6.
Recall, that the function 
$\Lambda:\Sigma^{\Bbb N}\to\Bbb R$ is defined by
$\Lambda(\bold i)=\log |DS_{i_{1}}(\pi S\bold i)|$
for
$\bold i=i_{1}i_{2}\ldots\in\Sigma^{*}$.
Define
$U:\Cal P(\Sigma^{\Bbb N})\to\Bbb R$ by
 $$
 U\mu
 =
\int f\,d\mu\,.
\tag6.30
 $$
and note that if $\bold i\in\Sigma^{*}$, then
 $$
 UL_{|\bold i|}[\bold i]
 =
 \Bigg\{
 \frac{1}{|\bold i|}
 \sum_{k=0}^{|\bold i|-1}f(S^{k}(\bold i\bold u))
  \,\Bigg|\,
 \bold u\in\Sigma^{\Bbb N}
 \Bigg\}\,.
 $$
It therefore follows that
 $$
 \align
 \zeta_{C}^{\dyn,U}(\varphi;z)
&=
\sum_{n}
\,\,
\frac{z^{n}}{n}
\,\,
\left(
 \sum
 \Sb
 |\bold i|=n\\
 {}\\
 \forall
 \bold u\in\Sigma^{\Bbb N}
 \,\,:\,\,
 \frac{1}{n}
 \sum_{k=0}^{n-1}f(S^{k}(\bold i\bold u))
 \,
\in
C
 \endSb
   \sup_{\bold u\in[\bold i]}
  \exp
  \sum_{k=0}^{n-1}
  \varphi S^{k}\bold u
  \right)\,.
 \tag6.31
 \endalign 
 $$
 In order to prove Theorem 6.6, we first
 prove  the following  auxiliary result.

\bigskip

\proclaim{Proposition 6.7}
Let $U$ 
be defined by (6.30).
Fix a continuous function $\varphi:\Sigma^{\Bbb N}\to\Bbb R$.
\roster
\item"(1)"
There is a sequence $(\Delta_{n})_{n}$ with $\Delta_{n}>0$ 
for all $n$ and 
$\Delta_{n}\to 0$ such that
for all closed subsets $C$ of $\Bbb R$
and
for
all
$n\in\Bbb N$, $\bold i\in\Sigma^{n}$ and $\bold u\in\Sigma^{\Bbb N}$, we have
 $$
 \align
 \dist
 \Bigg(
 \,
 \frac{1}{n}
 \sum_{k=0}^{n-1}f(S^{k}(\bold i\bold u))
 \,,\,
 C
 \,
 \Bigg)
&\le
\dist
\Bigg(
 \,
 \frac{1}{n}
 \sum_{k=0}^{n-1}f(S^{k}(\overline{\bold i})) \,,\,
 C
 \,
 \Bigg)
\,\,\,\,\,\,
+
\,\,
\Delta_{n}\,,\\
\dist
\Bigg(
 \,
 \frac{1}{n}
 \sum_{k=0}^{n-1}f(S^{k}(\overline{\bold i})) \,,\,
 C
 \,
 \Bigg)
&\le
 \dist
 \Bigg(
 \,
 \frac{1}{n}
 \sum_{k=0}^{n-1}f(S^{k}(\bold i\bold u))
 \,,\,
 C
 \,
 \Bigg)
 \,\,
+
\,\,
\Delta_{n}\,;
\endalign
$$
recall, that for $\bold i\in\Sigma^{*}$, we write
$\overline{\bold i}=\bold i\bold i\bold i\ldots$.

\item"(2)"
Let $C$ be a closed subset of $\Bbb R^{M}$.
We have
 $$
 \lim_{r\searrow0}
 \sigma_{\radius}
 \big(
 \,
 \zeta_{B(C,r),f}^{\dyn\text{-}\erg}(\varphi;\cdot) 
 \,
 \big)
 =
 \lim_{r\searrow0}
 \sigma_{\radius}
 \big(
 \,
 \zeta_{B(C,r)}^{\dyn,U}(\varphi;\cdot)
 \,
 \big)\,.
$$

\endroster
\endproclaim
\noindent{\it Proof}\newline
\noindent
(1)
Let $\Lip(f)$ denote the Lipschitz constant of $f$. 
It is clear that
for
all $n\in\Bbb N$,
$\bold i\in\Sigma^{n}$ and $\bold u\in\Sigma^{\Bbb N}$, we have
 $$
 \align
 \Bigg|
  \frac{1}{n}
 \sum_{k=0}^{n-1}f(S^{k}(\overline{\bold i})) 
 -
  \frac{1}{n}
 \sum_{k=0}^{n-1}f(S^{k}(\bold i\bold u)) 
 \Bigg|
 &\le
 \frac{1}{n}
 \sum_{k=0}^{n-1}
 |f(S^{k}(\overline{\bold i}))  - f(S^{k}(\bold i\bold u)) |\\
&\le
\Lip(f)
 \frac{1}{n}
 \sum_{k=0}^{n-1}
 \distance_{\gamma}
 \big(
 \,
 S^{k}(\overline{\bold i}),S^{k}(\bold i\bold u)
 \,
 \big) \\ 
&\le
\Lip(f)
 \frac{1}{n}
 \sum_{k=0}^{n-1}\gamma^{k}\\
&\le
\Lip(f)
 \frac{1}{n(1-\gamma)}\,.
 \tag6.32
  \endalign
 $$
It is not difficult to see that the desired result follows from (6.32).

\noindent
(2)
This statement follows from Part (1)
by an argument very similar to the proofs of Part (2) and Part (3)
in
Proposition 6.5, and the proof is therefore omitted.
\hfill$\square$

\bigskip

We can now prove Theorem 6.6.

\bigskip

\noindent{\it Proof of Theorem 6.6}\newline
(1) This statement
follows immediately from 
Proposition 5.2 and Proposition 6.7.

\noindent
(2)  This statement
follows immediately from 
Part (1)
using
Corollary 5.4.
\hfill$\square$

\noindent
(3)  This statement
follows immediately from 
Part (2)
using
 Theorem F.
\hfill$\square$

  \bigskip


\heading
{
7. Proofs. Preliminary results:  the modified multifractal pressure}
\endheading

In this section we introduce our main technical tool, namely, the
modified multifractal pressure;
see definition (7.2) below.
The two main results is this section 
are Theorem 7.3
providing a variational principle for the modified multifractal pressure
and 
Theorem 
7.5
showing that
the 
multifractal pressure and the modified multifractal pressure
are (almost)
comparable.
Both Theorem 7.3 and Theorem 7.5 play
major roles in the 
in the proof 
of Theorem 5.3 in Section 8
and 
in the proof of Theorem 5.5 in Section 9.

We first  define the modified multifractal pressure.
We start by introducing some notation.
If $\bold i\in\Sigma^{*}$, then we define 
$\overline{\bold i}\in\Sigma^{\Bbb N}$ by
$\overline{\bold i}=\bold i\bold i\ldots$.
We also
define $M_{n}:\Sigma^{\Bbb N}\to\Cal P_{S}(\Sigma^{\Bbb N})$ by
 $$
 \align
 M_{n}\bold i
&=
 L_{n}\left(\,\overline{\bold i|n}\,\right)
 = 
 \frac{1}{n}\sum_{k=0}^{n-1}\delta_{S^{k}(\,\overline{\bold i|n}\,)}
 \tag7.1
 \endalign
 $$
for $\bold i\in\Sigma^{\Bbb N}$;
recall, that
the map
$L_{n}:\Sigma^{\Bbb N}\to\Cal P(\Sigma^{\Bbb N})$ is defined in (5.1).
Furthermore,
note that
if
$\bold i\in\Sigma^{\Bbb N}$,
then
$M_{n}\bold i$ is shift invariant,
i\.e\.
$M_{n}$ 
maps 
$\Sigma^{\Bbb N}$ into $\Cal P_{S}(\Sigma^{\Bbb N})$
as claimed.
Next, let $P$ denote the probability measure on $\Sigma^{\Bbb N}$ given by
 $$
 P={\underset{\Bbb N}\to \X}\,\,\sum_{i=1}^{N}\frac{1}{N}\delta_{i}\,.
 $$
 For a continuous function $\varphi:\Sigma^{\Bbb N}\to\Bbb R$,
 we
define $F_{\varphi}:\Cal P_{S}(\Sigma^{\Bbb N})\to\Bbb R$ by
 $$
 F_{\varphi}(\mu)=\int\varphi\,d\mu\,.
 $$
Observe that since $\varphi$ is bounded, 
i\.e\. $\|\varphi\|_{\infty}<\infty$,
we conclude that
$\|F_{\varphi}\|_{\infty}
\le
\|\varphi\|_{\infty}<\infty$.
 Next, for a positive integer $n$, define
probability measures
$P_{n},Q_{\varphi,n}\in\Cal P\big(\,\Cal P_{S}(\Sigma^{\Bbb N})\,\big)$ by
 $$
 \align
 {}\\
 P_{n}
&=
 P\circ M_{n}^{-1}\,,\\
 {}\\
 Q_{\varphi,n}(E)
&=
 \frac
 {\int_{E}\exp(nF_{\varphi})\,dP_{n}}
 {\int\exp(nF_{\varphi})\,dP_{n}}
 \quad
 \text{for Borel subsets $E$ of $\Cal P_{S}(\Sigma^{\Bbb N})$.}\\
 {}
 \endalign
 $$
Finally, we define  
modified multifractal pressures as follows.
Namely, for $C\subseteq X$,
we define 
 the 
 modified lower and upper 
  mutifractal pressure
of 
 $\varphi$
associated with the space $X$ and the map $U$ and  by
 $$
 \aligned
  \underline Q_{C}^{U}(\varphi)
 &=
 \,
  \liminf_{n}
 \,
  \,\,
  \frac{1}{n}
 \,\,
 \log
   \sum
  \Sb
  |\bold i|=n\\
  {}\\
  UM_{n}[\bold i]\subseteq C
  \endSb
  \sup_{\bold u\in[\bold i]}
  \,\,
   \exp
 \,\,
 \sum_{k=0}^{n-1}\varphi S^{k}\bold u\,,\\
   \overline Q_{C}^{U}(\varphi)
 &=
  \limsup_{n}
  \,\,
  \frac{1}{n}
 \,\,
 \log
   \sum
  \Sb
 |\bold i|=n\\
  {}\\
  UM_{n}[\bold i]\subseteq C
  \endSb
   \sup_{\bold u\in[\bold i]}
   \,\,
   \exp
 \,\,
 \sum_{k=0}^{n-1}\varphi S^{k}\bold u\,.
 \endaligned
 \tag7.2
 $$

We now turn towards the proof of the first main result in
  this section, namely, Theorem 7.3 providing 
  a variational principle for the modified multifractal pressure.
The proof of Theorem 7.3 is based on
large deviation theory. 
In particular, 
we need Varadhan's [Va]
large deviation theorem (Theorem 7.1.(i) below),
and a non-trivial application of this (namely Theorem 7.1.(ii) below) 
providing first order asymptotics of certain
\lq\lq Boltzmann distributions".
However, we begin with a definition.

\bigskip

\proclaim{Definition}
Let $X$
be a complete separable metric space
and
let $(P_{n})_{n}$ be a sequence of probability measures on $X$. 
Let
$(a_{n})_{n}$ be a sequence of positive numbers with
$a_{n}\to\infty$ 
and let
$I:X\to[0,\infty]$ be a lower semicontinuous function
with compact level sets.
The sequence $(P_{n})_{n}$ is said to have the large deviation property
with constants $(a_{n})_{n}$ and rate function $I$ if the following 
two condistions
hold:
\roster
\item"(i)" For each closed subset $K$ of $X$, we have
 $$
 \limsup_{n}\frac{1}{a_{n}}\log P_{n}(K)\le-\inf_{x\in K}I(x)\,;
 $$
\item"(ii)" For each open subset $G$ of $X$, we have
 $$
 \liminf_{n}\frac{1}{a_{n}}\log P_{n}(G)\ge-\inf_{x\in G}I(x)\,.
 $$
\endroster
\endproclaim

\bigskip

\proclaim{Theorem 7.1}
Let $X$ be a complete separable metric space
and
let $(P_{n})_{n}$ be a sequence of probability measures on $X$.
Assume that the sequence $(P_{n})_{n}$
has the large deviation property
with constants $(a_{n})_{n}$ and rate function $I$.
Let $F:X\to\Bbb R$ be a continuous function
satisfying the following two conditions:
\roster
\item"(i)"
For all $n$, we have
 $$
 \int\exp(a_{n}F)\,dP_{n}
 <
 \infty\,.
 $$
\item"(ii)"
We have
 $$
 \lim_{M\to\infty}\,\,
 \limsup_{n}\,\,
 \frac{1}{a_{n}}
 \log
 \int_{\{M\le F\}}
 \exp(a_{n}F)\,dP_{n}
 =
 -\infty\,.
 $$
 \endroster
(Observe that the  Conditions (i)--(ii)
 are satisfied if $F$ is bounded.)
Then the following statements hold.
\roster
\item"(1)" We have
 $$
 \lim_{n}\,\,
 \frac{1}{a_{n}}
 \log
 \int
 \exp(a_{n}F)\,dP_{n}
 =
 -\inf_{x\in X}(I(x)-F(x))\,.
 $$
\item"(2)"
For each $n$ define a probability measure $Q_{n}$ on $X$ by
 $$
 Q_{n}(E)
 =
 \,\,
 \frac
 {\int_{E}\exp(a_{n}F)\,dP_{n}}
 {\int\exp(a_{n}F)\,dP_{n}}\,.
 $$
Then the sequence $(Q_{n})_{n}$
has the large deviation property with constants $(a_{n})_{n}$
and rate function
$(I-F)-\inf_{x\in X}(I(x)-F(x))$.
\endroster
\endproclaim 
\noindent{\it Proof}\newline
\noindent 
Statement (1) follows from [El, Theorem II.7.1] or [DeZe, Theorem 4.3.1], and
statement (2) follows from [El, Theorem II.7.2].
\hfill$\square$

 \bigskip

 \noindent
 Before stating and proving Theorem 7.3, we establish the following auxiliary result.

\bigskip

\proclaim{Theorem 7.2}
Let $X$ be a metric space and let $U:\Cal P(\Sigma^{\Bbb N})\to X$ be 
continuous with respect to the weak topology.
Let $C\subseteq X$ be a  subset of $X$.
Fix a continuous function $\varphi:\Sigma^{\Bbb N}\to\Bbb R$.
Then there is a constant $c$
such that
for all positive integers $n$, we have
 $$
 \align
 \sum
   \Sb
   |\bold k|=n\\
   {}\\
   UM_{n}[\bold k]\subseteq C
	 \endSb
 \sup_{\bold u\in[\bold k]}
 \,\,
 \exp
 \sum_{k=0}^{n-1}\varphi S^{k}\bold u
&\le
c
\,\,
N^{n}
\,\,
Q_{\varphi,n}\Big(\{U\in C\}\Big)
\,\,
\int\exp(nF_{\varphi})\,dP_{n}\,,\\
{}\\
 \sum
   \Sb
   |\bold k|=n\\
   {}\\
   UM_{n}[\bold k]\subseteq C
	 \endSb
 \sup_{\bold u\in[\bold k]}
 \,\,	 
 \exp
 \sum_{k=0}^{n-1}\varphi S^{k}\bold u
&\ge
\frac{1}{c}
\,\,
N^{n}
\,\,
Q_{\varphi,n}\Big(\{U\in C\}\Big)
\,\,
\int\exp(nF_{\varphi})\,dP_{n}\,.
\endalign
 $$
\endproclaim
\noindent{\it Proof}
\newline
\noindent 
For each positive integer $n$ and each $\bold i$ with $|\bold i|=n$, 
we write
$s_{\bold i}
=
\sup_{\bold u\in[\bold i]}
\exp
 \sum_{k=0}^{n-1}\varphi S^{k}\bold u$
 for sake of brevity.
Let $C$ be a subset of $X$.
For each positive integer $n$, we clearly have
 $$
 \align
 &\int
 \limits_{
 \big\{
 \bold j\in\Sigma^{\Bbb N}
 \,\big|\,
 U M_{n}[\bold j|n]\subseteq C
 \big\}
 }
s_{\bold i|n}
 \,dP(\bold i)\\
&\qquad\qquad
=
\qquad\qquad
\,\,\,\,
 \sum_{|\bold k|=n}
 \,\,\,\,
 \int\limits_{
 [\bold k]
 \,\,\cap\,\,
 \big\{
 \bold j\in\Sigma^{\Bbb N}
 \,\big|\,
 U M_{n}[\bold j|n]\subseteq C
  \big\}
 }
 s_{\bold i|n}
 \,dP(\bold i)\\
\allowdisplaybreak 
&\qquad\qquad
=
\qquad\qquad
\,\,\,\,
 \sum_{|\bold k|=n}
 s_{\bold k}
 \,\,
 P
 \Big(
 \,
  [\bold k]
 \,\,\cap\,\,
 \Big\{
 \bold j\in\Sigma^{\Bbb N}
 \,\Big|\,
 U M_{n}[\bold j|n]\subseteq C
 \Big\}
 \,
 \Big)\\
\allowdisplaybreak 
&\qquad\qquad
=
 \sum
   \Sb
   |\bold k|=n\\
   {}\\
   [\bold k]
 \,\,\cap\,\,
 \big\{
 \bold j\in\Sigma^{\Bbb N}
 \,\big|\,
 U M_{n}[\bold j|n]\subseteq C
 \big\}
 \not=\varnothing
	 \endSb
 s_{\bold k}
 \,\,
 P
 \Big(
 \,
  [\bold k]
 \,\,\cap\,\,
 \Big\{
 \bold j\in\Sigma^{\Bbb N}
 \,\Big|\,
 U M_{n}[\bold j|n]\subseteq C
  \Big\}
 \,
 \Big)\,.\\
& \tag7.3
 \endalign
 $$
Now observe that if
$\bold k\in\Sigma^{*}$ with $|\bold k|=n$ 
and
 $[\bold k]
 \cap
 \{
 \bold j\in\Sigma^{\Bbb N}
 \,|\,
 U M_{n}[\bold j|n]\subseteq C
 \}
 \not=\varnothing$,
 then there is 
 $\bold u\in
 [\bold k]
 \cap
 \{
 \bold j\in\Sigma^{\Bbb N}
 \,|\,
 U M_{n}[\bold j|n]\subseteq C
 \}$.
 Since
  $\bold u\in
 [\bold k]$, we conclude that
 $\bold u=\bold k\bold v$ for some $£\bold v\in\Sigma^{\Bbb N}$.
 Next, since
 also
 $\bold k\bold v 
 =
 \bold u
 \in
 \{
 \bold j\in\Sigma^{\Bbb N}
 \,|\,
 U M_{n}[\bold j|n]\subseteq C
 \}$, we conclude
 that
 $U M_{n}[\bold k]
 =
 U M_{n}[(\bold k\bold v)|n]
 =
 U M_{n}[\bold u|n]
 \subseteq 
 C$. This shows that
 $$
 \align
& \sum
   \Sb
   |\bold k|=n\\
   {}\\
   [\bold k]
 \,\,\cap\,\,
 \big\{
 \bold j\in\Sigma^{\Bbb N}
 \,\big|\,
 U M_{n}[\bold j|n]\subseteq C
 \big\}
 \not=\varnothing
	 \endSb
 s_{\bold k}
 \,\,
 P
 \Big(
 \,
  [\bold k]
 \,\,\cap\,\,
 \Big\{
 \bold j\in\Sigma^{\Bbb N}
 \,\Big|\,
 U M_{n}[\bold j|n]\subseteq C
  \Big\}
 \,
 \Big)\\
 &\qquad\qquad
 =
 \sum
   \Sb
   |\bold k|=n\\
   {}\\
   [\bold k]
 \,\,\cap\,\,
 \big\{
 \bold j\in\Sigma^{\Bbb N}
 \,\big|\,
 U M_{n}[\bold j|n]\subseteq C
 \big\}
 \not=\varnothing\\
 {}\\
 UM_{n}[\bold k]\subseteq C
	 \endSb
 s_{\bold k}
 \,\,
 P
 \Big(
 \,
  [\bold k]
 \,\,\cap\,\,
 \Big\{
 \bold j\in\Sigma^{\Bbb N}
 \,\Big|\,
 U M_{n}[\bold j|n]\subseteq C
  \Big\}
 \,
 \Big)\,.\\
 \tag7.4
 \endalign
 $$
Combining (7.3) and (7.4)  gives
 $$
 \align
&\int
 \limits_{
 \big\{
 \bold j\in\Sigma^{\Bbb N}
 \,\big|\,
 U M_{n}[\bold j|n]\subseteq C
 \big\}
 }
s_{\bold i|n}
 \,dP(\bold i)\\
&\qquad\qquad
=
 \sum
   \Sb
   |\bold k|=n\\
   {}\\
   [\bold k]
 \,\,\cap\,\,
 \big\{
 \bold j\in\Sigma^{\Bbb N}
 \,\big|\,
 U M_{n}[\bold j|n]\subseteq C
 \big\}
 \not=\varnothing\\
  {}\\
 UM_{n}[\bold k]\subseteq C
	 \endSb
 s_{\bold k}
 \,\,
 P
 \Big(
 \,
  [\bold k]
 \,\,\cap\,\,
 \Big\{
 \bold j\in\Sigma^{\Bbb N}
 \,\Big|\,
 U M_{n}[\bold j|n]\subseteq C
  \Big\}
 \,
 \Big)\,.\\
 &\qquad\qquad
=
\qquad\quad\,\,\,\,
 \sum
   \Sb
   |\bold k|=n\\
  {}\\
 UM_{n}[\bold k]\subseteq C
	 \endSb
 s_{\bold k}
 \,\,
 P
 \Big(
 \,
  [\bold k]
 \,\,\cap\,\,
 \Big\{
 \bold j\in\Sigma^{\Bbb N}
 \,\Big|\,
 U M_{n}[\bold j|n]\subseteq C
  \Big\}
 \,
 \Big)\,.\\
 \tag7.5
 \endalign
 $$

However,
if
$\bold k\in\Sigma^{*}$ with
 $|\bold k|=n$
 and
  $ U M_{n}[\bold k]\subseteq C$,
  then  it is clear that
  $[\bold k]\subseteq
 \{
 \bold j\in\Sigma^{\Bbb N}
 \,|\,
  U M_{n}[\bold j|n]\subseteq C
 \}$,
 whence
 $[\bold k]\cap
 \{
 \bold j\in\Sigma^{\Bbb N}
 \,|\,
  U M_{n}[\bold j|n]\subseteq C
 \}
 =
 [\bold k]$.
 This and (7.5) now imply that
 $$
 \align
  \int
 \limits_{
 \big\{
 \bold j\in\Sigma^{\Bbb N}
 \,\big|\,
 U M_{n}[\bold j|n]\subseteq C
 \big\}
 }
&s_{\bold i|n}
 \,dP(\bold i)\\
&=
 \sum
   \Sb
   |\bold k|=n\\
   {}\\
   U M_{n}[\bold k]\subseteq C
	 \endSb
 s_{\bold k}
 \,\,
 P
 \Big(
 \,
  [\bold k]
 \,\,\cap\,\,
 \Big\{
 \bold j\in\Sigma^{\Bbb N}
 \,\Big|\,
 U M_{n}[\bold j|n]\subseteq C
  \Big\}
 \,               
 \Big)\\ 
&=
 \sum
   \Sb
   |\bold k|=n\\
   {}\\
 U M_{n}[\bold k]\subseteq C
	 \endSb
 s_{\bold k}
 \,\,
 P\big(\,
 [\bold k]
 \,\big)\\
\allowdisplaybreak 
&=
 \sum
   \Sb
  |\bold k|=n\\
  {}\\
   U M_{n}[\bold k]\subseteq C
	 \endSb
 s_{\bold k}\,\frac{1}{N^{n}}\,,
 \endalign
 $$
whence
 $$
 \align
 \sum
   \Sb
  |\bold k|=n\\
  {}\\
   U M_{n}[\bold k]\subseteq C
	 \endSb
 s_{\bold k}
&=
 N^{n}
 \int
 \limits_{
 \big\{
 \bold j\in\Sigma^{\Bbb N}
 \,\big|\,
 U M_{n}[\bold j|n]\subseteq C
 \big\}
 }
s_{\bold i|n}
 \,dP(\bold i)\,.
 \tag7.6
 \endalign
 $$

It follows from the 
Principle of
Bounded Distortion 
(see, for example, [Bar,Fa2])
that there is a constant $c>0$
such that
if $n\in\Bbb N$, $\bold i\in\Sigma^{n}$
and $\bold u,\bold v\in[\bold i]$, then
$\frac{1}{c}
\le
\frac
{
\exp\sum_{k=0}^{n-1}\varphi S^{k}\bold u
}
{
\exp\sum_{k=0}^{n-1}\varphi S^{k}\bold v
}
\le
c$.
In particular, this implies that
for all 
$n\in\Bbb N$ and for all $\bold i\in\Sigma^{n}$, we have
 $$
 \frac{1}{c}
\exp
\sum_{k=0}^{n-1}\varphi S^{k}\overline{\bold i}
\le
s_{\bold i}
\le
c
\exp
\sum_{k=0}^{n-1}\varphi S^{k}\overline{\bold i}\,.
\tag7.7
$$

\medskip

{\it Claim 1. For all positive integers $n$, we have
 $$
 \align
  \sum
   \Sb
  |\bold k|=n\\
  {}\\
   U M_{n}[\bold k]\subseteq B(C, r)
	 \endSb
 s_{\bold k}
&\le
\,
 c
 \,
 N^{n}
 \int\limits_{
 \big\{
 \bold j\in\Sigma^{\Bbb N}
 \,\big|\,
 U M_{n}[\bold j|n]\subseteq C
 \big\}
 }
 \exp\left(n F_{\varphi}(M_{n}\bold i)\right)
 \,dP(\bold i)\,,
 \tag7.8\\
  \sum
   \Sb
  |\bold k|=n\\
  {}\\
   U M_{n}[\bold k]\subseteq B(C, r)
	 \endSb
 s_{\bold k}
&\ge
 \frac{1}{c}
 \,
 N^{n}
 \int\limits_{
 \big\{
 \bold j\in\Sigma^{\Bbb N}
 \,\big|\,
 U M_{n}[\bold j|n]\subseteq C
 \big\}
 }
 \exp\left(n F_{\varphi}(M_{n}\bold i)\right)
 \,dP(\bold i)\,.
 \tag7.9
 \endalign
 $$
}

\noindent
{\it Proof of Claim 1.}
It follows from (7.6) and (7.7) that
if $n$ is a positive integer, then
we have
 $$
 \align
  \sum
   \Sb
  |\bold k|=n\\
  {}\\
   U M_{n}[\bold k]\subseteq C
	 \endSb
 s_{\bold k}^{t}
&\le 
 c
 \,
 N^{n}
 \int\limits_{
 \big\{
 \bold j\in\Sigma^{\Bbb N}
 \,\big|\,
U M_{n}[\bold j|n]\subseteq C
 \big\}
 }
 \exp
 \Bigg(
 \sum_{k=0}^{n-1}\varphi S^{k}\left(\,\overline{\bold i|n}\,\right)
 \Bigg)
 \,dP(\bold i)\\
\allowdisplaybreak 
&=
 c
 \,
 N^{n}
 \int\limits_{
 \big\{
 \bold j\in\Sigma^{\Bbb N}
 \,\big|\,
U M_{n}[\bold j|n]\subseteq C
 \big\}
  }
 \exp
 \Bigg(
 n\int\varphi\,d(M_{n}\bold i)
 \Bigg)
 \,dP(\bold i)\\
\allowdisplaybreak 
&=
 c
 \,
 N^{n}
 \int\limits_{
 \big\{
 \bold j\in\Sigma^{\Bbb N}
 \,\big|\,
 U M_{n}[\bold j|n]\subseteq C
 \big\}
 }
 \exp\left(n F_{\varphi}(M_{n}\bold i)\right)
 \,dP(\bold i)\,.
 \endalign
 $$
This proves inequality (7.8).
Inequality (7.9) is proved similarly.
This completes the proof of Claim 1.

\medskip

{\it Claim 2.
For all positive integers $n$, we have
$\{
 \bold j\in\Sigma^{\Bbb N}
 \,|\,
 U M_{n}[\bold j|n]\subseteq C
 \}
=
\{
 \bold j\in\Sigma^{\Bbb N}
 \,|\,
 U M_{n}\bold j\subseteq C
 \}$.
}
 
\noindent
{\it Proof of Claim 2.}
Indeed, it is clear that
$\{
 \bold j\in\Sigma^{\Bbb N}
 \,|\,
 U M_{n}[\bold j|n]\subseteq C
 \}
\subseteq
\{
 \bold j\in\Sigma^{\Bbb N}
 \,|\,
 U M_{n}\bold j\in C
 \}$.
We will now prove the reverse inclusion.
We therefore fix  
$ \bold j\in\Sigma^{\Bbb N}$ with
$U M_{n}\bold j\in C$.
We must now prove that
 $U M_{n}[\bold j|n]\subseteq C$.
 In order to do this, we let $\bold u\in[\bold j|n]$.
 Since $\bold u\in[\bold j|n]$, 
 we conclude that
 $\bold u|n=\bold j|n$,
 whence
 $U M_{n}\bold u
 =
 U L_{n}(\,\overline{\bold u|n}\,)
 =
 U L_{n}(\,\overline{\bold j|n}\,)
 =
 U M_{n}\bold j
 \in
 C$.
 This completes the proof of Claim 2.

 \medskip

 For all positive integers $n$,
 we now  deduce from 
 Claim 1 and Claim 2
 that
 $$
 \align
  \sum
   \Sb
  |\bold k|=n\\
  {}\\
   U M_{n}[\bold k]\subseteq C
	 \endSb
 s_{\bold k}
&\le
 c
 \,
 N^{n}
 \int\limits_{
 \big\{
 \bold j\in\Sigma^{\Bbb N}
 \,\big|\,
 U M_{n}[\bold j|n]\subseteq C
 \big\}
 }
 \exp\left(n F_{\varphi}(M_{n}\bold i)\right)
 \,dP(\bold i)\\ 
&=
 c
 \,
 N^{n}
 \int\limits_{
 \big\{U M_{n}\in C\big\}
 }
 \exp\left(nF_{\varphi}(M_{n}\bold i)\right)
 \,dP(\bold i)\\
\allowdisplaybreak 
&=
c
 \,
 N^{n}
 \int\limits_{
 \big\{U\in C\big\}
 }
 \exp\left(nF_{\varphi}\right)
 \,dP_{n}\\
\allowdisplaybreak 
&=
 c
 \,
 N^{n}
 \,\,
 Q_{\varphi,n}\Big(\{U\in C\}\Big)
 \,\,
 \int\exp\left(nF_{\varphi}\right)\,dP_{n}\,.
 \endalign
 $$

Similarly, we prove that
for all positive integers $n$, we have
  $$
 \align
  \sum
   \Sb
  |\bold k|=n\\
  {}\\
   U M_{n}[\bold k]\subseteq C
	 \endSb
 s_{\bold k}
&\ge
\frac{1}{c}
 \,
 N^{n}
 \,\,
 Q_{\varphi,n}\Big(\{U\in C\}\Big)
 \,\,
 \int\exp\left(nF_{\varphi}\right)\,dP_{n}\,.
 \endalign
 $$
This completes the proof of Theorem 7.2.
\hfill$\square$

\bigskip

\noindent
We can now state and prove the first main result in this section, namely, Theorem 7.3.

\bigskip

\proclaim{Theorem 7.3. 
The variational principle for the modified
multifractal pressure.}
Let $X$ be a metric space and let $U:\Cal P(\Sigma^{\Bbb N})\to X$ be 
continuous with respect to the weak topology.
Let $C\subseteq X$ be a  subset of $X$.
Fix a continuous function $\varphi:\Sigma^{\Bbb N}\to\Bbb R$.
\roster
\item"(1)"
If $G$ is an open subset of $X$, then
 $$
\underline Q_{G}^{U}(\varphi)
\ge
\sup
\Sb
\mu\in\Cal P_{S}(\Sigma^{\Bbb N})\\
{}\\
U\mu\in G
\endSb
\Bigg(
h(\mu)+\int\varphi\,d\mu
\Bigg)\,.
\tag7.10
$$
\item"(2)"
If $K$ is a closed subset of $X$, then
$$
\overline Q_{K}^{U}(\varphi)
\le
\sup
\Sb
\mu\in\Cal P_{S}(\Sigma^{\Bbb N})\\
{}\\
U\mu\in K
\endSb
\Bigg(
h(\mu)+\int\varphi\,d\mu
\Bigg)\,.
\tag7.11
$$ 
\endroster
\endproclaim

\noindent{\it Proof}\newline
\noindent
We introduce the 
simplified notation from the proof of Theorem 7.2, i\.e\.
for each positive integer $n$ and each $\bold i$ with $|\bold i|=n$, 
we write
$s_{\bold i}
=
\sup_{\bold u\in[\bold i]}
\exp
 \sum_{k=0}^{n-1}\varphi S^{k}\bold u$.
First
note that it follows immediately from Theorem 7.2 that
 $$
 \aligned
 \liminf_{n}
 \frac{1}{n}
 \log
 \sum
   \Sb
   |\bold i|=n\\
   {}\\
   U M_{n}[\bold i]\subseteq G
	 \endSb
 s_{\bold i}
&\ge
 \log N
 \,+\,
 \limsup_{n}
 \frac{1}{n}
 \log Q_{\varphi,n}\Big(\{U\in G \}\Big)\\
&\qquad\qquad
  \qquad\qquad
 \,+\,
 \limsup_{n}
 \frac{1}{n}
 \log 
 \int
 \exp
 \left(
 nF_{\varphi}
 \right)
 \,dP_{n}\,,\\
\limsup_{n}
 \frac{1}{n}
 \log
 \sum
   \Sb
   |\bold i|=n\\
   {}\\
   U M_{n}[\bold i]\subseteq K
	 \endSb
 s_{\bold i}
&\le
 \log N
 \,+\,
 \limsup_{n}
 \frac{1}{n}
 \log Q_{\varphi,n}\Big(\{U\in K \}\Big)\\
&\qquad\qquad
  \qquad\qquad
 \,+\,
 \limsup_{n}
 \frac{1}{n}
 \log 
 \int
 \exp
 \left(
 nF_{\varphi}
 \right)
 \,dP_{n}\,.
 \endaligned
 \tag7.12
 $$

Next, we observe that it follows from [El]
that the sequence 
$(P_{n}=P\circ M_{n}^{-1})_{n}
\subseteq
\Cal P\big(\,\Cal P_{S}(\Sigma^{\Bbb N})\,\big)$ 
has the large deviation property with respect to
the sequence
$(n)_{n}$ and rate function 
$I:\Cal P_{S}(\Sigma^{\Bbb N})\to\Bbb R$ given by
$I(\mu)=\log N-h(\mu)$.
We therefore 
conclude from Part (1) of Theorem 7.1 that
 $$
  \lim_{n}
 \frac{1}{n}
 \log 
 \int
 \exp
 \left(
 nF_{\varphi}
 \right)
 \,dP_{n}
 =
 -
  \inf_{\nu\in\Cal P_{S}(\Sigma^{\Bbb N})}(I(\nu)-F_{\varphi}(\nu))\,.
  \tag7.13
 $$
Also,
since the sequence 
$(P_{n}=P\circ M_{n}^{-1})_{n}
\subseteq
\Cal P\big(\,\Cal P_{S}(\Sigma^{\Bbb N})\,\big)$ 
has the large deviation property with respect to
the sequence
$(n)_{n}$ and rate function 
$I:\Cal P_{S}(\Sigma^{\Bbb N})\to\Bbb R$ given by
$I(\mu)=\log N-h(\mu)$, 
we conclude from Part (2) of Theorem 7.1
that the sequence $(Q_{\varphi,n})_{n}$ has the large deviation property with 
respect to 
the sequence
$(n)_{n}$ and rate function
$(I-F_{\varphi})-\inf_{\nu\in\Cal P_{S}(\Sigma^{\Bbb N})}(I(\nu)-F_{\varphi}(\nu))$.
As the set
$\{U\in G\}
=
U^{-1}G$
is open
and
the set $\{U\in K\}
=
U^{-1}K$
is closed,
 it therefore follows from the large deviation property that
 $$
 \aligned
 \limsup_{n}
 \frac{1}{n}
&\log Q_{\varphi,n}\Big(\{U\in G \}\Big)\\
&\ge
-
 \inf
  \Sb
	\mu\in\Cal P_{S}(\Sigma^{\Bbb N})\\
	{}\\
	U\mu\in G
	\endSb
 \Bigg(
 (I(\mu)-F_{\varphi}(\mu))
 -
 \inf_{\nu\in\Cal P_{S}(\Sigma^{\Bbb N})}
 (I(\nu)-F_{\varphi}(\nu))
 \Bigg)\,.\\
 \limsup_{n}
 \frac{1}{n}
&\log Q_{\varphi,n}\Big(\{U\in K \}\Big)\\
&\le
-
 \inf
  \Sb
	\mu\in\Cal P_{S}(\Sigma^{\Bbb N})\\
	{}\\
	U\mu\in K
	\endSb
 \Bigg(
 (I(\mu)-F_{\varphi}(\mu))
 -
 \inf_{\nu\in\Cal P_{S}(\Sigma^{\Bbb N})}
 (I(\nu)-F_{\varphi}(\nu))
 \Bigg)\,. 
 \endaligned
 \tag7.14
 $$

Combining (7.12). (7.13) and (7.14) now yields
 $$
 \align
 \limsup_{n}
 \frac{1}{n}
 \log
 \sum
   \Sb
   |\bold i|=n\\
   {}\\
   U M_{n}[\bold i]\subseteq G
	 \endSb
 s_{\bold i}
&\ge
 \log N
 \,+\,
 \limsup_{n}
 \frac{1}{n}
 \log Q_{\varphi,n}\Big(\{U\in G \}\Big)\\
&\qquad\qquad
  \qquad\qquad
 \,+\,
 \limsup_{n}
 \frac{1}{n}
 \log 
 \int
 \exp
 \left(
 nF_{\varphi}
 \right)
 \,dP_{n}\\
&\ge
 \log N\\
&\qquad
-
 \inf
  \Sb
	\mu\in\Cal P_{S}(\Sigma^{\Bbb N})\\
	{}\\
	U\mu\in G
	\endSb
 \Bigg(
 (I(\mu)-F_{\varphi}(\mu))
 -
 \inf_{\nu\in\Cal P_{S}(\Sigma^{\Bbb N})}
 (I(\nu)-F_{\varphi}(\nu))
 \Bigg)\\
&\qquad\qquad
  \qquad\qquad
 \,-\,
 \inf_{\nu\in\Cal P_{S}(\Sigma^{\Bbb N})}(I(\nu)-F_{\varphi}(\nu))\\
&=
 \log N
 \,+\,
 \sup
  \Sb
	\mu\in\Cal P_{S}(\Sigma^{\Bbb N})\\
	{}\\
	U\mu\in G
	\endSb
 (F_{\varphi}(\mu)-I(\mu))\\
&=
 \sup
  \Sb
	\mu\in\Cal P_{S}(\Sigma^{\Bbb N})\\
	{}\\
	U\mu\in G
	\endSb
 \left(\int\varphi\,d\mu+h(\mu)\right)\,. 
 \endalign
 $$
This completes the proof of inequality (7.10). 
Inequality (7.11) is proved similarly. 
\hfill$\square$

\bigskip

We now turn towards the second main result in this section, namely, 
Theorem 7.5 showing that
 the multifractal pressure
 and the modified multifractal pressure are (almost)
comparable.
We first 
prove a small auxiliary lemma.

\bigskip

\proclaim{Lemma 7.4}
Let $(X,\distance)$ be a metric space
and let $U:\Cal P(\Sigma^{\Bbb N})\to X$ be continuous with respect to the weak topology.
Let $C$ be a subset of $X$
and $r>0$.
\roster
\item"(1)"
For all $n$, we have
 $$
 \Big\{
 \bold u\in\Sigma^{n}
 \,\Big|\,
 UL_{u}[\bold u]\subseteq C
\Big\}
\subseteq
\Big\{
 \bold u\in\Sigma^{n}
 \,\Big|\,
 UM_{n}[\bold u]\subseteq C
\Big\}\,.
$$

\item"(2)"
There is a positive integer $N_{r}$
such that
if $n\ge N_{r}$,
$\bold u\in\Sigma^{n}$ and $\bold k,\bold l\in\Sigma^{\Bbb N}$, then
we have
$\distance
 \big(
 \,
 UL_{n}(\bold u\bold k)
 \,,\,
 UL_{n}(\bold u\bold l)
 \,
 \big)
 \le
 r$.
 \item"(3)"
There is a positive integer $N_{r}$
such that
if $n\ge N_{r}$,
then
we have
 $$
\Big\{
 \bold u\in\Sigma^{n}
 \,\Big|\,
 UM_{n}[\bold u]\subseteq C
\Big\}
\subseteq
 \Big\{
 \bold u\in\Sigma^{n}
 \,\Big|\,
 \,\,
 UL_{u}[\bold u]\subseteq B(C,r)\Big\}\,.
$$
\endroster
\endproclaim
\noindent{\it Proof}\newline
\noindent
(1) This statement follows immediately from the fact that
 if $\bold u\in\Sigma^{n}$, then
$M_{n}[\bold u]
 =
 \{\,L_{n}\overline{\bold u}\,\}$.

\noindent
(2)
Fix $\gamma\in(0,1)$
and let $\distance_{\gamma}$
denote the metric on $\Sigma^{\Bbb N}$
introduced in Section 6.3.
For a function $f:\Sigma^{\Bbb N}\to\Bbb R$, we let
$\Lip(f)$ denote the Lipschitz constant of $f$
with respect to the metric
$\distance_{\gamma}$, 
i\.e\.
$\Lip(f)
=
\sup_{
\bold i,\bold j\in\Sigma^{\Bbb N},
\bold i\not=\bold j}
\frac
{|f(\bold i)-f(\bold j)|}
{\distance_{\gamma}(\bold i,\bold j)}$
and
define the metric $\LDistance$ in $\Cal P(\Sigma^{\Bbb N})$ by
$$
\LDistance(\mu,\nu)
=
\sup
\Sb
f:\Sigma^{\Bbb N}\to\Bbb R\\
\Lip(f)\le 1
\endSb
\Bigg|
\int f\,d\mu-\int f\,d\nu
\Bigg|;
$$
we note that it is well-known that $\LDistance$ is a metric
and  that $\LDistance$ induces the weak topology.
Since 
$U:\Cal P(\Sigma^{\Bbb N})\to X$ is continuous and 
$\Cal P(\Sigma^{\Bbb N})$ is compact, we
conclude that
$U:\Cal P(\Sigma^{\Bbb N})\to X$
is uniformly continuous.
This implies that we can choose $\delta>0$
such that
all measures
$\mu,\nu\in\Cal P(\Sigma^{\Bbb N})$
satisfy the following implication:
 $$
 \LDistance(\mu,\nu)<\delta
 \,\,\,\,
 \Rightarrow
 \,\,\,\,
 \distance(U\mu,U\nu)<r\,.
 \tag7.15
 $$

Next,
choose a positive integer $N_{r}$ such that
 $$
 \frac{1}{N_{r}(1-\gamma)}
 <
 \delta\,.
 \tag7.16
 $$

If $n\ge N_{r}$,
$\bold u\in\Sigma^{n}$ and $\bold k,\bold l\in\Sigma^{\Bbb N}$, then
it follows from 
(7.16) that
 $$
 \align
\LDistance
 \big(
 \,
 L_{n}(\bold u\bold k)
 \,,\,
 L_{n}(\bold u\bold l)
 \,
 \big)
&=
\sup
\Sb
f:\Sigma^{\Bbb N}\to\Bbb R\\
\Lip(f)\le 1
\endSb
\Bigg|
\int f\,d(L_{n}(\bold u\bold k))-\int f\,d(L_{n}(\bold u\bold l))
\Bigg|\\
&=
\sup
\Sb
f:\Sigma^{\Bbb N}\to\Bbb R\\
\Lip(f)\le 1
\endSb
\Bigg|
\frac{1}{n}\sum_{i=0}^{n-1} f(S^{i}(\bold u\bold k))
-
\frac{1}{n}\sum_{i=0}^{n-1} f(S^{i}(\bold u\bold l))
\Bigg|\\
&\le
\sup
\Sb
f:\Sigma^{\Bbb N}\to\Bbb R\\
\Lip(f)\le 1
\endSb
\frac{1}{n}\sum_{i=0}^{n-1} 
|
f(S^{i}(\bold u\bold k))
-
 f(S^{i}(\bold u\bold l))
|\\
&\le
\frac{1}{n}\sum_{i=0}^{n-1} 
\distance_{\Sigma^{\Bbb N}}
\big(
\,
S^{i}(\bold u\bold k)
\,,\,
S^{i}(\bold u\bold l)
\,
\big)\\
&=
\frac{1}{n}\sum_{i=0}^{n-1} 
\frac{1}
{N^{
|S^{i}(\bold u\bold k)
\wedge
S^{i}(\bold u\bold l)|
}
}\\
&\le
\frac{1}{N_{r}}\sum_{i=0}^{n-1} 
\frac{1}{N^{n-i}}\\
&\le
 \frac{1}{N_{r}(1-\gamma)}\\
&<
\delta\,,
\endalign
$$
and we therefore conclude from (7.15) that
$\distance
 (
 \,
 UL_{n}(\bold u\bold k)
 \,,\,
 UL_{n}(\bold u\bold l)
 \,
 )
 <
r$.

\noindent
(3)
It follows from Part (2) that
there is a positive integer $N_{r}$
such that
if $n\ge N_{r}$,
$\bold u\in\Sigma^{n}$ and $\bold k,\bold l\in\Sigma^{\Bbb N}$, then
 $\distance
 (
 \,
 UL_{n}(\bold u\bold k)
 \,,\,
 UL_{n}(\bold u\bold l)
 \,
 )
<
 r$.
We now claim that
if $n\ge N_{r}$,
then
 $$
\Big\{
 \bold u\in\Sigma^{n}
 \,\Big|\,
 UM_{n}[\bold u]\subseteq C
\Big\}
\subseteq
 \Big\{
 \bold u\in\Sigma^{n}
 \,\Big|\,
 \,\,
 UL_{u}[\bold u]\subseteq B(C,r)\Big\}\,.
$$
In order to prove this inclusion, we fix
$n\ge N_{r}$
and
$\bold u\in\Sigma^{n}$
with
$UM_{n}[\bold u]\subseteq C$.
We must now prove that
$UL_{n}[\bold u]\subseteq B(C,r)$.
Fix $\bold i\in[\bold u]$.
Since $\bold i\in[\bold u]$,
we can now find a (unique) $\bold k\in\Sigma^{\Bbb N}$ such that
$\bold i=\bold u\bold k$, whence
 $$
 \align
 \dist
 \big(
 \,
 UL_{n}\bold i
 \,,\,
 C
 \,
 \big)
&\le
\distance
\big(
\,
UL_{n}\bold i
\,,\,
UL_{n}\overline{\bold u}
\,
\big)
+
 \dist
 \big(
 \,
 UL_{n}\overline{\bold u}
 \,,\,
 C
 \,
 \big)\\
 &=
\distance
\big(
\,
UL_{n}(\bold u\bold k)
\,,\,
UL_{n}(\bold u\overline{\bold u})
\,
\big)
+
 \dist
 \big(
 \,
 UL_{n}\overline{\bold u}
 \,,\,
 C
 \,
 \big)\,.
 \tag7.17
\endalign
$$
However, since
$n\ge N_{r}$
and $\bold u\in\Sigma^{n}$,
we conclude that
 $\distance
 (
 \,
 UL_{n}(\bold u\bold k)
 \,,\,
 UL_{n}(\bold u\overline{\bold u})
 \,
 )
<
r$.
Also,
$UL_{n}\overline{\bold u}
=
UM_{n}\overline{\bold u}
\in
UM_{n}[\bold u]
\subseteq
 C$, whence
$\dist
 (
 \,
 UL_{n}\overline{\bold u}
 \,,\,
 C
 \,
 )
 =
0$.
It therefore follows from (7.17) that
$$
 \align
 \dist
 \big(
 \,
 UL_{n}\bold i
 \,,\,
 C
 \,
 \big)
 &=
\distance
\big(
\,
UL_{n}(\bold u\bold k)
\,,\,
UL_{n}(\bold u\overline{\bold u})
\,
\big)
+
 \dist
 \big(
 \,
 UL_{n}\overline{\bold u}
 \,,\,
 C
 \,
 \big)\\
&<
r\,. 
\endalign
$$
This completes the proof.
\hfill$\square$

\bigskip

We can now state and prove the second main result in this section, namely,
Theorem 7.5.

\bigskip

\proclaim{Theorem 7.5}
Let $X$ be a metric space and let $U:\Cal P(\Sigma^{\Bbb N})\to X$ be 
continuous with respect to the weak topology.
Let $C\subseteq X$ be a  subset of $X$ and $r>0$.
Fix a continuous function $\varphi:\Sigma^{\Bbb N}\to\Bbb R$.
Then we have
 $$
 \gather
 \underline P_{C}^{U}(\varphi)
 \le
 \underline Q_{C}^{U}(\varphi)
 \le
 \underline P_{B(C,r)}^{U}(\varphi)\,,\\
 \overline P_{C}^{U}(\varphi)
 \le
 \overline Q_{C}^{U}(\varphi)
 \le
 \overline P_{B(C,r)}^{U}(\varphi)\,.
 \endgather
 $$
\endproclaim
\noindent{\it Proof}\newline
\noindent
This follows immediately from Lemma 7.4.
\hfill$\square$

  \bigskip


\heading
{
8. Proof of Theorem 5.3}
\endheading

The purpose of this section is to prove Theorem 5.3.

\bigskip

\proclaim{Lemma 8.1}
Let $X$ be  a metric space and let $F:X\to\Bbb R$ be an upper 
semi-continuous function.
Let $K_{1},K_{2},\ldots\subseteq X$ be non-empty compact subsets of 
$X$ with
$K_{1}\supseteq K_{2}\supseteq\ldots$. Then
 $$
 \inf_{n}\sup_{x\in K_{n}}F(x)
 =
 \sup_{x\in\bigcap_{n}K_{n}}F(x)\,.
 $$
\endproclaim
\noindent{\it Proof}\newline
First note that it is clear that
$\inf_{n}\sup_{x\in K_{n}}F(x)
 \ge
 \sup_{x\in\cap_{n}K_{n}}F(x)$.
 We will now prove the reverse inequality, namely,
$\inf_{n}\sup_{x\in K_{n}}F(x)
 \le
 \sup_{x\in\cap_{n}K_{n}}F(x)$. 
Let $\varepsilon>0$.
For each $n$, we can choose $x_{n}\in K_{n}$ such that
$F(x_{n})\ge \sup_{x\in K_{n}}F(x)-\varepsilon$.
Next, since $K_{n}$ is compact for all $n$ and
$K_{1}\supseteq K_{2}\supseteq\ldots$, we can find a subsequence $(x_{n_{k}})_{k}$ and a 
point $x_{0}\in\cap_{n}K_{n}$ such that $x_{n_{k}}\to x_{0}$.
Also, since
$K_{n_{1}}\supseteq K_{n_{2}}\supseteq\ldots$, we conclude that
$\sup_{x\in K_{n_{1}}}F(x)
\ge
\sup_{x\in K_{n_{2}}}F(x)
\ge
\ldots$, 
whence
$\inf_{k}\sup_{x\in K_{n_{k}}}F(x)
=
\limsup_{k}\sup_{x\in K_{n_{k}}}F(x)$.
This implies that
$\inf_{n}\sup_{x\in K_{n}}F(x)
\le
\inf_{k}\sup_{x\in K_{n_{k}}}F(x)
=
\limsup_{k}\sup_{x\in K_{n_{k}}}F(x)
\le
\limsup_{k}F(x_{n_{k}})+\varepsilon$.
However, since 
$x_{n_{k}}\to x_{0}$, we deduce from the 
upper semi-continuity of the function $F$, that
$\limsup_{k}F(x_{n_{k}})\le F(x_{0})$.
Consequently
$\inf_{n}\sup_{x\in K_{n}}F(x)
\le
\limsup_{k}F(x_{n_{k}})+\varepsilon
\le
F(x_{0})+\varepsilon
\le
 \sup_{x\in\cap_{n}K_{n}}F(x)+\varepsilon$.
 Finally, letting $\varepsilon\searrow 0$ gives the desired result.
\hfill$\square$

\bigskip

\noindent
We can now prove Theorem 5.3.

\bigskip

\noindent{\it Proof of Theorem 5.3}\newline
\noindent (1)
We must prove the following two inequalities, namely,
 $$
 \gather
  \sup
\Sb
\mu\in\Cal P_{S}(\Sigma^{\Bbb N})\\
{}\\
U\mu\in\overline  C
\endSb
\Bigg(
h(\mu)+\int\varphi\,d\mu
\Bigg)
\le
\inf_{r>0}
\,\,
\underline P_{B(C,r)}^{U}(\varphi)\,,
\tag8.1\\
 \inf_{r>0}
\,\,
\overline P_{B(C,r)}^{U}(\varphi)
\le
 \sup
\Sb
\mu\in\Cal P_{S}(\Sigma^{\Bbb N})\\
{}\\
U\mu\in\overline  C
\endSb
\Bigg(
h(\mu)+\int\varphi\,d\mu
\Bigg)\,.
\tag8.2
 \endgather
 $$

{\it Proof of (8.1).}
Since $B(C,r)$ is open with 
$\overline C\subseteq B(C,r)$, 
we conclude from Theorem 7.3 that
 $$
 \align
 \sup
\Sb
\mu\in\Cal P_{S}(\Sigma^{\Bbb N})\\
{}\\
U\mu\in\overline  C
\endSb
\Bigg(
h(\mu)+\int\varphi\,d\mu
\Bigg)
&\le
\sup
\Sb
\mu\in\Cal P_{S}(\Sigma^{\Bbb N})\\
{}\\
U\mu\in B(C,r)
\endSb
\Bigg(
h(\mu)+\int\varphi\,d\mu
\Bigg)\\
&{}\\
&\le
\underline Q_{B(C,r)}^{U}(\varphi)\,.
\tag8.3
\endalign
$$
Taking infimum over all $r>0$ in (8.3) gives
  $$
 \align
 \sup
\Sb
\mu\in\Cal P_{S}(\Sigma^{\Bbb N})\\
{}\\
U\mu\in \overline C
\endSb
\Bigg(
h(\mu)+\int\varphi\,d\mu
\Bigg)
&\le
\inf_{r>0}
\,\,
\underline Q_{B(C,r)}^{U}(\varphi)\,.
\tag8.4
\endalign
$$
Next, we note that it follows from Theorem 7.5
that
$\underline Q_{B(C,r)}^{U}(\varphi)
 \le
 \underline P_{B(\,B(C,r)\,,\,r\,)}^{U}(\varphi)$.
Combining this inequality with and (8.4) and using the fact that 
$B(\,B(C,r)\,,\,r\,)\subseteq B(C,2r)$, we now conclude that
 $$
 \align
 \sup
\Sb
\mu\in\Cal P_{S}(\Sigma^{\Bbb N})\\
{}\\
U\mu\in \overline C
\endSb
\Bigg(
h(\mu)+\int\varphi\,d\mu
\Bigg)
&\le
\inf_{r>0}
\,\,
\underline Q_{B(C,r)}^{U}(\varphi)\\
&\le
\inf_{r>0}
\,\,
\underline P_{B(\,B(C,r)\,,\,r\,)}^{U}(\varphi)\\
&\le
\inf_{r>0}
\,\,
\underline P_{B(C,2r)}^{U}(\varphi)\\
&\le
\inf_{s>0}
\,\,
\underline P_{B(C,s)}^{U}(\varphi)\,.
\endalign
$$
This completes the proof of inequality (8.1).
 
\bigskip

{\it Proof of (8.2).}
Since
$\overline{B(C,r)}$ is closed,
we conclude from Theorem 7.3
that
  $$
  \align
  \inf_{r>0}
\,\,
\overline P_{B(C,r)}^{U}(\varphi)
&\le
  \inf_{r>0}
\,\,
\overline P_{\overline{B(C,r)}}^{U}(\varphi)\\
&\le
  \inf_{r>0}
\,\,
 \sup
\Sb
\mu\in\Cal P_{S}(\Sigma^{\Bbb N})\\
{}\\
U\mu\in \overline{B(C,r)}
\endSb
\Bigg(
h(\mu)+\int\varphi\,d\mu
\Bigg)\,.
\endalign
$$
Letting $U_{S}:\Cal P_{S}(\Sigma^{\Bbb N})\to X$ denote the restriction of $U$ to 
$\Cal P_{S}(\Sigma^{\Bbb N})$,
the above inequality can be written as
    $$
  \align
  \inf_{r>0}
\,\,
\overline P_{B(C,r)}^{U}(\varphi)
&\le
  \inf_{r>0}
\,\,
 \sup
\Sb
\mu\in U_{S}^{-1}\overline{B(C,r)}
\endSb
\Bigg(
h(\mu)+\int\varphi\,d\mu
\Bigg)\\
&=
  \inf_{n}
\,\,
 \sup
\Sb
\mu\in U_{S}^{-1}\overline{B(C,\frac{1}{n})}
\endSb
\Bigg(
h(\mu)+\int\varphi\,d\mu
\Bigg)\,.
\tag8.5
\endalign
$$

Next, note that
since
$  \overline{B(C,\frac{1}{n})}$ is closed and $U_{S}$ is continuous,
the set
$U_{S}^{-1}\overline{B(C,\frac{1}{n})}$
is a closed 
  subset of 
 $\Cal P_{S}(\Sigma^{\Bbb N})$.
 As  $\Cal P_{S}(\Sigma^{\Bbb N})$ is compact, we therefore
 deduce that
  $U_{S}^{-1}\overline{B(C,\frac{1}{n})}$
  is compact.
 Also, note that it follows from [Wa]
 that the
 entropy map $h:\Cal P_{S}(\Sigma^{\Bbb N})\to\Bbb R$
 is upper semi-continuous.
 We conclude from this
 that the
 map $F:\Cal P_{S}(\Sigma^{\Bbb N})\to\Bbb R$
 defined by
 $F(\mu)
 =
  h(\mu)+\int\varphi\,d\mu$ is upper semi-continuous.
  Finally, since the sets 
  $K_{n}=U_{S}^{-1}\overline{B(C,\frac{1}{n})}$
  are compact with 
  $K_{1}\supseteq K_{2}\supseteq K_{3}\supseteq\ldots$
  and
  $F$ is upper semi-continuous, we deduce from Lemma 8.1 that
   $$
  \align
  \inf_{n}
  \,\,
  \sup
\Sb
\mu\in U_{S}^{-1}\overline{B(C,\frac{1}{n})}
\endSb
\Bigg(
h(\mu)+\int\varphi\,d\mu
\Bigg)
&=
\inf_{n}
\sup
\Sb
\mu\in K_{n}
\endSb
F(\mu)\\
&=
\sup
\Sb
\mu\in \bigcap_{n}K_{n}
\endSb
F(\mu)\\
&=
\sup
\Sb
\mu\in \bigcap_{n}U_{S}^{-1}\overline{B(C,\frac{1}{n})}
\endSb
\Bigg(
h(\mu)+\int\varphi\,d\mu
\Bigg)\,.
\tag8.6
\endalign
$$

 Observe that
  $\bigcap_{n}U_{S}^{-1}\overline{B(C,\frac{1}{n})}
  \subseteq
  U_{S}^{-1}(\,\bigcap_{n}\overline{B(C,\frac{1}{n})})
  =
  U_{S}^{-1}\overline C$,
  whence
 $$
 \align
\sup
\Sb
\mu\in \bigcap_{n}U_{S}^{-1}\overline{B(C,\frac{1}{n})}
\endSb
\Bigg(
h(\mu)+\int\varphi\,d\mu
\Bigg)
&\le
\sup
\Sb
\mu\in \bigcap_{n}U_{S}^{-1}\overline C
\endSb
\Bigg(
h(\mu)+\int\varphi\,d\mu
\Bigg)\\
&=
\sup
\Sb
\mu\in \Cal P_{S}(\Sigma^{\Bbb N})\\
{}\\
U\mu\in\overline C
\endSb
\Bigg(
h(\mu)+\int\varphi\,d\mu
\Bigg)\,.
\tag8.7
\endalign
 $$

  Finally, combining (8.5), (8.6) and (8.7) gives inequality (8.2).

  \noindent
  (2)
  This part follows immediately from Part (1) and Proposition 5.1.
  \hfill$\square$

  \bigskip


\heading
{9. Proof of Theorem 5.5}
\endheading

The purpose of this section is to prove Theorem 5.5.
We first prove two small lemmas.

\bigskip

\proclaim{Lemma 9.1}
let
$\Delta:\Cal P(\Sigma^{\Bbb N})\to \Bbb R$
be continuous 
with
$\Delta(\mu)\not=0$
for all $\mu\in\Cal P(\Sigma^{\Bbb N})$.
The either  $\Delta<0$ or $\Delta>0$.
\endproclaim
\noindent{\it Proof}\newline
\noindent
Assume, in order to reach a contradiction,
that there are $\mu_{-},\mu_{+}\in\Cal P(\Sigma^{\Bbb N})$
such that
$\Delta(\mu_{-})<0$ and $\Delta(\mu_{+})>0$.
For $t\in[0,1]$,
 let
$\mu_{t}=t\mu_{-}+(1-t)\mu_{+}\in\Cal P(\Sigma^{\Bbb N})$
and define
$f:[0,1]\to\Bbb R$ by
$f(t)=\Delta(\mu_{t})$.
The function $f$ is clearly continuous 
with
$f(0)=\Delta(\mu_{+})>0$ and
$f(1)=\Delta(\mu_{-})<0$,
and we therefore conclude from the 
intermediate  value theorem that there is 
a number
$t_{0}\in(0,1)$ such that
$\Delta(\mu_{t_{0}})=f(t_{0})=0$.
However, this clearly contradicts the fact that
$\Delta(\mu)\not=0$
for all $\mu\in\Cal P(\Sigma^{\Bbb N})$.
\hfill$\square$

\bigskip

\proclaim{Lemma 9.2}
Let $X$ be a normed vector space.
Let $\Gamma:\Cal P(\Sigma^{\Bbb N})\to X$
be continuous and affine
and let
$\Delta:\Cal P(\Sigma^{\Bbb N})\to \Bbb R$
be continuous and affine
with
$\Delta(\mu)\not=0$
for all $\mu\in\Cal P(\Sigma^{\Bbb N})$.
Define 
$U:\Cal P(\Sigma^{\Bbb N})\to X$
by
$U=\frac{\Gamma}{\Delta}$.
Let $C$ be a closed and convex subset of $X$ and assume that
 $$
 \overset{\,\circ}\to{C}
 \cap
 \,
 U\big(\,\Cal P_{S}(\Sigma^{\Bbb N})\,\big)
 \not=
 \varnothing\,.
 $$
Then
$$
\sup
\Sb
\mu\in \Cal P_{S}(\Sigma^{\Bbb N})\\
{}\\
U\mu\in C
\endSb
\Bigg(
h(\mu)+\int\varphi\,d\mu
\Bigg)
=
\sup
\Sb
\mu\in \Cal P_{S}(\Sigma^{\Bbb N})\\
U\mu\in \overset{\,\circ}\to{C}
\endSb
\Bigg(
h(\mu)+\int\varphi\,d\mu
\Bigg)\,.
$$
\endproclaim
\noindent{\it Proof}\newline
\noindent
For brevity define 
$F:\Cal P_{S}(\Sigma^{\Bbb N})\to\Bbb R$ by
$F(\mu)
=
h(\mu)+\int\varphi\,d\mu$.
It clearly suffices to show that
$$
\sup
\Sb
\mu\in \Cal P_{S}(\Sigma^{\Bbb N})\\
{}\\
U\mu\in C
\endSb
F(\mu)
\le
\sup
\Sb
\mu\in \Cal P_{S}(\Sigma^{\Bbb N})\\
U\mu\in \overset{\,\circ}\to{C}
\endSb
F(\mu)\,.
\tag9.1
$$

We will now prove inequality (9.1).
Write
$s
 =
\sup_{
\mu\in \Cal P_{S}(\Sigma^{\Bbb N})\,,\,
U\mu\in C
}
F(\mu)$.
Fix $\varepsilon>0$.
It follows from the definition of $s$ that we can choose 
$\lambda\in\Cal P_{S}(\Sigma^{\Bbb N})$ with
$U\lambda\in C$
and
$F(\lambda)>s-\varepsilon$.
Also, since
$\overset{\,\circ}\to{C}
 \cap
 \,
 U\big(\,\Cal P_{S}(\Sigma^{\Bbb N})\,\big)
 \not=
 \varnothing$, we can find
 $\nu\in\Cal P_{S}(\Sigma^{\Bbb N})$, with
 $U\nu\in \overset{\,\circ}\to{C}$.
For $t\in(0,1)$ we now define $\gamma_{t}\in\Cal P_{S}(\Sigma^{\Bbb N})$
by
$\gamma_{t}=t\nu+(1-t)\lambda$.
Next, we prove the following three claims.

\medskip

{\it Claim 1. 
For all $t\in(0,1)$, we have
$U\gamma_{t}\in \overset{\,\circ}\to{C}$.}

\noindent
{\it Proof of Claim 1.}
Fix $t\in(0,1)$.
Write
$a
=
\frac
{t\Delta(\nu)}
{t\Delta(\nu)+(1-t)\Delta(\lambda)}$
and 
$b
=
\frac
{(1-t)\Delta(\lambda)}
{t\Delta(\nu)+(1-t)\Delta(\lambda)}$.
We now make a few observations.
We first observe that
 it follows from Lemma 9.1
 that
either  $\Delta<0$ or $\Delta>0$.
This clearly implies that
 $a,b\in(0,1)$.
Next, we  note that
 $
  U\gamma_{t}
=
 \frac{\Gamma(t\nu+(1-t)\lambda)}{\Delta(t\nu+(1-t)\lambda)}
 =
 \frac
 {t\Gamma(\nu)+(1-t)\Gamma(\lambda)}
 {t\Delta(\nu)+(1-t)\Delta(\lambda)}
  =
  \frac
 {t\Gamma(\nu)}
 {t\Delta(\nu)+(1-t)\Delta(\lambda)} 
 +
  \frac
 {(1-t)\Gamma(\lambda)}
 {t\Delta(\nu)+(1-t)\Delta(\lambda)} 
 =
  \frac
 {t\Delta(\nu)}
 {t\Delta(\nu)+(1-t)\Delta(\lambda)} U\nu
 +
  \frac
 {(1-t)\Delta(\lambda)}
 {t\Delta(\nu)+(1-t)\Delta(\lambda)} U\lambda
 =
 aU\nu+bU\lambda$.
 We can now prove that
 $U\gamma_{t}\in \overset{\,\circ}\to{C}$.
 Indeed,
 since
 $a,b\in(0,1)$ with $a+b=1$
 and
 $U\lambda\in C$
and
$U\nu\in \overset{\,\circ}\to{C}$,
we conclude from [Co, p\. 102, Proposition 1.11]
 that
  $
  U\gamma_{t}
 =
 aU\nu+bU\lambda\in\overset{\,\circ}\to{C}$.
 This completes the proof of Claim 1.

\medskip

{\it Claim 2. There is $t_{0}\in(0,1)$ such that
$F(\gamma_{t_{0}})>s-\varepsilon$.}

\noindent
{\it Proof of Claim 2.}
Since
the entropy function $h:\Cal P_{S}(\Sigma^{\Bbb N})$ is affine (see [Wa]),
we conclude that
 $F$ is affine, and so
 $F(\gamma_{t})
 =
 F(t\nu+(1-t)\lambda)
=
tF(\nu)+(1-t)F(\lambda)
\to
F(\lambda)
>
s-\varepsilon$.
This implies that there is $t_{0}\in(0,1)$
with
$F(\gamma_{t_{0}})>s-\varepsilon$.
This completes the proof of Claim 2.

\medskip

{\it Claim 3. There is $\pi\in \Cal P_{S}(\Sigma^{\Bbb N})$ with 
$U\pi\in \overset{\,\circ}\to{C}$
such that
$F(\pi)>s-\varepsilon$.}

\noindent
{\it Proof of Claim 3.}
It follows
from
Claim 2 that there is $t_{0}\in(0,1)$ such that
$F(\gamma_{t_{0}})>s-\varepsilon$
and Claim 1 implies that
$U\gamma_{t_{0}}\in\overset{\,\circ}\to{C}$.
We now put $\pi=\gamma_{t_{0}}$.
This completes the proof of Claim 3.

\medskip

We can now prove inequality (9.1). It follows from Claim 3 that there is 
$\pi\in \Cal P_{S}(\Sigma^{\Bbb N})$ with 
$U\pi\in \overset{\,\circ}\to{C}$
such that
$F(\pi)>s-\varepsilon$,
whence
$s-\varepsilon
<
F(\pi)
\le
\sup_{
\mu\in \Cal P_{S}(\Sigma^{\Bbb N})\,,\,U\mu\in \overset{\,\circ}\to{C}
}
F(\mu)$.
Finally, letting $\varepsilon\searrow 0$ gives
$s
\le
\sup_{
\mu\in \Cal P_{S}(\Sigma^{\Bbb N})
\,,\,
U\mu\in \overset{\,\circ}\to{C}
}
F(\mu)$.
\hfill$\square$

\bigskip

\noindent
We can now prove Theorem 5.5.

\bigskip

\noindent{\it Proof of Theorem 5.5}\newline
\noindent
In view of Lemma 9.2, 
it suffices
to prove the following two inequalities, namely,
 $$
\gather
\sup
\Sb
\mu\in \Cal P_{S}(\Sigma^{\Bbb N})\\
{}\\
U\mu\in \overset{\,\circ}\to{C}
\endSb
\Bigg(
h(\mu)+\int\varphi\,d\mu
\Bigg)
\le
\underline P_{C}^{U}(\varphi)\,,
\tag9.2\\
{}\\
\overline P_{C}^{U}(\varphi)
\le
\sup
\Sb
\mu\in \Cal P_{S}(\Sigma^{\Bbb N})\\
{}\\
U\mu\in C
\endSb
\Bigg(
h(\mu)+\int\varphi\,d\mu
\Bigg)\,.
\tag9.3
\endgather
$$

{\it Proof of inequality (9.2).}
For $r>0$, let
$G_{r}
=
\{x\in C\,|\,\dist(x,X\setminus C)>r\}$,
and
note that
$G_{r}$ is open 
with
$B(G_{r},\rho)\subseteq C$
for all $0<\rho<r$.
We therefore conclude from Theorem 7.3
and Theorem 7.5
 that if $0<\rho<r$, then
 $$
 \align
 \underline P_{C}^{U}(\varphi)
&\ge
  \underline P_{ B(G_{r},\rho)}^{U}(\varphi)\\
 &\ge
  \underline Q_{ G_{r}}^{U}(\varphi)
  \qquad\qquad
  \qquad\qquad
  \qquad\,\,\,\,\,
  \text{[by Theorem 7.5]}\\
 &\ge
\sup
\Sb
\mu\in \Cal P_{S}(\Sigma^{\Bbb N})\\
{}\\
U\mu\in G_{r}
\endSb
\Bigg(
h(\mu)+\int\varphi\,d\mu
\Bigg)\,.
  \qquad
  \text{[by Theorem 7.3]}
  \tag9.4
 \endalign
 $$
Taking supremum over all $r>0$ in (9.4)  yields
 $$
 \align
 \underline P_{C}^{U}(\varphi)
 &\ge
 \sup_{r>0}
 \,\,
\sup
\Sb
\mu\in \Cal P_{S}(\Sigma^{\Bbb N})\\
{}\\
U\mu\in G_{r}
\endSb
\Bigg(
h(\mu)+\int\varphi\,d\mu
\Bigg)\,.
 \endalign
 $$
Letting $U_{S}:\Cal P_{S}(\Sigma^{\Bbb N})\to(0,\infty)$ denote the restriction of $U$ to
$\Cal P_{S}(\Sigma^{\Bbb N})$,
the previous inequality can be written as
 $$
 \align
 \underline P_{C}^{U}(\varphi)
 &\ge
 \sup_{r>0}
 \,\,
\sup
\Sb
\mu\in U_{S}^{-1} G_{r}
\endSb
\Bigg(
h(\mu)+\int\varphi\,d\mu
\Bigg)\\
 &=
\sup
\Sb
\mu\in \bigcup_{r>0}U_{S}^{-1} G_{r}
\endSb
\Bigg(
h(\mu)+\int\varphi\,d\mu
\Bigg)\,.
\tag9.5
 \endalign
 $$
However, it is easily seen that
$\bigcup_{r>0} G_{r}
\supseteq
\overset{\,\circ}\to{C}$,
whence
$\bigcup_{r>0}U_{S}^{-1} G_{r}
=
U_{S}^{-1}(\,\bigcup_{r>0} G_{r})
\supseteq
U_{S}^{-1}\overset{\,\circ}\to{C}$.
We conclude from this inclusion and  inequality (9.5) that
  $$
 \align
 \underline P_{C}^{U}(\varphi)
 &\ge
\sup
\Sb
\mu\in U_{S}^{-1}\overset{\,\circ}\to{C}
\endSb
\Bigg(
h(\mu)+\int\varphi\,d\mu
\Bigg)\\
&=
\sup
\Sb
\mu\in \Cal P_{S}(\Sigma^{\Bbb N})\\
U\mu\in \overset{\,\circ}\to{C}
\endSb
\Bigg(
h(\mu)+\int\varphi\,d\mu
\Bigg)\,.
 \endalign
 $$
This proves inequality (9.2).

\bigskip

{\it Proof of inequality (9.3).}
Since $C$ is closed we immediately 
 conclude from Theorem 7.3
and Theorem 7.5
 that
 $$
 \align
\overline P_{C}^{U}(\varphi)
&\le
\overline Q_{C}^{U}(\varphi)
  \qquad\qquad
  \qquad\qquad
  \qquad\,\,\,\,\,\,\,\,
  \text{[by Theorem 7.5]}\\
&\le
\sup
\Sb
\mu\in \Cal P_{S}(\Sigma^{\Bbb N})\\
{}\\
U\mu\in C
\endSb
\Bigg(
h(\mu)+\int\varphi\,d\mu
\Bigg)\,.
  \qquad
  \text{[by Theorem 7.3]}
\endalign
$$
This proves inequality (9.3).
\hfill$\square$

\end